\documentclass[12pt,a4paper,reqno]{amsart}

\usepackage{a4wide}
\usepackage{amssymb}
\usepackage{graphicx}

\theoremstyle{plain}
\newtheorem{theorem}{Theorem}[section]
\newtheorem{maintheorem}{Theorem}
\newtheorem{lemma}[theorem]{Lemma}
\newtheorem{proposition}[theorem]{Proposition}
\newtheorem{corollary}[theorem]{Corollary}

\newtheorem*{claim}{Claim}

\theoremstyle{definition}
\newtheorem{definition}[theorem]{Definition}

\theoremstyle{remark}
\newtheorem{remark}[theorem]{Remark}

\def\R{\ensuremath{\mathbb R}}
\def\N{\ensuremath{\mathbb N}}

\def\Z{\ensuremath{\mathbb Z}}
\def\I{\ensuremath{{\bf 1}}}
\def\e{\ensuremath{\text e}}
\def\U{\ensuremath{\mathcal U}}
\def\B{\ensuremath{\mathcal BC}}
\def\C{\ensuremath{\mathcal C}}
\def\L{\ensuremath{\mathcal L}}
\def\P{\ensuremath{\mathcal P}}
\def\O{\ensuremath{\mathcal O}}

\newcommand{\qand}{\quad\text{and}\quad}

\numberwithin{equation}{section}

\begin{document}


\title[Statistical Stability for H\'{e}non maps]{Statistical Stability for H\'{e}non maps of the Benedicks-Carleson type}
\author[J. F. Alves]{Jos\'{e} F. Alves}
\address{Jos\'{e} F. Alves\\ Centro de Matem\'{a}tica da Universidade do Porto\\ Rua do Campo Alegre 687\\ 4169-007 Porto\\ Portugal}
\email{jfalves@fc.up.pt} \urladdr{http://www.fc.up.pt/cmup/jfalves}

\author[M. Carvalho]{Maria Carvalho}
\address{Maria Carvalho\\ Centro de Matem\'{a}tica da Universidade do Porto\\ Rua do
Campo Alegre 687\\ 4169-007 Porto\\ Portugal}
\email{mpcarval@fc.up.pt}

\author[J. M. Freitas]{Jorge Milhazes Freitas}
\address{Jorge M. Freitas\\ Centro de Matem\'{a}tica da Universidade do Porto\\ Rua do
Campo Alegre 687\\ 4169-007 Porto\\ Portugal}
\email{jmfreita@fc.up.pt}
\urladdr{http://www.fc.up.pt/pessoas/jmfreita}

\date{\today}
\thanks{Work partially supported by
FCT through CMUP and POCI/MAT/61237/2004. } 

\begin{abstract}
We consider the family of H\'enon maps in the plane and  show that
the SRB measures vary continuously in the weak* topology within
the set of Benedicks-Carleson parameters.

\end{abstract}

\maketitle

\setcounter{tocdepth}{1}

\tableofcontents 

\section{Introduction}

At the end of the 19th century, Poincar\'{e} addressed the problem
of evolution and stability of the solar system, which led to many
surprising questions and  gave birth to the modern theory of
dynamical systems as a qualitative study of the dynamics. The main
goal of this theory is the description of the \emph{typical
behavior} of orbits, and the understanding of how this behavior
changes when we perturb the system or to which extent it is
robust. In the present work we are specially concerned with the
problem of the \textit{stability} of systems.

The first fundamental concept of robustness, \textit{structural
stability}\index{structural stability}, was formulated in the late
30's by Andronov and Pontryagin. It requires the persistence of
the orbit topological structure under small perturbations,
expressed in terms of the existence of a homeomorphism sending
orbits of the initial system onto orbits of the perturbed one.
This concept is tied with the notion of \textit{uniform
hyperbolicity} introduced by Smale in the mid 60's. A stronger
connection was conjectured by Palis and Smale in 1970: a
diffeomorphism is structurally stable if and only if it is
uniformly hyperbolic and satisfies the so-called transversality
condition. During this decade the ``if'' part of the conjecture
was solved due to the contributions of Robbin, de Melo and
Robinson; in the 80's Ma\~{n}\'{e} settled the $C^1$-stability
conjecture. The flow case was solved by Aoki and Hayashi,
independently, in the 90's, also in the $C^1$-topology.

In spite of these astonishing successes, structural stability
proved to be somewhat restrictive. Several important models, such
as Lorenz flows, H\'{e}non maps and other non-uniformly hyperbolic
systems fail to present structural stability, although some key
aspects of a statistical nature persist after small perturbations.
The contributions of Kolmogorov, Sinai, Ruelle, Bowen, Oseledets,
Pesin, Katok, Ma\~n\'e and many others turned the attention of the
study of dynamical systems from a topological perspective to a
more statistical approach, and Ergodic Theory experienced an
unprecedent development. In trying to capture this statistical
persistence of phenomena, Alves and Viana \cite{AV02} proposed the
notion of \textit{statistical stability}\index{statistical
stability}, which expresses the continuous variation of physical
measures as a function of the evolution law governing the systems.
A \textit{physical measure} for a smooth map $f:M\rightarrow M$ on
a manifold $M$ is a Borel probability measure $\mu$ on $M$ for
which there is a positive Lebesgue measure set of points $x\in M$,
whose union forms the \emph{basin} of~$\mu$, such that
\begin{equation}
\label{eq:def-physical-measure}
\lim_{n\rightarrow\infty}\frac{1}{n}\sum_{j=0}^{n-1}\varphi
\left(f^j(x)\right)=\int\varphi\, d\mu,
\end{equation}
for any continuous function $\varphi:M\rightarrow\R$. 
Physical measures are intimately connected with
\textit{Sinai-Ruelle-Bowen measures}\index{SRB measures} (SRB for
short). An $f$-invariant Borel probability measure $\mu$ is said
to be SRB if it has a positive Lyapunov exponent and the
conditional measures of $\mu$ on unstable leaves are absolutely
continuous with respect to the Riemannian measure induced on those
leaves; see Section~\ref{subsec:SRB-measures} for a precise
definition. The existence of SRB measures for general dynamical
systems is usually a difficult problem. However,  Sinai, Ruelle
and Bowen established the existence of SRB measures  for Axiom A
attractors which qualify as physical measures. Moreover, Axiom A
diffeomorphisms are statistically stable.

The existence of SRB measures for a large set of one-dimensional
quadratic maps exhibiting non-uniformly expanding behavior has
been established in the pioneer paper of Jakobson
\cite{Ja81}. 
Additionally, the work of Collet-Eckmann \cite{CE80a,CE80b,CE83}
and Benedicks-Carleson \cite{BC85} became a major breakthrough in
that direction and allowed a well succeeded approach to higher
dimensional maps. A key ingredient is the exponential growth of
the derivative along the critical orbit for a positive Lebesgue
measure set $\B_1$ of parameters.
Regarding statistical stability, Freitas \cite{Fr05} showed that
the SRB measures vary continuously within the parameter set
$\B_1$; see also \cite{Ts96} and \cite{RS97} where Tsujii, Rychlik
and Sorets obtained related results. Notice that, by the work of
Thunberg \cite{Th01}, one cannot expect statistical stability on a
full Lebesgue measure set of parameters for the quadratic family.

H\'enon \cite{He76} proposed the two-parameter family of maps
\index{H\'enon maps}
\[ f_{a,b}:
\begin{array}[t]{ccc}
\R^2 & \longrightarrow & \R^2\\
(x,y)& \longmapsto & (1-ax^2+y,bx)
\end{array}
\]
as a model for non-linear two-dimensional dynamical systems, which
can be thought as a simplified discrete version of the Lorenz
flow. Based on numerical experiments for the parameters $a = 1.4$
and $b = 0.3$,  H\'enon  conjectured that this dynamical system
should have a {\em strange attractor}. In principle, most initial
points could be attracted to a periodic cycle, so it was not at
all \emph{a priori} clear that the attractor detected by H\'enon
in his experiments was not a long stable periodic orbit. However,
Benedicks and Carleson \cite{BC91} managed to prove that H\'enon's
conjecture was true for small $b>0$, 
showing that for a positive Lebesgue measure set
$\B$ 
the  map $f_{a,b}$ with $(a,b)\in\B$ exhibits a {\em
non-hyperbolic attractor}. Afterwards, Benedicks and Young
\cite{BY93} proved that each of these non-hyperbolic attractors
supports a unique SRB measure $\nu_{a,b}$, which is also a
physical measure and whose main statistical features were studied
by Benedicks, Viana and Young in
\cite{BV01, BV06,BY00}. 
Thus, a natural question is: are the H\'{e}non maps of the
Benedicks-Carleson type statistically stable?  The main result of
this work gives a positive answer to this question.

\begin{maintheorem}
The map which associates to each $(a,b)\in \B  $ the SRB measure $
\nu_{a,b} $ is continuous with respect to the weak* topology in
the space of probability measures.
\end{maintheorem}\index{Theorem!F}

Despite being metrically robust, the strange attractors appearing
for the Benedicks-Carleson parameters are very fragile. In fact,
Ures \cite{Ur95} showed that the Benedicks-Carleson parameters can
be approximated by other parameters for which the H\'enon map has
a homoclinic tangency associated to a fixed point. Hence,
according to Newhouse's results \cite{Ne74,Ne79}, one may deduce
the appearance of infinitely many attractors in the neighborhood
of the H\'enon attractor. Moreover, Ures \cite{Ur96} proved that
the SRB measures $\nu_{a,b}$ corresponding to $(a,b)\in\B$ can be
approximated by Dirac measures supported on sinks.  Nevertheless,
Benedicks and Viana \cite{BV06} showed that the H\'enon maps in
$\B$ are \textit{stochastically stable}. This notion  was
introduced by Kolmogorov and Sinai in the 70's and in broad terms
asserts that time-averages of continuous functions are only
slightly affected when iteration by the dynamics is perturbed by a
small random noise. Stochastic stability may imply statistical
stability if we allow a deterministic noise. However, the proof of
the stochastic stability in \cite{BV06} uses strongly the absolute
continuity of the random noise, which prevents us to deduce the
statistical stability from the stochastic one.

\section{Insight into the reasoning}

We consider a sequence of parameters $(a_n,b_n)_{n\in\N}\in\B$
converging to $(a_0,b_0)\in\B$. Let $(\nu_n)_{n\in\N}$ and $\nu_0$
denote the respective SRB measures. Our goal is to show that
$\nu_n$ converges to $\nu_0$ in the weak* topology. We prove this
by showing that every subsequence $(\nu_{n_i})_{i\in\N}$ contains
a subsequence convergent to $\nu_0$. Let us give some details on
how to find this convergent subsequence.

The main problem we have to overcome is the need of comparing
measures supported on different attractors. Our strategy is to
look for a common ground where the construction of the SRB measure
for every parameter is rooted. To do so, we start by noting that
each of these maps admits a horseshoe $\Lambda_{a,b}$ with
infinitely many branches and variable return times (we will drop
the indices when we refer to properties that apply to all these
objects) obtained by intersecting two transversal families of
local stable and unstable curves.
Besides, $\Lambda$ intersects each local unstable curve in a
positive Lebesgue measure Cantor set, and for each $z\in\Lambda$
it is possible to assign a positive integer $R(z)$ defining the
return time function $R:\Lambda\rightarrow\N$ which indicates that
$z$ returns to $\Lambda$ after $R(z)$ iterates. The hyperbolic
properties of $\Lambda$ and the good behavior of $R$ allow us to
build a Markov extension that organizes the dynamics of these
H\'{e}non maps. Thus, one needs to show first that for nearby
parameters the corresponding horseshoes are also close. We remark
that for each parameter there is not a unique horseshoe with the
required properties. Therefore, what we can establish is that for
a given parameter $(a,b)$ and a chosen horseshoe $\Lambda_{a,b}$,
if we consider a small perturbation $(a',b')$, then it is possible
to build a horseshoe $\Lambda_{a',b'}$ with the desired hyperbolic
properties and which is close to $\Lambda_{a,b}$.

These horseshoes play an important role in a  construction of the
SRB measures that suits our purposes. Actually,
$f^R:\Lambda\rightarrow\Lambda$ preserves a measure $\tilde{\nu}$
with absolutely continuous conditional measures on local unstable
curves with respect to the Lebesgue measure on each curve; the
good behavior of the function $R$ ensures that the saturation of
$\tilde{\nu}$ is an SRB measure, and by uniqueness it follows that
the saturation of $\tilde{\nu}$ is the SRB measure. To prove the
continuous dependence of these SRB measures on the parameter,
$\Lambda$ is collapsed along stable curves yielding a quotient
space $\bar{\Lambda}$, which can be thought inside a fixed local
unstable curve $\hat{\gamma}^u$, and whose elements are
represented by the intersection of the corresponding stable curve
with $\hat{\gamma}^u$. This way our task is reduced to analyze
$\overline{f^R}:\bar{\Lambda}\rightarrow \bar{\Lambda}$. This map
is piecewise uniformly expanding and its Perron-Frobenius operator
has a spectral gap under the usual
aperiodicity conditions; so there is 
an $\overline{f^R}$-invariant density with respect to Lebesgue
measure on~$\hat{\gamma}^u$. As $\hat\gamma^u$ is nearly
horizontal, we can think of $\bar{\rho}$ as a function defined on
a subset of the $x$-axis. 
The advantage of this perspective is that it gives us the desired
common domain for these densities, providing the first step in the
verification of the continuity.

Therefore, the steps for the construction of the convergent
subsequence are the following:
\begin{itemize}
\item Fix a parameter $(a_0,b_0)\in\B$ and a respective horseshoe
$\Lambda_0$.

\item Pick any sequence of parameters $(a_n,b_n)\in\B$
such that $(a_n,b_n)\rightarrow(a_0,b_0)$ as $n\rightarrow\infty$
and consider $f_n=f_{a_n,b_n}$ for all $n\in\N_0$.

\item Construct for every $n\in\N$ an horseshoe
$\Lambda_n$ adequate to $f_n$ and such that it gets closer to
$\Lambda_0$ as $n\rightarrow\infty$.

\item Collapse $\Lambda_n$ and consider the $\overline{f_n^R}$-invariant
densities $\bar{\rho}_n$. Realize them as functions defined on an
interval of the $x$-axis and belonging to a closed disk of
$L^\infty$. Apply Banach-Alaoglu Theorem to derive a convergent
subsequence $\bar{\rho}_{n_i}\rightarrow\bar{\rho}_{\infty}$.

\item Employ a technique used by Bowen in \cite{Bo75} to lift the
$\overline{f^R}$-invariant measure from the quotient space
$\bar{\Lambda}$ to an $f^R$-invariant measure on the horseshoe
$\Lambda$. This way we obtain measures $\tilde{\nu}_{n_i}$  and
$\tilde{\nu}_{\infty}$, defined on $\Lambda_{n_i}$ and
$\Lambda_0$, respectively.

\item Verify that all the measures
$\tilde{\nu}_{n_i}$ and $\tilde{\nu}_{\infty}$ desintegrate into
conditional absolutely continuous measures on unstable leaves.

\item Saturate the measures
$\tilde{\nu}_{n_i}$ and $\tilde{\nu}_{\infty}$. These saturations
are $f_{n_i}$-invariant  and $f_0$-invariant, respectively, and
have absolutely continuous conditional measures on unstable
leaves. The uniqueness of the SRB measures ensures that the
saturation of $\tilde{\nu}_{n_i}$ is $\nu_{n_i}$ (the
$f_{n_i}$-invariant SRB measure) and that of $\tilde{\nu}_\infty$
is $\nu_0$ (the $f_0$-invariant SRB measure).

\item Finally, show that this construction yields
$\nu_{n_i}\rightarrow\nu_0$ in the weak* topology.

\end{itemize}

\section{Dynamics of H\'{e}non maps on Benedicks-Carleson parameters}
\label{sec:dynamics-henon-maps}

In this section we provide information regarding the dynamical
properties of the H\'{e}non maps $f=f_{a,b}$, corresponding to the
Benedicks-Carleson parameters $(a,b)\in\B$. We do not intend to
give an exhaustive description but rather a brief summary of the
most relevant features whose main ideas are scattered through the
papers \cite{BC91,BY93,MV93, BY00}. We recommend the summary in
\cite{BY93} and Chapter 4 of \cite{BDV05} where the reader can
find a comprehensive description of the techniques and results
regarding H\'{e}non-like maps, including a revision of the
referred papers; both texts inspired our summary. The survey
\cite{LV03} provides a deep discussion about the exclusion of
parameters which are the basis of Benedicks-Carleson results.
Concerning the 1-dimensional case we also refer the paper
\cite{Fr05} in which a description of the Benedicks-Carleson
techniques in the phase space setting can be found.

\subsection{One-dimensional model}
\label{subsec:1-d-model}

The pioneer work of Jakobson \cite{Ja81} establishing the
existence of a positive Lebesgue measure set of parameters where
the logistic family presents chaotic behavior paved the way for a
better understanding of the dynamics beyond the non-hyperbolic
case. The analysis of the H\'{e}non maps made by Benedicks and
Carleson, triggered by the work of Collet-Eckmann
\cite{CE80a,CE80b} and Benedicks-Carleson \cite{BC85} themselves,
was a major breakthrough in that direction.  A key idea is the
exponential growth of the derivative along the critical orbit,
introduced in \cite{CE83}. In their remarkable paper \cite{BC91},
Benedicks and Carleson manage to establish, in a very creative
fashion, a parallelism between the estimates for the 1-dimensional
quadratic maps and the H\'{e}non maps. This connection supports
the use of 1-dimensional language in the present paper and compels
us to remind the results in Section 2 of \cite{BC91}. In there, it
is proved the existence of a positive Lebesgue measure set of
parameters, say $\B_1$ \index{aabc1@$\B_1$},
\index{Benedicks-Carleson!set of parameters!quadratic maps} within
the family $f_a:[-1,1]\rightarrow[-1,1]$, given by $f_a(x)=1-ax^2$
verifying
\begin{enumerate}

\item
\label{item:EG-1-dim} there is $c>0$ ($c\approx \log2$) such that
$|Df_a^n(f_a(0))|\geq \e^{cn}$ for all $n\geq0$; \index{aac@$c$}

\item
\label{item:BA-1-dim}
 there is a small $\alpha>0$ such that $|f_a^n(0)|\geq\e^
{-\alpha n}$ for all $n\geq 1$.\index{aaalpha@$\alpha$}

\end{enumerate}

The idea, roughly speaking, is that while the orbit of the
critical point is outside a critical region we have expansion (see
Subsection \ref{subsubsec:free-period-1d}); when it returns we
have a serious setback in the expansion but then, by continuity,
the orbit repeats its early history regaining expansion on account
of \eqref{item:EG-1-dim}. To arrange for \eqref{item:EG-1-dim} one
has to guarantee that the losses at the returns are not too
drastic hence, by parameter elimination, \eqref{item:BA-1-dim} is
imposed. The argument is mounted in a very intricate induction
scheme that guarantees both the conditions for the parameters that
survive the exclusions.

We focus on the maps corresponding to Benedicks-Carleson
parameters and study the growth of $Df_a^n(x)$ for $x\in[-1,1]$
and $a\in\B_1$. For that matter we split the orbit in free periods
and bound periods. During the former we are certain that the orbit
never visits the critical region. The latter begin when the orbit
returns to the critical region and initiates a bound to the
critical point, accompanying its early iterates. \index{free
period!quadratic maps} We describe the behavior of the derivative
during these periods in Subsections \ref{subsubsec:free-period-1d}
and \ref{subsubsec:bound-period-1d}.

The critical region is the interval
$(-\delta,\delta)$,\index{aadelta@$\delta$, $\Delta$} where
$\delta=\e^{-\Delta}>0$ is chosen small but much larger than
$2-a$. This region is partitioned into the intervals
\index{critical region!quadratic maps}
\[
(-\delta,\delta)=\bigcup_{m\geq\Delta} I_m,
\]
where $I_m=(\e^{-(m+1)},\e^{-m}]$ for $m>0$ and $I_{-m}=-I_m$ for
$m<0$; then each $I_m$ is further subdivided into $m^2$ intervals
$\{I_{m,j}\}$ \index{aaim@$I_m$, $I_{m,j}$}of equal length
inducing the partition $\P$\index{aaPg@$\P$} of $[-1,1]$ into
\begin{equation}
\label{eq:partition} [-1,-\delta)\cup
\bigcup_{m,j}I_{m,j}\cup(-\delta,1].
\end{equation}
Given $J\in\P$, we let $nJ$ denote the interval $n$ times the
length of $J$ centered at $J$.

\subsubsection{Expansion outside the critical region}
\label{subsubsec:free-period-1d} There is $c_0>0$ and $M_0\in\N$
such that \index{aac0@$c_0$}
\begin{enumerate}

\item \label{item:free-period-M0-1d} If
$x,\ldots,f_a^{k-1}(x)\notin(-\delta,\delta)$ and $k\geq M_0$,
then $|Df_a^k(x)|\geq\e^{c_0k}$;

\item \label{item:free-period-return-1d} If
$x,\ldots,f_a^{k-1}(x)\notin(-\delta,\delta)$ and
$f_a^k(x)\in(-\delta,\delta)$, then $|Df_a^k(x)|\geq\e^{c_0k}$;

\item \label{item:free-period-1d} If
$x,\ldots,f_a^{k-1}(x)\notin(-\delta,\delta)$, then
$|Df_a^k(x)|\geq\delta\e^{c_0k}$.

\end{enumerate}

\subsubsection{Bound period definition and properties}
\label{subsubsec:bound-period-1d} Let
$\beta=14\alpha$\index{aabeta@$\beta$}. For $x\in(-\delta,\delta)$
define $p(x)$ to be the largest integer $p$\index{aap@$p$} such
that
\begin{equation}
\label{eq:def-bound-period} |f_a^k(x)-f_a^k(0)|<\e^{-\beta
k},\qquad \forall k<p. \index{bound period!quadratic maps}
\end{equation}Then

\begin{enumerate}

\item \label{item:bound-period-bounds-1d}$\frac{1}{2}|m|\leq p(x)\leq3|m|$,
for each $ x\in I_m$;

\item \label{item:bound-period-derivative-1d} $|Df_a^p(x)|\geq
\e^{c'p}$, where $c'=\frac{1-4\beta}{3}>0$.

\end{enumerate}
The orbit of $x$ is said to be bound to the critical point during
the period $0\leq k<p$. We may assume that $p$ is constant on each
$I_{m,j}$.

\subsubsection{Distortion of the derivative}
\label{subsubsec:distortion-1d} The partition $\P$ is designed so
that if $\omega\subset[-1,1]$ is such that, for all $k<n$,
$f^k(\omega)\subset 3J$ for some $J\in\P$, then there exists a
constant $C$ independent of $\omega$, $n$ and the parameter so
that for every $x,y\in\omega$,
\[
\frac{|Df_a^n(x)|}{|Df_a^n(y)|}\leq C.
\]
\index{bounded distortion!quadratic maps}
\subsubsection{Derivative estimate}
\label{subsubsec:derivative-estimate-1d} Suppose that
\begin{equation}
\label{eq:SAn-1d} |f_a^j(x)|\geq\delta\e^{-\alpha j},\qquad
\forall j<n.
\end{equation}
Then there is a constant $c_2>0$ \index{aac2@$c_2$}such that
\begin{equation}
\label{eq:EEn-1d} |Df_a^n(x)|\geq\delta\e^{c_2n}.
\end{equation}
A proof of this fact can be found in  \cite[Section 3]{Fr05} where
it is also shown that there is $\kappa>0$ such that
\[
|\{x\in[-1,1]:\,|f_a^j(x)|\geq \e^{-\alpha j},\,\forall j<n\}|\geq
2-\mbox{const}\,\e^{-\kappa n}.
\] \index{derivative estimate!quadratic maps}
As an easy consequence, it is deduced that Lebesgue almost every
$x$ has a positive Lyapunov exponent. Moreover, we have a positive
Lebesgue measure set of points $x\in[-1,1]$ satisfying
\eqref{eq:SAn-1d}, and so \eqref{eq:EEn-1d}, for all $n\in\N$.

\subsection{General description of the H\'{e}non attractor}

The following facts are elementary for $f=f_{a,b}$ with $(a,b)$
inside an open set of parameters.

Each $f$ has a unique fixed point in the first quadrant
$z^*\approx\left(\tfrac{1}{2},\tfrac{1}{2}b\right)$
\index{aaz*@$z^*$}. This fixed point is hyperbolic with an
expanding direction presenting a slope of order $-b/2$ and a
contractive direction with a slope of approximately $2$. The
respective eigenvalues are approximately $-2$ and $b/2$. In
\cite{BC91} it is shown that if we choose $a_0<a_1<2$ with $a_0$
sufficiently near 2, then there exists $b_0$ sufficiently small
when compared to $2-a_0$ such that for all $(a,b)\in
[a_0,a_1]\times(0,b_0]$, the unstable manifold of $z^*$, say $W$,
never leaves a bounded region. Moreover, its closure
$\overline{W}$ is an attractor in the sense that there is an open
neighborhood $U$ of $\overline{W}$ such that for every $z\in U$ we
have $f^n(z)\rightarrow \overline{W}$ as $n\rightarrow\infty$.
\index{aaWg@$W$}

\subsubsection{Hyperbolicity outside the critical region}
\label{subsubsec:hyperb-out-crit-region} Let $\delta$ be at least
as small as in our 1-dimensional analysis and assume that
$b_0\ll2-a_0\ll\delta$. The critical region is now
$(-\delta,\delta)\times\R$.\index{critical region!H\'enon maps} A
simple calculation shows that outside the critical region $Df$
preserves the cones $\{|s(v)|\leq\delta\}$ (see \cite{BY93}
Subsection 1.2.3), where $s(v)$ denotes the slope of the vector
$v$. For $z=(x,y)\notin(-\delta,\delta)\times\R$ and a unit vector
$v$ with $s(v)\leq\delta$, we have essentially the same estimates
as in 1-dimension. That is, there is $c_0>0$ and $M_0\in\N$ such
that
\begin{enumerate}

\item \label{item:free-period-M0} If
$z,\ldots,f^{k-1}(z)\notin(-\delta,\delta)\times\R$ and $k\geq
M_0$ then $|Df^k(z)v|\geq\e^{c_0k}$;

\item \label{item:free-period-return} If
$z,\ldots,f^{k-1}(z)\notin(-\delta,\delta)\times\R$ and
$f^k(z)\in(-\delta,\delta)\times\R$ then
$|Df^k(z)v|\geq\e^{c_0k}$;

\item \label{item:free-period} If
$z,\ldots,f^{k-1}(z)\notin(-\delta,\delta)\times\R$ then
$|Df^k(z)v|\geq\delta\e^{c_0k}$.
\end{enumerate}

\subsection{The contractive vector field}
\label{subsec:contractive-vector-field} For $A\in\mbox{GL}(2,\R)$
and a unit vector $v$, if $v\mapsto|Av|$ is not constant, let
$e(A)$ denote the unit vector maximally contracted by $A$. We will
write $e_n(z):=e(Df^n(z))$ whenever it makes sense.
\index{aaen0@$e_n(z)$}Observe that if we have some sort of
expansion in $z$, say $|Df^n(z)v|>1$ for some vector $v$, then
$e_n(z)$ is defined and $|Df^n(z)e_n(z)|\leq b^n$ since
$\det(Df^n(z))=(-b)^n$. \index{contractive vector field}

The following general perturbation lemma is stated in \cite{BY00}
and clarifies the assertions of Lemma 5.5 and Corollary 5.7 in
\cite{BC91}, where the proofs can be found. Given
$A_1,A_2,\ldots$, we write $A^n:=A_n\ldots A_1$; all the matrices
below are assumed to have determinant equal to $b$.

\begin{lemma}[Matrix Perturbation Lemma] Given $\kappa\gg b$,
exists $ \lambda$ with $b\ll\lambda<\min(1,\kappa)$
\index{aalambda@$\lambda$} such that if $A_1\ldots,A_n$,
$A_1'\ldots,A_n'\in\mbox{GL}(2,\R)$ and $v\in\R^2$ satisfy
\[
|A^iv|\geq\kappa^i\qquad \mbox{and} \qquad
\|A_i-A_i'\|<\lambda^i\quad\forall i\leq n,
\]
then we have, for all $i\leq n$:
\begin{itemize}

\item $|A'^iv|\geq\frac{1}{2}\kappa^i$;

\item $\sphericalangle(A^iv,A'^iv)\leq\lambda^{\tfrac{i}{4}}$.
\end{itemize}
\end{lemma}\index{matrix perturbation lemma}
From the Matrix Perturbation Lemma, it follows that  if for some
$\kappa$ and $v$, we have $|Df^j(z_0)v|\geq \kappa^j$ for all $
j\in\{0,\ldots,n\}$, then there is a ball of radius
$(\lambda/5)^n$ about $z_0$ on which $e_n$ is defined and $|Df^n
e_n|\leq2(b/\kappa)^n$. Assuming that $\kappa$ is fixed and $e_n$
is defined in a ball $B_n$ around $z_0$ the following facts hold
(see \cite[Section 5]{BC91}, \cite[Section 1.3.4]{BY93} or
\cite[Section 1.5]{BY00}):
\begin{enumerate}\index{contractive vector field!properties}

\item \label{item:cont-field-e1} $e_1$ is defined everywhere and has
slope equal to $2ax+\O(b)$;

\item \label{item:cont-field-Lip-2d} there is a constant $C>0$ such that for all
 $z_1,z_2\in B_n$,
\[
|e_n(z_1)-e_n(z_2)|\leq C|z_1-z_2|;
\]

\item \label{item:cont-field-Lip-x-dir} for $z_1=(x_1,y_1)$,
$z_2=(x_2,y_2)\in B_n$ with $|y_1-y_2|\leq|x_1-x_2|$
\[
|e_n(z_1)-e_n(z_2)|=(2a+\O(b))|x_1-x_2|;
\]

\item \label{item:cont-field-em-en} for $m<n$, $|e_n-e_m|\leq \O(b^m)$
on $B_n$.
\end{enumerate}

From this point onward we restrict ourselves to the H\'{e}non maps
of the Benedicks-Carleson type, that is, we are considering
$f=f_{a,b}$ for $(a,b)\in\B$.

\subsection{Critical points}
\label{subsec:Critical-points} The cornerstone of
Benedicks-Carleson strategy is the critical set in $W$ denoted by
$\C$ \index{aaCg@$\C$}, that plays the role of the critical point
$0$ in the 1-dimensional model. The critical points correspond to
homoclinic tangencies of Pesin stable and unstable manifolds. For
$z\in W$, let $\tau(z)\in T_z\R^2$ \index{aatau@$\tau$} denote a
unit vector tangent to W at $z$. For each $\zeta\in\C$, the vector
$\tau(\zeta)$ is contracted by both forward and backward iterates
of the derivative. In fact, we have
$\lim_{n\rightarrow\infty}e_n(\zeta)=\tau(\zeta)$, which can be
thought as the moral equivalent to $Df(0)=0$ in 1-dimension. The
following subsections refer to \cite{BC91}, mostly Sections 5 and
6 (see also \cite[Section 1.3.1]{BY93}).\index{critical points}

\subsubsection{Rules for the construction of the critical set}
\label{subsubsec:rules-critical-set} \index{critical points!rules
of construction}The critical set $\C$ is
located in $W\cap(-10b,10b)\times\R$. 
There is a unique $z_0\in\C$ on the roughly horizontal segment of
$W$ containing the fixed point $z^*$. The part of $W$ between
$f^2(z_0)$ and $f(z_0)$ is denoted by $W_1$ and called the leaf of
generation 1.\index{aaWg1@$W_1$} Leaves of generation $g\geq2$ are
defined by $W_g:=f^{g-1}W_1\setminus\bigcup_{j\leq g-1}W_j$. We
assume that $(a,b)$ is sufficiently near $(2,0)$ so that
$\bigcup_{g\leq27}W_g$ consists of $2^{26}$ roughly horizontal
segments linked by sharp turns near $x=\pm1$, $y=0$, and that
$\bigcup_{g\leq27}W_g\cap(-\delta,\delta)\times\R$ consists of
$2^{26}$ curves whose slope and curvature are $\leq10b$ -- in
\cite{BC91} such a curve is called $C^2(b)$. In each of them there
is a unique critical point\index{generation}

For $g>27$, assume that all critical points of generation $\leq
g-1$ are already defined. Consider a maximal piece of $C^2(b)$
curve $\gamma\subset W_g$. If $\gamma$ contains a segment of
length $2\varrho^g$ centered at $z=(x,y)$, where $\varrho$
\index{aarho2@$\varrho$} verifies $b\ll\varrho\ll\e^{-72}$, and
there is a critical point $\tilde{z}=(\tilde{x},\tilde{y})$ of
generation $\leq g-1$ with $x=\tilde{x}$ and $|y-\tilde{y}|\leq
b^{g/540}$, then a unique critical point $z_0\in\C\cap\gamma$ of
generation $g$ is created satisfying the condition $|z_0-z|\leq
|y-\tilde{y}|^{1/2}$. These are the only critical points of
generation $g$.

Observe that the exact position of a critical point is
unaccessible since its definition depends on the limiting relation
$\lim_{n\rightarrow\infty}e_n(\zeta)=\tau(\zeta)$. So the strategy
in \cite{BC91} is to produce approximate critical points $\zeta^n$
of increasing order which are solutions of the equation
$e_n(z)=\tau(z)$. Once an approximate critical point is born,
parameters are excluded to ensure that a critical point
$\zeta\in\C$ is created nearby. Moreover,
$|\zeta^n-\zeta|=\O(b^n)$.\index{critical points!approximate
critical points}

\subsubsection{Dynamical properties of the critical set}
\label{subsubsec:dynamic-prop-crit-set} The parameter exclusion
procedure leading to $\B$ is designed so that every $z\in\C$ has
the following properties:

\begin{itemize}
\item
 there is $c\approx\log2$ and $C$
independent of $b$ such that for all $n\in\N_0$,
\begin{equation}
\tag{$UH$}\label{eq:UH}\index{critical points!uniform
hyperbolicity (UH)}
|Df^n(z)\bigl(\begin{smallmatrix}0\\1\end{smallmatrix}\bigr)
|\geq\e^{cn} \qand
|Df^n(z)\tau|\leq(Cb)^n;
\end{equation}

\item there is a small number $\alpha>0$, say $\alpha=10^{-6}$, such
that for all $n\in\N$
\begin{equation}
\tag{$BA$} \label{eq:BA} \mbox{dist}(f^n(z),\C)\geq\e^{-\alpha
n}.\index{critical points!basic assumption (BA)}
\end{equation}

\end{itemize}
The precise meaning of ``$\mbox{dist}$'' in the last equation will
be described in Section \ref{subsubsec:tang-pos-dist-cs}.
 The uniform hyperbolicity expressed in \eqref{eq:UH} is analogous
to condition \eqref{item:EG-1-dim} in Subsection
\ref{subsec:1-d-model} while the basic assumption stated in
\eqref{eq:BA} is the surface analogue to condition
\eqref{item:BA-1-dim} of the 1-dimensional model.

One of the reasons why the Benedicks-Carleson proof is so involved
is that in order to $e_n$ be defined in the vicinity of critical
points, one has to require some amount of hyperbolicity which is
exactly what one wants to achieve (see \eqref{eq:UH} above). This
difficulty is overcome by working with finite time approximations
and imposing slow recurrence in a very intricate induction scheme.
Once an approximate critical point $\zeta^n$ of order $n$  is
designated one studies its orbit. When it comes near $\C$, there
is a near-interchange of stable and unstable directions -- hence a
setback in hyperbolicity. But then the orbit of $\zeta^n$ follows
for some time the orbit of some $\tilde{\zeta}\in\C$ of earlier
generation and regains hyperbolicity on account of \eqref{eq:UH}
for $\tilde{\zeta}$. To arrange \eqref{eq:UH} at time $n+1$ for
$\zeta^n$, it is necessary to keep the orbits from switching
stable and unstable directions too fast, so by parameter exclusion
we impose \eqref{eq:BA}. At this stage it is possible to define
$e_{n+1}$ and thus find a critical point approximation of order
$n+1$, denoted by $\zeta^{n+1}$. The information is updated and
the process is repeated. Fortunately, a positive Lebesgue measure
set of parameters survives the exclusions.

\subsection{Binding to critical points}
\label{subsec:binding-crit-point} The critical point $0$ in the
1-dimensional context plays a dual role. Firstly, the distance to
the critical point is a measure of the norm of the derivative,
which is the reason why a recurrence condition like
\eqref{item:BA-1-dim} of Subsection \ref{subsec:1-d-model} can be
used to bound the loss of expansion when an orbit comes near the
critical point and to obtain the exponential growth expressed in
\eqref{item:EG-1-dim} of Subsection \ref{subsec:1-d-model}.
Secondly, during the bound, period information of the early
iterates of the critical point is passed through continuity to the
points returning to the critical region. In order to replicate
this in the H\'{e}non family, for every return time $n$ of the
orbit of $z\in W$ ($z$ may belong to $\C$) we must associate a
suitable binding critical point for $f^n(z)$ so that we can have
some meaning of the distance of $f^n(z)$ to the critical set.
The suitability depends on the validity of two requirements:
tangential position and correct splitting.

\subsubsection{Tangential position and distance to the critical
set} \label{subsubsec:tang-pos-dist-cs} Let $z\in W$ and $n$ be
one of its return time to the critical region. Let $\zeta\in\C$.
Essentially we say that $f^n(z)$ is in tangential position with
respect to $\zeta$ if its horizontal distance to $\zeta$ is much
larger than the vertical distance. In fact we will use the notion
of generalized tangential positions introduced in \cite[Section
1.6.2]{BY93}
 instead of the original one from \cite{BC91} (see
\cite[Section 1.4.1]{BY93}). For $z\in W$ we say that $(x',y')$ is
the natural coordinate system at $z$ if $(0,0)$ is at $z$, the
$x'$ axis is aligned with $\tau(z)$ and the $y'$ axis with
$\tau(z)^\perp$.
\begin{definition}
\label{def:tangential-pos} Let $c>0$ be a small number much less
than $2a$, say $c=10^{-2}$, and let $\zeta\in\C$. A point $z$ is
said to be in {\em tangential position} with respect to $\zeta$,
if $z=(x',y')$ with $|y'|\leq cx'^2$, in the natural coordinate
system at $\zeta$.
\end{definition}\index{tangential position}

In \cite[Section 7.2]{BC91} it is arranged that for every
$\zeta\in\C$ and any $n$-th return to the critical region,
 there is a critical point $\hat{\zeta}$ of earlier generation
with respect to which $f^n(\zeta)$ is in tangential position. This
is done through an argument known as the capture procedure (see
also \cite[ Section 2.2.2]{BY93}) which essentially consists in
showing that when a critical orbit $\zeta\in\C$ experiences a free
return at time $n$, then $f^n(\zeta)$ is surrounded by a fairly
regular collection of $C^2(b)$ segments $\{\gamma_j\}$ of $W$
which are relatively long and of earlier generations. In fact, we
have $\mbox{gen}(\gamma_j)\approx 3^j$,
$\mbox{length}(\gamma_j)\approx\varrho^{3^j}$ and
$\mbox{dist}(f^n(z),\gamma_j)\approx b^{3^j}$, where $3^j<\theta
n$ and $\theta\approx\frac{1}{|\log b|}$. Some (maybe all) of
these captured segments will have critical points and most
locations of $f^n(\zeta)$ will be in tangential position with
respect to one of these critical points. Bad locations of
$f^n(\zeta)$ correspond to deleted parameters. This is another
subtlety of Benedicks-Carleson proof: every time a critical point
is created it causes a certain amount of parameters to be
discarded so we cannot afford to have too many critical points;
however, we must have enough critical points so that a convenient
one, in tangential position, may be found every time a return
occurs.

In \cite{BY93} it is shown that this kind of control when a
critical orbit returns can be extended to all points in $W$. Thus,
for any return of the orbit of $z\in W$ to the critical region
there is an available binding critical point with respect to which
the tangential position requirement holds. In fact \cite[Lemma
7]{BY93} guarantees that one can systematically assign to each
maximal free segment $\gamma\subset W$ intersecting the critical
region a critical point $\tilde{z}(\gamma)$ with respect to which
each $z\in\gamma$ are in tangential position. When the orbit of
$z\in W$ returns to the critical region, say at time $n$, we
denote by $z(f^n(z))\in\C$ a critical point with respect to which
$f^n(z)$ is in tangential position.

These facts lead us to the notion of distance to the critical set.
We do not intend to give a formal definition but rather introduce
a concept that gives an indication of closeness to the critical
set. In \cite{BC91} and \cite{BY93} two different perspectives of
distance to the critical set have been introduced. In
\cite[Section 2]{BY00} this notion is cleaned up and these two
different perspectives are seen to translate essentially the same
geometrical facts. Let $z\in W$. \index{distance to the critical
set}If $z=(x,y)\notin (-\delta,\delta)\times \R$ we consider that
$\mbox{dist}(z,\C)=|x|$; if $z\in (-\delta,\delta)\times\R$ then
we pick any critical point $\zeta\in\C$ with respect to which $z$
is in tangential position and let $\mbox{dist}(z,\C)=|z-\zeta|$.
In order to this notion make sense one has to verify that if
$\hat{\zeta}\in\C$ is a different critical point with respect to
which $z$ is also in tangential position then
$|z-\zeta|\approx|z-\hat{\zeta}|$. This is exactly the content
 of \cite[Lemma 1']{BY00}, where it is proved that
${|z-\zeta|}/{|z-\hat{\zeta}|}=1+\O(\max(b,d^2))$, for
$d=\min(|z-\zeta|,|z-\hat{\zeta}|)$. As observed in \cite{BY00}
for a better understanding of the distance of a given point $z\in
W\cap(-\delta,\delta)\times\R$ to the critical set, one should
look at the angle between $\tau(z)$ and $e_m(z)$, the most
contracted vector at $z$ of a convenient order $m$. The reason for
this is that, at the critical points, this angle is extremely
close to $0$; actually  it tends to $0$ if we let $m$ go to
infinity.\index{distance to the critical set}

\subsubsection{Bound period and fold period}
\label{subsubsec:bound-fold-period} Let $z\in
W\cap(-\delta,\delta)\times\R$ be in tangential position with
respect to $\zeta\in\C$. Then $z$ initiates a binding to $\zeta$
of length $p$, where $p=p(z,\zeta)$ \index{aap@$p$} is the largest
$k$ such that
\[
|f^j(z)-f^j(\zeta)|<\e^{-\beta j},\qquad \forall j<k
\]\index{bound period!H\'enon maps}
where $\beta=14\alpha$. We say that in the next $p$ iterates, $z$
is bounded to $\zeta$. It is convenient to modify slightly the
above definition of $p$ so that the bound periods become nested.
This means that if the orbit of $z$ returns to the critical region
before $p$ then the bound period initiated at that time must cease
before the end of the bound relation to $\zeta$. This is done in
 \cite[Section 6.2]{BC91}. It is further required that if the
 bound relation between $z$ and
$\zeta$ is still in effect at time $n$, which is a return time for
both, then $z(f^n(\zeta))=z(f^n(z))$.

An additional complication arises in the  H\'{e}non maps: the
folding. To illustrate it, let $\gamma\subset W$ be a $C^2(b)$
segment containing a critical point $\zeta$. The practically
horizontal vector $\tau(\zeta)$ will be sent by $Df$ into an
approximately vertical direction, which is the typical contracting
direction of the system, and will be contracted forever. After few
iterations $\gamma$ develops very sharp bends at the iterates of
$\zeta$, which induce an unstable setting near the bends. In fact,
if we pick a point $z\in\gamma$ very close to $\zeta$, its
iterates diverge very fast from the bends which means that after
some time, say $n$, depending on how close $z$ and $\zeta$ are,
the vector $\tau(f^n(z))$ will be practically aligned with the
horizontal direction again, which, on the contrary, is the typical
expanding direction of the system. The interval of time that the
tangent direction takes to be horizontal again is called the fold
period. \index{fold period}

The actual definition of fold period is given  in \cite[Sections
6.2 and 6.3]{BC91}; here, we stick to the previous heuristic
motivation and to the following properties. If $z\in W$ has a
return at time $n$, the fold period of $f^n(z)$ with respect to
$z(f^n(z))\in\C$ is a positive integer $l=l(f^n(z),z(f^n(z)))$
\index{aal@$l$}such that
\begin{enumerate}

\item \label{item:fold-period-bounds} $2m\leq l\leq 3m$, where
$(5b)^m\leq|f^n(z)-z(f^n(z))|\leq(5b)^{m-1}$;

\item \label{item:fold-period-compar-bp} $l/p\leq \mbox{const}/|\log
b|$, that is the fold period associated to a return is very short
when compared to the bound period initiated at that time.

\end{enumerate}

\subsubsection{Correct splitting and controlled orbits}
\label{subsubsec:correct-split}In order to duplicate the
1-dimensional behavior not only one assigns a binding critical
point every time a return to the critical region occurs but also
one would like to guarantee that the loss of hiperbolicity due to
the return is in some sense proportional to the distance to the
critical set. This is achieved through the notion of correct
splitting.
\begin{definition}
\label{def:correct-split} Let $z\in W$, $v\in T_z\R^2$, $n\in\N$
be a return time for $z$ and consider $z(f^n(z))\in\C$ with
respect to which $f^n(z)$ is in tangential position. We say that
the vector $Df^n(z)v$ {\em splits correctly} with respect to
$z(f^n(z))\in\C$ if and only if we have that\index{correct
splitting}
\[
3|f^n(z)-z(f^n(z))|\leq\sphericalangle(Df^n(z)v,e_l(f_n(z)))
\leq5|f^n(z)-z(f^n(z))|,
\]
where $l$ is the fold period associated to the return.
\end{definition}
Now we are in condition of defining controlled orbits.
\begin{definition}
Let $z\in W$ and $v\in T_z\R^2$ and $N\in\N$. We say that the pair
$(z,v)$ is controlled on the time interval $[0,N)$ if for every
return $n\in[0,N)$ of the orbit of $z$ to the critical region,
there is $z(f^n(z))\in\C$ with respect to which $f^n(z)$ is in
tangential position and $Df^n(z)v$ splits correctly with respect
to $z(f^n(z))\in\C$. We say that the pair $(z,v)$ is {\em
controlled} during the time interval $[0,\infty)$ if it is
controlled on $[0,N)$ for every $N\in\N$.\index{controlled orbits}
\end{definition}

One of the most important properties of $f$ proved in \cite{BC91}
is that for every $\zeta\in\C$, the pair
$(\zeta,\bigl(\begin{smallmatrix}0\\1\end{smallmatrix}\bigr))$ is
controlled during the time interval $[0,\infty)$. This fact
supports the validity of the 1-dimensional estimates in the
surface case.

We say that the pair $(z,v)$ is controlled on $[j,0)$ with
$-\infty<j<0$, if $(f^{j}(z),Df^j(z)v)$ is controlled on $[0,-j)$
and that $(z,v)$ is controlled on $(-\infty,0)$ if it is
controlled on $[j,0)$ for all $j<0$. In  \cite[Proposition
1]{BY93} it is proved that if the orbit of $z\in W$ never hits the
critical set $\C$ then the pair $(z,\tau(z))$ is controlled in the
time interval $(-\infty,\infty)$.

\subsection{Dynamics in $W$}
\label{subsec:Dynamics-in-W}

As referred, \cite[Proposition 1]{BY93} shows that every orbit of
$z\in W$ can be controlled using those of $\C$, just as it was
done for critical orbits in \cite{BC91}. This means that each
orbit in $W$ can be organized into free periods and bound periods.
To illustrate, consider $z$ belonging to a small segment of $W$
around the fixed point $z^*$. By definition $z$ is considered to
be free at this particular time. The first forward iterates of $z$
are also in a free state \index{free period!H\'enon maps}, until
the first return to the critical region occurs, say at time $n$.
Then since the pair $(z,\tau(z))$ is controlled there is
$z(f^n(z))$ with respect to which $f^n(z)$ is in tangential
position and $Df^n(z)\tau$ splits correctly. During the next $p$
iterates we say that $z$ is bound to the critical point
$z(f^n(z))$. If $f^n(z)\in\C$ then the bound period is infinite;
otherwise, after the time $n+p$ the iterates of $z$ are said to be
in free state once again and history repeats itself.

This division of the orbits into free periods, bound periods and
the special design of the control of orbits through the tangential
position and correct splitting requirements allowed \cite{BY93} to
recover the one dimensional estimates. In fact, the loss of
expansion at the returns is somehow proportional to the distance
to the binding critical point and it is completely overcome at the
end of the bound period.

The following estimates, unless otherwise mentioned, are proved in
 \cite[Corollary~1]{BY93}.

\begin{enumerate}

\item \label{item:free-period-estimates} Free period estimates.
    \begin{enumerate}

    \item \label{item:fpe-slope} Every free segment $\gamma$ has
    slope less than
    $2b/\delta$, and $\gamma\cap(-\delta,\delta)\times\R$ is a $C^2(b)$
    curve (Lemmas 1 and 2 of \cite{BY93});

    \item \label{item:fpe-M0}There is $c_0>0$ and $M_0\in\N$ such that
    if $z$ is free and $z,\ldots,f^{k-1}(z)\notin(-\delta,\delta)\times\R$
    with $k\geq
    M_0$ then $|Df^k(z)\tau|\geq\e^{c_0k}$;

    \item \label{item:fpe-return}There is $c_0>0$ such that if
    $z$ is free,
    $z,\ldots,f^{k-1}(z)\notin(-\delta,\delta)\times\R$ and
    $f^k(z)\in(-\delta,\delta)\times\R$ then
    $|Df^k(z)\tau|\geq\e^{c_0k}$.

    \end{enumerate}

\item \label{item:bound-period-estimates} Bound period
estimates.\newline
    There is $c\approx\log2$ such that if
    $z\in(-\delta,\delta)\times\R$ is free and initiates a binding
    to $\zeta\in\C$ with bound period $p$, then
    \begin{enumerate}

    \item \label{item:bpe-length} If
    $\e^{-m-1}\leq|z-\zeta|\leq\e^{-m}$, then
    $\frac{1}{2}m\leq p\leq5m$;

    \item \label{item:bpe-deriv-at-middle}
    $|Df^j(z)\tau|\geq|z-\zeta|\e^{cj}$ for $0<j<p$.

    \item \label{item:bpe-deriv-at-end}
    $|Df^p(z)\tau|\geq\e^{c\frac{p}{3}}$.

    \end{enumerate}

\item \label{item:orbits-end-free-states} Orbits ending in free states.
\newline There exists $c_1>\frac{1}{3}\log2$ \index{aac1@$c_1$}
such that if $z\in W\cap(-\delta,\delta)\times\R$ is in a free
state, then $|Df^{-j}(z)\tau|\leq\e^{-c_1j}$, for all $j\geq0$
(\cite[Lemma 3]{BY93}).

\end{enumerate}

\subsubsection{Derivative estimate}
\label{subsubsec:derivative-estimate}

The next derivative estimate can be found in
\cite[Section~1.4]{BY00}. It is the 2-dimensional analogue to the
1-dimensional derivative estimate expressed in Subsection
\ref{subsubsec:derivative-estimate-1d}. Consider $n\in\N$ and a
point $z$ belonging to a free segment of $W$ and satisfying, for
every $j<n$
\begin{equation}
\tag{$SA$} \label{eq:SAn} \mbox{dist}(f^j(z),\mathcal{C})\geq
\delta \e^{-\alpha j}.
\end{equation} Then there is a constant
$c_2>0$ such that
\begin{equation}
\tag{$EE$} \label{eq:EGn} |Df^n(z)\tau|\geq\delta\e^{c_2n}.
\index{aac2@$c_2$}
\end{equation}\index{derivative estimate!H\'enon maps}
Essentially this estimate is saying that if we have slow
approximation to the critical set (or, in other words, a
\eqref{eq:BA} type property), then we have exponential expansion
along the tangent direction to $W$.

\subsubsection{Bookkeeping and bounded distortion}
\label{subsubsec:bounded-distortion} For $x_0\in\R$, we let
$\P_{[x_0]}$ \index{aaPgx0@$\P_{[x_0]}$} denote the partition $\P$
defined in \eqref{eq:partition} after being translated from $0$ to
$x_0$. Similarly, if $\gamma$ is a roughly horizontal curve in
$\R^2$ and $z_0=(x_0,y_0)\in\gamma$, we let $\P_{[z_0]}$
\index{aaPgz0@$\P_{[z_0]}$}denote the partition of $\gamma$ that
projects vertically onto $\P_{[x_0]}$ on the $x$-axis. Once
$\gamma$ and $z_0$ are specified, we will use $I_{m,j}$ to denote
the corresponding subsegment of $\gamma$.

Let $\gamma\subset(-\delta,\delta)\times\R$ be a segment of $W$.
We assume that the entire segment has the same itinerary up to
time $n$ in the sense that:

\begin{itemize}

\item all $z\in\gamma$ are bound or free simultaneously at any
moment;

\item if $0=t_0<t_1<\ldots<t_q$ are the consecutive free return
times before $n$, then for all $j\leq q$ the entire segment
$f^{t_j}\gamma$ has a common binding point $\zeta_j\in\C$ and
$f^{t_j}\gamma\subset 5 I_{m,k}^j$ for some
$I_{m,k}^j\in\P_{[\zeta_j]}$.

\end{itemize}

Then there exists $C_1>0$ \index{aaCg1@$C_1$} independent of
$\gamma$ and $n$ such that for all $z_1,z_2\in\gamma$
\[
\frac{|Df^n(z_1)\tau|}{|Df^n(z_2)\tau|}\leq C_1.
\]\index{bounded distortion!H\'enon maps}
This result can be found in \cite[Proposition 2]{BY93}.

\subsection{Dynamical and geometric description of the critical set}
\label{subsec:dynamical-geometric-description-cs} The construction
of the critical set seems to be done according to a quite
discretionary set of rules. However, as observed in \cite{BY93}
there are certain intrinsic characterizations of $\C$.
 Corollary 1
of \cite{BY93} gives the following dynamical description of $\C$.
Let $z\in W$. Then
\[
z \mbox{ lies on a critical orbit
}\quad\Leftrightarrow\quad\limsup_{n\rightarrow\infty}|Df^n(z)\tau|<\infty
\quad\Leftrightarrow\quad\limsup_{n\rightarrow\infty}|Df^n(z)\tau|=0.
\]
In fact, $z\in\C$ if and only if $|Df^j(z)\tau|\leq\e^{-c_1|j|}$,
for all $j\in\Z$, i.e. the critical points correspond to the
tangencies of Pesin stable manifolds with $W$ which endow an
homoclinic type behavior.

The critical set $\C$ has also a nice geometric characterization.
Given $\zeta\in W$, $\kappa(\zeta)$ denotes the curvature of $W$
at $\zeta$. From the curvature computations in \cite[Section
7.6]{BC91} (see also \cite[Section 2.1.3]{BY93}) one gets that
$$z\in\C \quad\Leftrightarrow\quad \kappa(z)\ll 1 \qand
\kappa(f^n(z))>b^{-n},\; \forall n\in\N.$$ This means that one can
look at the critical points as the points that are sent into the
folds of $W$.

\subsection{SRB measures}
\label{subsec:SRB-measures}  We begin by giving a formal
definition of \emph{Sinai-Ruelle-Bowen measures} (SRB
\emph{measures}). Let $f:M\rightarrow M$ be an arbitrary $C^2$
diffeomorphism of a finite dimensional manifold and let $\mu$ be
an $f$ invariant probability measure on $M$ with compact support.
We will assume that $\mu$-a.e. point, there is a strictly positive
Lyapunov exponent. Under these conditions, the unstable manifold
theorem of Pesin \cite{Pe78} or Ruelle \cite{Ru79} asserts that
passing through $\mu$-a.e. $z$ there is an unstable manifold which
we denote by $\gamma^u(z)$.

 A measurable partition $\L$ of $M$ is said to be
\emph{subordinate} to $\gamma^u$ (with respect to the measure
$\mu$) if at $\mu$-a.e. $z$, $\L(z)$ is contained in $\gamma^u(z)$
and contains an open neighborhood of $z$ in $\gamma^u(z)$, where
$\L(z)$ denotes the atom of $\L$ containing $z$. By Rokhlin's
desintegration theorem there exists a family $\{\mu_z^\L\}$ of
conditional measures of $\mu$ with respect to the partition $\L$
(see for example \cite[Appendixes C.4 and C.6]{BDV05}).

\begin{definition}
\label{def:SRB} Let $f:M\rightarrow M$ and $\mu$ be as above. We
say that $\mu$ is an SRB probability measure if for every
measurable partition $\L$ subordinate to $\gamma^u$, we have that
$\{\mu_z^\L\}$ is absolutely continuous with respect to Lebesgue
measure in $\gamma^u(z)$ for $\mu$-a.e. $z$.\index{SRB measures}
\end{definition}

In \cite{BY93} it is proved that $f_{a,b}$ admits an SRB measure
$\nu_{a,b}$, for every $(a,b)\in\B$. Moreover, $\nu_{a,b}$ is
unique (hence ergodic), it is a physical measure, its support is
$\overline W_{a,b}$
 and $(f_{a,b},\nu_{a,b})$ is
isomorphic to a Bernoulli shift.

\section{A horseshoe with positive measure}
\label{sec:horseshoe}

In order to obtain decay of correlations for H\'{e}non maps of the
Benedicks-Carleson type, Benedicks and Young build, in
\cite{BY00}, a set $\Lambda$ of positive SRB-measure with good
hyperbolic properties. $\Lambda$ has hyperbolic product structure
and it may be looked at as a horseshoe with infinitely many
branches and unbounded return times; it is obtained by
intersecting two families of $C^1$ stable and unstable curves.
Dynamically, $\Lambda$ can be decomposed into a countable union of
$s$-sublattices, denoted $\Xi_i$, crossing $\Lambda$ completely in
the stable direction, with a Markov type property: for each
$\Xi_i$ there is $R_i\in\N$ such that $f^{R_i}(\Xi_i)$ is an
$u$-sublattice of $\Lambda$, crossing $\Lambda$ completely in the
unstable direction. The intersection of $\Lambda$ with every
unstable leaf is a positive 1-dimensional Lebesgue measure set.
Before continuing with an overview of the construction of such
horseshoes, we mention that Young \cite{Yo98} has extended the
argument in \cite{BY00} to a wider setting and observed that
similar horseshoes can be found in other situations. We will refer
to \cite{Yo98} for certain facts not specific to H\'{e}non maps.

Let $\Gamma^u$ and $\Gamma^s$ be two families of $C^1$ curves in
$\R^2$ such that
\begin{itemize}
\item the curves in $\Gamma^u$, respectively $\Gamma^s$, are
pairwise disjoint;

\item every $\gamma^u\in\Gamma^u$ meets every $\gamma^s\in\Gamma^s$
in exactly one point;

\item there is a minimum angle between $\gamma^u$ and $\gamma^s$ at
the point of intersection.
\end{itemize}
Then we define the lattice associated to $\Gamma^u$ and $\Gamma^s$
by
\[
\Lambda:=\{\gamma^u\cap\gamma^s:\,\gamma^u\in\Gamma^u,\,\gamma^s\in
\Gamma^s\}.\index{horseshoe}
\]\index{aaLambdag@$\Lambda$}
For $z\in\Lambda$ let $\gamma^u(z)$ and $\gamma^s(z)$ denote the
curves in $\Gamma^u$ and $\Gamma^s$ containing $z$, respectively.

We say that $\Xi$ is an {\em
$s$-sublattice}\index{s-sublattice@$s$-sublattice} (resp. {\em
$u$-sublattice}) \index{u-sublattice@$u$-sublattice}of $\Lambda$
if $\Lambda$ and $\Xi$ have a common defining family $\Gamma^u$
(resp. $\Gamma^s$) and the defining family $\Gamma^s$ (resp.
$\Gamma^u$) of $\Lambda$ contains that of~$\Xi$. A subset
$Q\subset\R^2$ \index{aaQg@$Q$} is said to be the {\em rectangle
spanned} \index{rectangle spanned} by $\Lambda$ if $\Lambda\subset
Q$ and $\partial Q$ is made up of two curves from $\Gamma^s$ and
two from $\Gamma^u$.

Next, we state Proposition A from \cite{BY00} which asserts the
existence of two lattices $\Lambda^+$ and $\Lambda^-$ with
essentially the same properties; for notation simplicity
statements about $\Lambda$ apply to both $\Lambda^+$ and
$\Lambda^-$.

\begin{proposition}
\label{prop:horseshoe-BY2000-propA} There are two lattices
$\Lambda^+$ and $\Lambda^-$ in $\R^2$ with the following
properties.
\begin{enumerate}

\item \label{item:propA-BY-top-struct}(Topological Structure)
$\Lambda$ is the disjoint union of $s$-sublattices $\Xi_i$,
$i=1,2\ldots$, where for each $i$, exists $R_i\in\N$ such that
$f^{R_i}(\Xi_i)$ is a $u$-sublattice of $\Lambda^+$ or
$\Lambda^-$.

\item \label{item:propA-BY-hyp-est} (Hyperbolic estimates)
    \begin{enumerate}

    \item Every $\gamma^u\in\Gamma^u$ is a $C^2(b)$ curve; and
     exists $\lambda_1>0$ such that for all $z\in\gamma^u\cap Q_i$,
    \[
    |Df^{R_i}(z)\tau|\geq\lambda_1^{R_i},
    \]
    where $\tau$ is the unit tangent
    vector to $\gamma^u$ at $z$ and $Q_i$ is the rectangle spanned
    by $\Xi_i$.

    \item For all $z\in\Lambda$, $\zeta\in \gamma^s(z)$ and $ j\geq
    1$
    we have
    \[
    |f^j(z)-f^j(\zeta)|<Cb^j,
    \]
    \end{enumerate}

\item (Measure estimate) $\mbox{Leb}(\Lambda\cap\gamma^u)>0$,
$\forall\gamma^u\in\Gamma^u$.

\item (Return time estimates) Let $R:\Lambda\rightarrow\N$ be
defined by $R(z)=R_i$ for $z\in\Xi_i$. Then there are $C_0>0$ and
$\theta_0<1$ such that on every $\gamma^u$
\[
\mbox{Leb}\{z\in\gamma^u:R(z)\geq n\}\leq C_0\theta_0^n,
\qquad\forall n\geq1.\index{aaCg0@$C_0$}
\index{aatheta0@$\theta_0$}
\]

\end{enumerate}
\end{proposition}

The proof of Proposition \ref{prop:horseshoe-BY2000-propA} can be
found in Sections 3 and 4 of \cite{BY00}. Since we will need to
prove the closeness of these horseshoes for nearby
Benedicks-Carleson parameters and this involves slight
modifications in the construction of the horseshoes itselves, we
will include, for the sake of completeness, the basic ideas of the
major steps leading to $\Lambda$.

 Consider the leaf of first
generation $W_1$ and the unique critical point
$z_0=(x_0,y_0)\in\C$ on it. Take the two outermost intervals of
the partition $\mathcal{P}_{[x_0]}$ as in Subsection
\ref{subsubsec:bounded-distortion} and denote them by $\Omega_0^+$
and $\Omega_0^-$; they support the construction of the lattices
$\Lambda^+$ and $\Lambda^-$, respectively. Again we use $\Omega_0$
to simplify notation and statements regarding to it apply to both
$\Omega_0^+$ and~$\Omega_0^-$.\index{aaOmegag0@$\Omega_0$}

Let $h:\Omega_0\rightarrow\R$ be a function whose graph is the
leaf of first generation $W_1$\index{aah@$h$}, when restricted to
the set $\Omega_0\times\R$ and $H:\Omega_0\rightarrow W_1$ be
given by $H(x)=(x,h(x))$.\index{aaHg@$H$}

\subsection{Leading Cantor sets}
\label{subsec:leading-Cantor-sets} The first step is to build the
Cantor set that constitutes the intersection of $\Lambda$ with the
leaf of first generation $W_1$. We build a sequence
$\Omega_0\supset\Omega_1\supset\Omega_2\ldots$ such that for every
$z\in H(\Omega_n)$,\index{aaOmegagn@$\Omega_n$}
$\mbox{dist}(f^j(z),\C)\geq\delta\e^{-\alpha j}$, for all
$j\in\{1,2,\ldots,n\}$. This is done by excluding from
$\Omega_{n-1}$ the points that at step $n$ fail to satisfy the
condition $\mbox{dist}(f^n(H(x)),\C)\geq\delta\e^{-\alpha n}$.
Then we define the Cantor set
$\Omega_\infty=\bigcap_{n\in\N}\Omega_n$.
\index{aaOmegaginf@$\Omega_\infty$} By the derivative estimate in
Subsection \ref{subsubsec:derivative-estimate}, on
$H(\Omega_\infty)$, the condition \eqref{eq:SAn} holds  and thus
$\left|Df^n(z)\tau(z)\right|>\e^{c_1n}$, for all
$n\in\N$.\index{leading Cantor sets}
\begin{remark}
We observe that there is a difference in the notation used in
\cite{BY00}: in here, the sets $\Omega_n$ (with
$n=0,1,\ldots,\infty$) are the vertical projections in the
$x$-axis of the corresponding sets in \cite{BY00}.
\end{remark}

\begin{remark}
\label{rem:non-uniqueness-omega-infty} We note that the procedure
leading to $\Omega_\infty$ is not unique. $\Omega_\infty$ is
obtained by successive exclusions of points from the set
$\Omega_0$. These exclusions are made according to the distance to
a suitable binding critical point every time we have a free return
to $[-\delta,\delta]\times\R$. Certainly, the choice for the
binding critical point in not unique which leads to different
exclusions. However, by the results referred in Subsection
\ref{subsubsec:tang-pos-dist-cs} all suitable binding points are
essentially the same and these possible differences in the
exclusions are insignificant in terms of the properties we want
$\Omega_\infty$ to have: slow approximation to the critical set
and expansion along the tangent direction to $W$.
\end{remark}

\subsection{Construction of long stable leaves}
\label{subsec:stable-curves} The next step towards building
$\Lambda$ involves the construction of long stable curves,
$\gamma^s(z)$, at every $z\in H(\Omega_\infty)$. This is done in
Lemma 2 of \cite{BY00}; let us review the inductive procedure used
there.

The contracting vector field of order 1, $e_1$, is defined
everywhere so we may consider the rectangle
$Q_0(\omega_0)=\cup_{z\in \omega_0}\gamma_1(z)$, where
$\gamma_1(z)$ denotes the $e_1$-integral curve segment $10b$ long
to each side of $z\in \omega_0$ and $\omega_0=H(\Omega_0)$. Let
also $Q^1_0(\omega_0)$ denote the $Cb$-neighborhood of
$Q_0(\omega_0)$ in $\R^2$. We observe that by
\eqref{item:cont-field-e1} of Section
\ref{subsec:contractive-vector-field} the $\gamma_1$ curves in
$Q_0(\omega_0)$ have slopes $\approx \pm 2a\delta$ depending on
whether $\Omega_0$ refers to $\Omega_0^+$ or $\Omega_0^-$.

Suppose that for every connected component $\omega\in
H(\Omega_{n-1})$ we have a strip foliated by integral curves of
$e_n$, $Q_{n-1}(\omega)=\cup_{z\in \omega}\gamma_n(z)$, where
$\gamma_n(z)$ \index{aagamman@$\gamma_n(z)$}denotes the
$e_n$-integral curve segment $10b$ long to each side of $z\in
\omega$. From \cite[Section 3.3]{BY00} one deduces that the vector
field $e_{n+1}$ is defined on a $3(Cb)^n$ neighborhood of each
curve $\gamma_n(z)$, if $z\in H(\Omega_n)$. Consider the $(Cb)^n$-
neighborhood of $Q_{n-1}(\omega)$ in $\R^2$, denoted by
$Q^1_{n-1}(\omega)$. If $\tilde{\omega}\subset\omega$ is a
connected component of $H(\Omega_n)$ then
$Q_{n}(\tilde{\omega})=\cup_{z\in \tilde{\omega}}\gamma_{n+1}(z)$
\index{aaQgn@$Q_n(\omega)$}is defined and
\begin{equation}
\label{eq:construction-gamma-infty} Q^1_{n}(\tilde{\omega})\subset
Q^1_{n-1}(\omega),
\end{equation}
where $Q^1_{n}(\tilde{\omega})$ is a $(Cb)^{n+1}$- neighborhood of
$Q_{n}(\tilde{\omega})$ in $\R^2$.\index{aaQgn1@$Q_n^1(\omega)$}

To fix notation, for some $\omega \subset H(\Omega_0)$ and $n\in
\N$, when defined, $Q_{n}(\omega)=\cup_{z\in
\omega}\gamma_{n+1}(z)$ denotes a rectangle foliated by integral
curves of $e_{n+1}$ passing through $z\in \omega$ and $10b$ long
to each side of $z$. Besides, $Q^1_{n}(\omega)$ is a $(Cb)^{n+1}$-
neighborhood of $Q_{n}(\omega)$ in $\R^2$.

To finish the construction of $\gamma^s(z)$, for each $z\in
H(\Omega_\infty)$, take the sequence of connected components
$\omega_i\subset H(\Omega_i)$ containing $z$. We have
$\{z\}=\cap_{i}\omega_i$. Let $z_n$ denote the right end point of
$\omega_{n-1}$. Then  $\gamma_n(z_n)$ converges in the $C^1$- norm
to a $C^1$-curve $\gamma^s(z)$\index{long stable
leaves}\index{aagammas0@$\gamma^s(z)$} with the properties stated
in Proposition \ref{prop:horseshoe-BY2000-propA}. The curve
$\gamma_n(z_n)$ acts as an approximate long stable leaf of order
$n$. Note that the choice of the right end point is quite
arbitrary; in fact any curve $\gamma_n(\zeta)$ with
$\zeta\in\omega_{n-1}$ suits as an approximate stable leaf of
order $n$.


\subsection{The families $\Gamma^u$ and $\Gamma^s$}
\label{subsec:def-Lambda} The final step in the construction of
$\Lambda$ is to specify the families $\Gamma^u$ and $\Gamma^s$.
Set
\[
\Gamma^s:=\{\gamma^s(z):\,z\in\Omega_\infty\},
\]\index{aaGammags@$\Gamma^s$}
where $\gamma^s(z)$ is obtained as described in Subsection
\ref{subsec:stable-curves}. Consider
$\tilde{\Gamma}^u:=\{\gamma\subset W:\, \gamma \mbox{ is a
$C^2(b)$ segment connecting $\partial^sQ_0$}\}$, where $Q_0$ is
the rectangle spanned by the family of curves $\Gamma^s$, i.e.,
$Q_0\supset\bigcup_{z\in H(\Omega_\infty)}\gamma^s(z)$ and
$\partial Q_0$ is made up from two curves of $\Gamma^s$.
\index{aaQg0@$Q_0$} Set
\[
\Gamma^u:=\{\gamma:\,\gamma\mbox{ is the pointwise limit of a
sequence in $\tilde{\Gamma}^u$}\}.
\]\index{aaGammagu@$\Gamma^u$}

\subsection{The $s$-sublattices and the return times}
\label{subsec:s-sublattice-return-times}Recall that we are
interested in two lattices $\Lambda^+$ and $\Lambda^-$.
Therefore, when we refer to return times we mean return times from
the set $\Lambda^+\cup\Lambda^-$ to itself; in particular, a point
in $\Lambda^+$ may return to $\Lambda^+$ or $\Lambda^-$. However,
in order to simplify we just write $\Lambda$.

We anticipate that the return time function
$R:\Lambda\rightarrow\N$ is constant in each
$\gamma^s\in\Gamma^s$, so $R$ needs only to be defined in
$\Lambda\cap H(\Omega_0)=H(\Omega_\infty)$. Moreover, since
$H:\Omega_0\rightarrow W_1$ is a bijection we may also look at $R$
as being defined on $\Omega_\infty$. We will build partitions on
subsets of $\Omega_0$ and use 1-dimensional language. For example,
$f^n(z)=\zeta$ for $z,\zeta\in H(\Omega_\infty)$ means that
$f^n(z)\in\gamma^s(\zeta)$; similarly, for subsegments
$\omega,\omega^*\subset H(\Omega_0)$, $f^n(\omega)=\omega^*$ means
that $f^n(\omega)\cap\Lambda$, when slid along $\gamma^s$ curves
back to $H(\Omega_0)$, gives exactly $\omega^*\cap\Lambda$. For an
interval $I\subset\Omega_{n-1}$ such that $f^n(H(I))$ intersects
the critical region, $\P|f^n(H(I))$ refers to $\P_{[\tilde{z}]}$
where $\tilde{z}\in\C$ is a suitable binding critical point for
all $f^n(H(I))$ whose existence is a consequence of Lemma 7 from
\cite{BY93}, mentioned in Subsection
\ref{subsubsec:tang-pos-dist-cs}.

We will construct sets $\tilde{\Omega}_n\subset\Omega_n$ and
partitions $\tilde{\P}_n$ of $\tilde{\Omega}_n$ so that
$\tilde{\Omega}_0\supset\tilde{\Omega}_1\supset\tilde{\Omega}_2\ldots$
and $z\in H(\tilde{\Omega}_{n-1}\setminus\tilde{\Omega}_n)$ if and
only if $R(z)=n$. Let $\hat{\P}$ \index{aaPghat@$\hat\P$}be the
partition of $H(\Omega_0\setminus\Omega_\infty)$ into connected
components. In what follows $\mathcal{A}\vee\mathcal{B}$ is the
join of the partitions $\mathcal{A}$ and $\mathcal{B}$, that is
$\mathcal{A}\vee\mathcal{B}=\{A\cap
B:\,A\in\mathcal{A},B\in\mathcal{B}\}$.

\begin{definition}
\label{def:regular-return} An interval $I\in\Omega_n$ is said to
make a {\em regular return} to $\Omega_0$ at time $n$ if
{\renewcommand{\theenumi}{\roman{enumi}}\index{regular return}
\begin{enumerate}

\item all of $f^n(H(I))$ is free;

\item\label{item:def-reg-ret-cond-2} $f^n(H(I))\supset 3H(\Omega_0)$.

\end{enumerate}
}

\begin{remark}The constant $3$ in the definition of regular return
 is quite arbitrary. In fact its purpose is to guarantee that
 $f^n(H(I))$ traverses $Q_0$ by wide margins. When $n$ is a regular
 return of a
 certain segment
 $I$ for a fixed parameter it may happen that $n$ does not
 classify as a regular return of a perturbed parameter even though
 the image of $I$ after $n$ iterates by the perturbed dynamics
 crosses $Q_0$ by wide margins. We overcome this detail
 simply by considering that if \eqref{item:def-reg-ret-cond-2}
 holds with $2$ instead of $3$ for any perturbed parameter
 then we consider $n$ as a regular return for the perturbed dynamics.
 Observe that no harm results from making this assumption since it is still
 guaranteed that $Q_0$ is traversed by wide margins.
 \end{remark}

\end{definition}

\subsubsection{Rules for defining $\tilde{\Omega}_n$, $\tilde{\P}_n$
and $R$}
\label{subsubsec:rules-def-R}\index{aaOmegagntilde@$\tilde\Omega_n$}
\index{aaPgntilde@$\tilde\P_n$}\index{aaRg@$R(z)$}
\begin{enumerate}
\setcounter{enumi}{-1}

\item $\tilde{\Omega}_0=\Omega_0$,
$\tilde{\P}_0=\{\tilde{\Omega}_0\}$.\newline Consider
$I\in\tilde{\P}_{n-1}$.

\item \label{item:rule1-def-R}If $I$ does not make a regular return
to $\Omega_0$ at time $n$, put $I\cap\Omega_n$ into
$\tilde{\Omega}_n$ and set
$\tilde{\P}_n|(I\cap\Omega_n)=H^{-1}\left((f^{-n}\P)|(H(I\cap\Omega_n))
\right)$.

\item \label{item:rule2-def-R}If $I$ makes a regular return at time $n$,
we put $\tilde{I}=H^{-1}\left(H(I)\setminus
f^{-n}(H(\Omega_\infty))\right)\cap\Omega_n$ in
$\tilde{\Omega}_n$, and let
$\tilde{\P}_n|\tilde{I}=H^{-1}\left((f^{-n}\P\vee
f^{-n}\hat{\P})|H(\tilde{I}) \right)$. For $z\in H(I)$ such that
$f^n(z)\in H(\Omega_\infty)$, we define $R(z)=n$.\index{return to
the horseshoe}


\item For $z\in H(\cap_{n\in\N_0}\tilde{\Omega}_n)$, set
$R(z)=\infty$.

\end{enumerate}

\subsubsection{Definition of the $s$-sublattices}
\label{subsubsec:def-s-sublattices} Each $\Xi_i$
\index{aaXigi@$\Xi_i$}in
Proposition~\ref{prop:horseshoe-BY2000-propA} is a sublattice
corresponding to a subset of $\Lambda\cap W_1$ of the form
$f^{-n}(H(\Omega_\infty))\cap\Lambda\cap H(I)$, where
$I\in\tilde{\P}_{n-1}$ makes a regular return at time $n$. We will
use the notation
$\Upsilon_{n,j}=H^{-1}\left(f^{-n}(H(\Omega_\infty))\cap\Lambda
\cap H(I)\right)$. \index{aaUpsilonn1@$\Upsilon_{n,j}$}
\index{s-sublattice@$s$-sublattice|textbf}Note that
$R(H(\Upsilon_{n,j}))=n$ and $\Upsilon_{n,j}$ determines
univocally the corresponding $s$-sublattice. For this reason we
allow some imprecision by referring ourselves to $\Upsilon_{n,j}$
as an $s$-sublattice.

In order to prove the assertions \eqref{item:propA-BY-top-struct}
and \eqref{item:propA-BY-hyp-est} of Proposition
\ref{prop:horseshoe-BY2000-propA} one needs to verify that
$f^{R_i}(\Xi_i)$ is an $u$-sublattice which requires to
demonstrate that $f^{R_i}(\Xi_i)$ matches completely with
$\Lambda$ in the horizontal direction. If $\Xi_i$ corresponds to
some $\Upsilon_{n,j}$, then the matching of the Cantor sets will
follow from the inclusion
\begin{equation}
\label{eq:match-Cantor-sets} f^n(H(I\cap\Omega_\infty))\supset
H(\Omega_\infty).
\end{equation}
It is obvious that $H(\Omega_\infty)\subset f^n(H(I))$ by
definition of regular return. Nevertheless,
\eqref{eq:match-Cantor-sets} is saying that if $z\in H(I)$ and
$f^n(z)$ hits $H(\Omega_\infty)$, after sliding along a $\gamma^s$
curve, then $z\in H(I)\cap H(\Omega_\infty)$. This is proved in
Lemma 3 of \cite{BY00}. In particular, we may write
$\Upsilon_{n,j}=H^{-1}\left(f^{-n}(H(\Omega_\infty))\cap
H(I)\right)$.\index{aaUpsilonn1@$\Upsilon_{n,j}$}

\subsection{Reduction to an expanding map}
\label{subsec:quotient-space} The H\'{e}non maps considered here
are perturbations of the map $f_{2,0}(x,y)=(1-2x^2,0)$ whose
action is horizontal. Also, as we have seen, the horizontal
direction is typically expanding. This motivates considering the
quotient space $\bar{\Lambda}$
\index{aaLambdagbar@$\bar\Lambda$}\index{quotient horseshoe}
obtained by collapsing the stable curves of $\Lambda$; i.e.
$\bar{\Lambda}=\Lambda/\sim$, where $z\sim z'$ if and only if
$z'\in\gamma^s(z)$. We define the natural projection
$\bar\pi:\Lambda\rightarrow\bar\Lambda$ given by
$\bar\pi(z)=\gamma^s(z)$\index{aapibar@$\bar\pi$}. As implied by
assertion \eqref{item:propA-BY-top-struct} of Proposition
\ref{prop:horseshoe-BY2000-propA}, $f^R:\Lambda\rightarrow\Lambda$
takes $\gamma^s$ leaves to $\gamma^s$ leaves (see Lemma 2 of
\cite{BY00} for a proof). Thus, we may define the quotient map
$\overline{f^R}:\bar{\Lambda}\rightarrow\bar{\Lambda}$. Observe
that each $\bar{\Xi}_i$\index{aaXigibar@$\bar\Xi_i$} is sent by
$\overline{f^R}$ homeomorphically onto $\bar{\Lambda}$. Besides we
may  define a \emph{reference measure} \index{reference measure}
$\bar{m}$ on $\bar{\Lambda}$, whose representative on each
$\gamma^u\in\Gamma^u$ is a finite measure equivalent to the
restriction of the 1-dimensional Lebesgue measure on
$\gamma^u\cap\Lambda$ and denoted by
$m_{\gamma^u}$.\index{aambar@$\bar
m$}\index{aamammau@$m_{\gamma^u}$}

One can look at $\overline{f^R}$ as an expanding Markov map (see
Proposition B of \cite{BY00} for precise statements and proofs).
Moreover, the corresponding transfer operator, relative to the
reference measure $\bar{m}$, has a spectral gap (see Section 3 of
\cite{Yo98}, specially Proposition A). It follows that
$\overline{f^R}$ has an absolutely continuous invariant measure
given by $\bar{\nu}=\bar{\rho}d\bar{m}$\index{aanubar@$\bar\nu$},
with $M^{-1}\leq \bar{\rho}\leq M$
\index{aarho1bar1@$\bar\rho$}for some $M>0$ \index{aaMg@$M$}(see
\cite[Lemma 2]{Yo98}).

\section{Proximity of critical points}
\label{sec:Proximity-Critical-Points}

In this section we show that up to a fixed generation we have
closeness of the critical points for nearby Benedicks-Carleson
parameters. This is the content of Proposition
\ref{prop:prox-crit-points} which summarizes this section. Its
proof involves a finite step induction scheme on the generation
level. We prepare it by proving first the closeness of critical
points of generation~1 in Lemma \ref{lem:prox-cp-gen-1}.
Afterwards, in Lemma \ref{lem:prox-cp-gen>27} we obtain the
closeness of critical points of higher generations using the
information available for lower ones.

Recall that since $f_{a,b}$ is $C^\infty$, then the unstable
manifold theorem ensures that $W$ is $C^r$  for any $r>0$.
Moreover, $W$ varies continuously in the $C^r$ topology with the
parameters in compact parts. As we are only considering parameters
in $\B$,  for each of these dynamics there is a unique critical
point $\hat{z}$ of generation~1 situated on the roughly horizontal
segment of $W$ containing the fixed point $z^*$.

\begin{lemma}
\label{lem:prox-cp-gen-1} Let $(a,b)\in\B$, $\varepsilon>0$ be
given and $\hat{z}$ be the critical point of generation 1 of
$f_{a,b}$. There exists a neighborhood $\U$ of $(a,b)$ such that,
if
 $(a',b')\in\U\cap\B$ and $\hat{z}'$ denotes the
critical point of $f_{a',b'}$ of generation $1$, then
$|\hat{z}-\hat{z}'|<\varepsilon$. Moreover, if $\tau(\hat{z})$ and
$\tau(\hat{z}')$ are the unit vectors tangent to $W$ and $W'$ at
$\hat{z}$ and $\hat{z}'$ respectively, then
$|\tau(\hat{z})-\tau(\hat{z'})|<\varepsilon$.
\end{lemma}

\begin{proof}
Consider the disk $\gamma=W_1\cap [-10b,10b]\times\R$. There is a
neighborhood $\U$ of $(a,b)$ such that for every $(a',b')\in \U$
there exists a disk $\gamma'\subset W'$ which is
$\varepsilon^2$-close to $\gamma$ in the $C^r$ topology. It is clear
that both $\gamma$ and $\gamma'$ are $C^2(b)$ curves and there are
$\hat{z}\in \gamma$ and $\hat{z}'\in \gamma'$ critical points of
$f_{a,b}$ and $f_{a',b'}$ respectively. Our goal is to show that
$|\hat{z}-\hat{z}'|<\varepsilon$. The strategy is to pick an
approximate critical point $\hat{z}_M$ of $\hat{z}$ and then prove
the existence of an approximate critical point $\hat{z}'_M$ of
$\hat{z}'$ sufficiently close to $\hat{z}_M$ in order to conclude
that, if we choose $M$ large enough, we get the desired closeness of
$\hat{z}$ and $\hat{z}'$ (see Figure~\ref{fig:relative-position-1}).
\begin{figure}[h]
\includegraphics[scale=1.5]{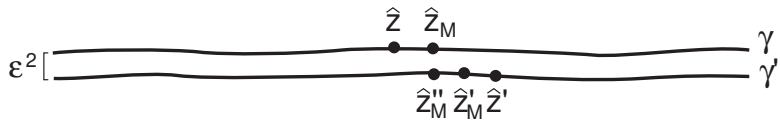}
\caption{Possible configuration of the critical points and their
approximates} \label{fig:relative-position-1}
\end{figure}
Take $M\in \N$ so that $b^M<\varepsilon^2\le b^{M-1}$. Let
$\hat{z}_M\in\gamma$ be such that
$e_M(\hat{z}_M)=\tau(\hat{z}_M)$. Note that
$|\hat{z}-\hat{z}_M|<Cb^M$. Let $\hat{z}_M''\in \gamma'$ be such
that $|\hat{z}_M-\hat{z}_M''|<\varepsilon^2$ and
$\left|\tau(\hat{z}_M)-\tau(\hat{z}_M'')\right|<\varepsilon^2$.
Now, $\hat{z}_M''$ may not be the approximate critical point
$\hat{z}_M'$ we are looking for, but we will show that it is very
close to $\hat{z}_M'$. In fact, we assert that the angle between
$e'_M(\hat{z}_M'')$ and $\tau(\hat{z}_M'')$ is of order
$\varepsilon^2$, which allows us to find a nearby $\hat{z}_M'$ as
a solution of $e'_M(z')=\tau(z')$, which ultimately is very close
to the critical point $\hat{z}'$.

Before we prove this last assertion we must guarantee that the
vector field $e_M'$ is defined in a neighborhood of $\hat{z}_M''$
and for that we must have some expansion. Since $\hat{z}$ is a
critical point of $f_{a,b}$, then
$\left|Df_{a,b}^M(\hat{z})\binom{0}{1}\right|>\e^{cM}$. If
necessary we tighten $\U$ so that for every $z$ in a compact set
of $\R^2$,
$\left|Df_{a,b}^M(z)\binom{0}{1}-Df_{a',b'}^M(z)\binom{0}{1}\right|$
is small enough for having
$\left|Df_{a',b'}^M(\hat{z})\binom{0}{1}\right|>\e^{cM/2}$, which
implies that $e'_M$ is well defined in a ball of radius
$3Cb^{M-1}>3C\varepsilon^2$ around $\hat{z}$. Note that $b\ll
\lambda$ and the Matrix Perturbation Lemma applies.

We take $\U$ sufficiently small so that
$\left|e'_M(\hat{z}_M'')-e_M(\hat{z}_M'')\right|<\varepsilon^2$.
This is possible because  $e'_M(z)$ and $e_M(z)$ are the maximally
contracted vectors of $Df_{a',b'}^M(z)$ and $Df_{a,b}^M(z)$,
respectively. Thus it is only a matter of making $Df_{a',b'}^M(z)$
very close to $Df_{a,b}^M(z)$, for every $z$ in a compact set.
Hence
\begin{equation*}
\begin{split}
\left|e'_M(\hat{z}_M'')-\tau(\hat{z}_M'')\right|&<\left|e'_M(\hat{z}_M'')-e_M(\hat{z}_M'')\right|+
\left|e_M(\hat{z}_M'')-e_M(\hat{z}_M)\right|+\left|e_M(\hat{z}_M)-\tau(\hat{z}_M)\right|\\
&\quad + \left|\tau(\hat{z}_M)-\tau(\hat{z}_M'')\right| \\
&<\varepsilon^2+ C|\hat{z}_M-\hat{z}_M''|+0+\varepsilon^2\\
&<C\varepsilon^2
\end{split}
\end{equation*}
Writing  $z=(x,y)$ and taking into account that $\gamma'$ is
nearly horizontal we may think of it as the graph of $\gamma'(x)$.
Let us also ease on the notation so that $\tau(x)$ and $e'_M(x)$
denote the slopes of the respective vectors at $z=\gamma'(x)$. We
know that $\left|{d\tau}/{dx}\right|<10b$,
$\left|{de'_M}/{dx}\right|=2a+\mathcal{O}(b)$ and
$\left|{d^2e'_M}/{dx^2}\right|<C$.
\begin{figure}[h]
\includegraphics[scale=0.9]{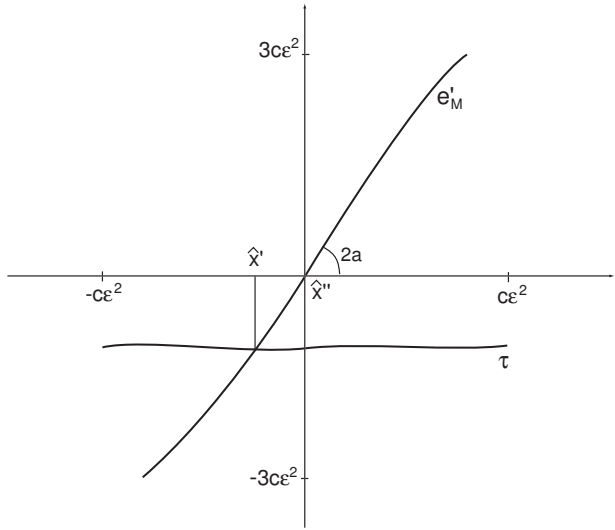}
\caption{Solution of $e'_M(z)=\tau(z)$} \label{fig:approx-cp}
\end{figure}
 As a consequence we obtain  $\hat{z}_M'$ such that $e'_M(\hat{z}_M')=\tau(\hat{z}_M)$ and
$|\hat{z}_M'-\hat{z}_M''|<{C}\varepsilon^2/3$ (see
Figure~\ref{fig:approx-cp}). Now since there is a unique critical
point $\hat{z}'$ in $\gamma'$ we must have
$|\hat{z}'-\hat{z}_M'|<C\varepsilon^2$, which yields
\[
|\hat{z}-\hat{z}'|\leq
|\hat{z}-\hat{z}_M|+|\hat{z}_M-\hat{z}_M''|+|\hat{z}_M''-\hat{z}_M'|+
|\hat{z}_M'-\hat{z}|<C\varepsilon^2<\varepsilon,
\]
as long as $\varepsilon$ is sufficiently small.

 Concerning the
inequality $|\tau(\hat{z})-\tau(\hat{z}')|<\varepsilon$, simply
observe that since $\gamma$ and $\gamma'$ are $C^2(b)$ curves we
have
\begin{equation*}
\begin{split}
|\tau(\hat{z})-\tau(\hat{z}')|&<|\tau(\hat{z})-
\tau(\hat{z}_M)|+|\tau(\hat{z}_M)-\tau(\hat{z}_M'')|
+|\tau(\hat{z}_M'')-\tau(\hat{z}_M')|+|\tau(\hat{z}_M')-\tau(\hat{z}')| \\
&<10b|\hat{z}-\hat{z}_M|+\varepsilon^2+10b|\hat{z}_M''-
\hat{z}_M'|+10b|\hat{z}_M'-\hat{z}'|\\
&<\varepsilon.
\end{split}
\end{equation*}
\end{proof}

As a consequence of Lemma \ref{lem:prox-cp-gen-1} we have that for
a sufficiently small  $\U$  we manage to make $W_1'$ (the leaf of
$W'$ of generation 1) to be as close to $W_1$ (the leaf of $W$ of
generation~1) as we want. This is important because the leaves of
higher generations are defined by successive iterations of the
first generation leaf. We also remark that by the rules of
construction of the critical set we may use the argument  of Lemma
\ref{lem:prox-cp-gen-1} to obtain proximity of the critical points
up to generation $27$. For higher generations we need the
following lemma.

\begin{lemma}
\label{lem:prox-cp-gen>27} Let $N\in\N$,  $(a,b)\in\B$ and
 $\varepsilon>0$ be given. Assume there exists a
neighborhood $\U$ of $(a,b)$ such that for each
$(a',b')\in\U\cap\B$ and any critical point $\hat{z}$ of $f_{a,b}$
of generation $g<N$, there is a critical point $\hat{z}'$ of
$f_{a',b'}$ of the same generation with
$|\hat{z}-\hat{z}'|<\varepsilon$.  If a critical point $\hat{z}$
of $f_{a,b}$ is created at step $g+1$, then we may tighten $\U$ so
that a critical point $\hat{z}'$ of generation $g+1$ is created
for $f_{a',b'}$ and $|\hat{z}-\hat{z}'|<\varepsilon$. Moreover, if
$\tau(\hat{z})$ and $\tau(\hat{z}')$ are the unit vectors tangent
to $W$ and $W'$ at $\hat{z}$ and $\hat{z}'$ respectively, then
$|\tau(\hat{z})-\tau(\hat{z'})|<\varepsilon$.
\end{lemma}

\begin{proof}
As we are only interested in arbitrarily small $\varepsilon$, we may
assume that $\varepsilon<b^N$. Suppose that a critical point
$\hat{z}$ of generation $g+1$ is created for $f_{a,b}$. Then, by the
rules of construction of critical points, there are $z=(x,y)$ lying
in a $C^2(b)$ segment $\gamma\subset W$ of generation $g+1$ with
$\gamma$ extending beyond $2\varrho^{g+1}$ to each side of $z$ and a
critical point $\tilde{z}=(x,\tilde{y})$ of generation not greater
than $g$ such that $|z-\tilde{z}|<b^{{(g+1)}/{540}}$. Moreover,
$|\hat{z}-z|<|z-\tilde{z}|^{{1}/{2}}$.

Taking  $\gamma$ as a compact disk of $W$, there is a neighborhood
$\U$ of $(a,b)$ such that for every $(a',b')\in \U$ we can find  a
disk $\gamma'\subset W'$ of generation $g+1$ which is
$\varepsilon^2$-close to $\gamma$ in the $C^r$-topology. It is
clear that $\gamma'$ is a $C^2(b)$ curve. Our aim is to show that
a critical point $\hat{z}'$ of $f_{a',b'}$ and generation $g+1$ is
created in the segment $\gamma'$ with
$|\hat{z}-\hat{z}'|<\varepsilon$.

By the inductive hypothesis there is
$\tilde{z}'=(\tilde{x}',\tilde{y}')$ a critical point of
$f_{a',b'}$ such that $|\tilde{z}-\tilde{z}'|<\varepsilon$. Let
$z'=(\tilde{x}',y')$ belonging to $\gamma'$. Since $\gamma'$ is
$\varepsilon^2$-close to $\gamma$ in the $C^r$ topology and
$\varepsilon<b^{N}$, which is completely insignificant when
compared to $\varrho^{g+1}<\varrho^{N}$ (recall that $\varrho \gg
b$), we may assume that $\gamma'$ extends more than
$2\varrho^{g+1}$ to both sides of $z'$. Moreover, letting
$\zeta'=(x,\eta')\in \gamma'$   we have
\begin{equation*}
\begin{split}
|\tilde{z}'-z'|&<|\tilde{z}'-\tilde{z}|+|\tilde{z}-z|+|z-\zeta'|
+|\zeta'-z'| \\
&<\varepsilon+ b^{\tfrac{g+1}{540}}+2\varepsilon^2+2\varepsilon\\
&\lesssim b^{\tfrac{g+1}{540}},
\end{split}
\end{equation*}
where we used the fact that  $\varepsilon<b^{N}\ll
b^{\tfrac{N}{540}}<b^{\tfrac{g+1}{540}}$ (see
Figure~\ref{fig:relative-position-2}).
\begin{figure}[h]
\includegraphics[scale=1.2]{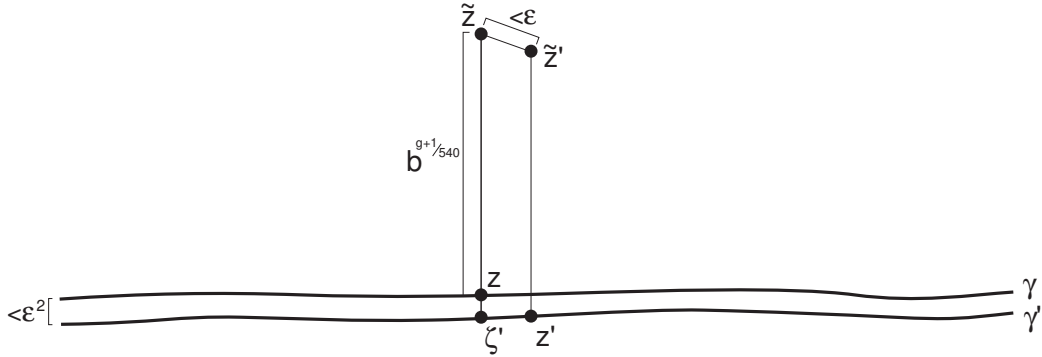}
\caption{Possible relative position of the critical points}
\label{fig:relative-position-2}
\end{figure}
By the rules of construction of critical points, a unique critical
point $\hat{z}'$ of generation $g+1$ is created in the segment
$\gamma'$. We are left to show that
$|\hat{z}-\hat{z}'|<\varepsilon$. For that we repeat the argument
in the proof of Lemma \ref{lem:prox-cp-gen-1}.
\end{proof}

\begin{corollary}
\label{prop:prox-crit-points} Let $N\in\N$, $(a,b)\in\B$ and
$\varepsilon>0$ be given. There is a neighborhood $\U$ of $(a,b)$
such that if $(a',b')\in\U\cap\B$ then, for any critical point
$\hat{z}$ of $f_{a,b}$ of generation smaller than $N$, there is a
critical point $\hat{z}'$ of $f_{a',b'}$ of the same generation
such that $|\hat{z}-\hat{z}'|<\varepsilon$. Moreover if
$\tau(\hat{z})$ and $\tau(\hat{z}')$ are the unit vectors tangent
to $W$ and $W'$ at $\hat{z}$ and $\hat{z}'$ respectively, then
$|\tau(\hat{z})-\tau(\hat{z'})|<\varepsilon$.
\end{corollary}

\begin{proof}
The proof is just a matter of collecting the information in the
Lemmas \ref{lem:prox-cp-gen-1} and \ref{lem:prox-cp-gen>27} and
organize it in a finite step induction scheme.

\begin{enumerate}
\item  First obtain the proximity
of the critical points of generation $1$, which has already been
done in Lemma \ref{lem:prox-cp-gen-1}.

\item Then realize that the same argument in the proof of Lemma
\ref{lem:prox-cp-gen-1} also gives the proximity of the $2^{26}$
critical points of generation smaller than $27$. (See the rules of
construction of critical points in Subsection
\ref{subsubsec:rules-critical-set}).

\item Apply the inductive step stated in Lemma \ref{lem:prox-cp-gen>27}
to obtain the proximity of critical points of higher and higher
generation.

\item Stop the process when the proximity of all
critical points of generation smaller than $N$ is achieved.
\end{enumerate}

Naturally every time we apply Lemma \ref{lem:prox-cp-gen>27} to
increase the generation level for which the conclusion of the
proposition holds, we may need to decrease the size of the
neighborhood $\U$. However, because the number of critical points
of a given generation is finite and the statement of the
proposition is up to generation $N$, at the end we still obtain a
neighborhood containing a non-degenerate ball around $(a,b)$ where
the proposition holds.
\end{proof}

\section{Proximity of leading Cantor sets}
\label{sec:Prox-Cantor-sets}

\begin{figure}[h]
\includegraphics[scale=1]{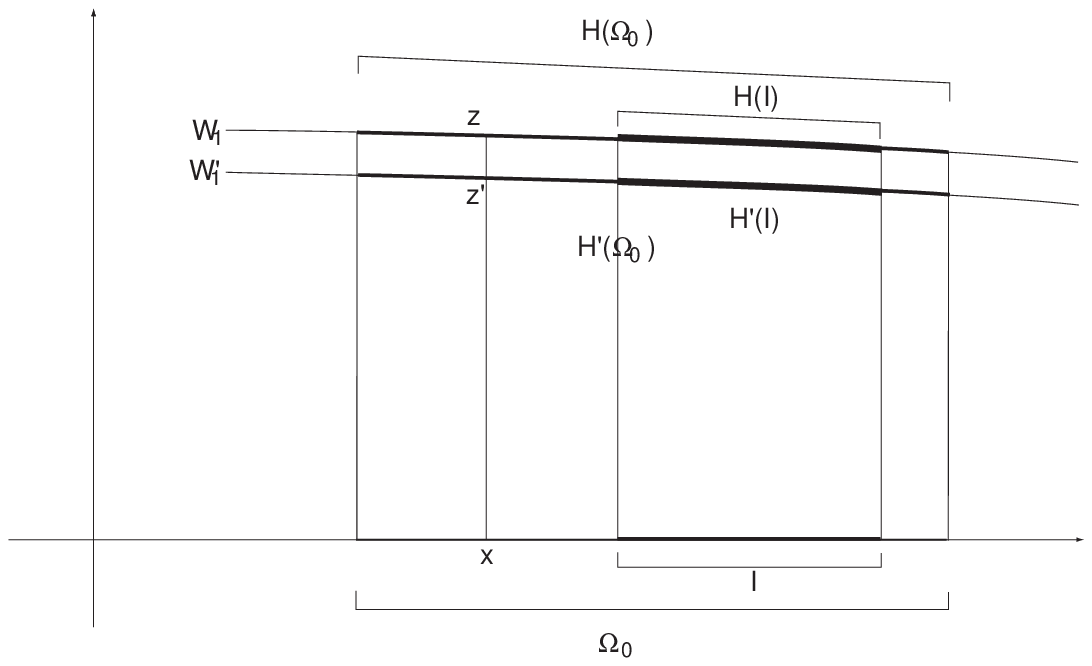}
\caption{} \label{fig:notation1}
\end{figure}

 Attending to Lemma
\ref{lem:prox-cp-gen-1}, we may assume  that
$\Omega_0=\Omega_0'$\index{aaOmegag0'@$\Omega_0'$}. Let
$h,h':\Omega_0\rightarrow\R$ \index{aah'@$h'$}be functions whose
graphs are the leaves of first generation $W_1$ and
$W_1'$\index{aaWg1'@$W_1'$} respectively, when restricted to the
set $\Omega_0\times\R$. Given an interval $I\subset\Omega_0$ the
segments $\omega=H(I)$ and $\omega'=H'(I)$ are respectively the
subsets of $W_1$ and $W_1'$ which correspond to the images in the
graph of $h$ and $h'$ of the interval $I$. Accordingly, if
$x\in\Omega_0$ then $z=H(x)=(x,h(x))$ and $z'=H'(x)=(x,h'(x))$.
See Figure~\ref{fig:notation1}.\index{aaHg'@$H'$}

Our goal in this section  is to show the proximity of the Cantor
sets $\Omega_\infty$ for close Benedicks-Carleson parameters. More
precisely, given any $\varepsilon>0$ we will exhibit a
neighborhood $\U$ of $(a,b)$ such that
$\left|\Omega_\infty\bigtriangleup\Omega_\infty'\right|<\varepsilon$
for all $(a',b')\in\U\cap\B$, where $\bigtriangleup$ represents
symmetric difference between two sets. In the process, we make a
modification in the first steps of the procedure described in
Subsection \ref{subsec:leading-Cantor-sets} to build
$\Omega_\infty'$\index{aaOmegaginf'@$\Omega_\infty'$}, which
carries only minor differences with respect to the set we would
obtain if we were to follow the rules strictly. Ultimately, this
affects the construction of the horseshoes $\Lambda'$. However,
the horseshoes are not uniquely determined and we will evince that
the modifications introduced leave unchanged the properties that
they are supposed to have.

\begin{lemma}
\label{lem:choice-N1} Given $\varepsilon>0$, there exists
$N_1\in\N$ such that
$|\Omega_n\setminus\Omega_\infty|<\varepsilon$ for every
$(a,b)\in\B$  and $n\geq N_1$.
\end{lemma}

\begin{proof}
This is a consequence of \cite[Lemma 4]{BY00} where it is proved
that
\begin{equation}
\label{eq:estimate-Omegan-Omegan-1}
\frac{|\Omega_{n-1}\setminus\Omega_{n}|}{|\Omega_{n-1}|}\leq
C_1\delta^{1-3\beta}\e^{-\alpha(1-3\beta)n}.
\end{equation}
This inequality follows from the fact that  any connected
component $\omega\in H(\Omega_{n-1})$ grows to reach a length
$|f^n(\omega)|\geq\delta^{3\beta}\e^{-3\alpha\beta n}$, while the
subsegment of $f^n(\omega)$ to be deleted in the construction of
$\Omega_n$ has length at most $4\delta\e^{-\alpha n}$; then,
simply take bounded distortion into consideration.

From \eqref{eq:estimate-Omegan-Omegan-1} one easily gets
\begin{equation*}
\begin{split}
|\Omega_{n}\setminus\Omega_\infty|&=\sum_{j=0}^{+\infty}|\Omega_{n+j}\setminus\Omega_{n+j+1}|\\
&\leq C_1\delta^{1-3\beta}\sum_{j=1}^{+\infty}\e^{-\alpha(1-3\beta)(n+j)}|\Omega_{n+j-1}|\\
&\leq
C_1\delta^{1-3\beta}|\Omega_n|\sum_{j=1}^{+\infty}\e^{-\alpha(1-3\beta)(n+j)}\\
&\leq
C_1\delta^{1-3\beta}\frac{\e^{-\alpha(1-3\beta)(n+1)}}{1-\e^{-\alpha(1-3\beta)}}
\end{split}
\end{equation*}
Hence, choose $N_1$ sufficiently large so that
\[
C_1\delta^{1-3\beta}\frac{\e^{-\alpha(1-3\beta)(N_1+1)}}
{1-\e^{-\alpha(1-3\beta)}}<\varepsilon.
\]
\end{proof}

Observe that, as a consequence of the unstable manifold theorem,
 for every $\varepsilon>0$ and $n\in\N$, there exists a neighborhood
$\U$ of $(a,b)$ such that for every $(a',b')\in\U$ we have
\begin{equation}
\label{eq:CR-omega0-proximity}
\max\left\{\left\|H-H'\right\|_r,\left\|f_{a,b}\circ
H-f_{a',b'}\circ H'\right\|_r,\ldots,\left\|f_{a,b}^{n}\circ
H-f_{a',b'}^{n}\circ H'\right\|_r \right\}<\varepsilon,
\end{equation}
where $r\geq2$ and $\|\cdot\|_r$ is the $C^r$-norm in $\Omega_0$.
In what follows $\Omega_\infty=\cap_{n\in\N}\Omega_n$ is built as
described in Section \ref{subsec:leading-Cantor-sets} for
$f=f_{a,b}$.

\begin{lemma}
\label{lem:suitable-cp} Let $n\in\N$ and $(a,b)\in\B$ be given and
$I$ be a connected component of $\Omega_{n-1}$. Suppose
$f_{a,b}^n(H(I))$ intersects $(-\delta,\delta)\times\R$. There is
a neighborhood $\U$ of $(a,b)$ such that for every
$(a',b')\in\U\cap\B$ and $x\in I\cap\Omega_n$, if
$f_{a,b}^n(H(x))\in(-\delta,\delta)\times\R$ and $\hat{z}$ is a
suitable binding critical point, then there exists a binding
critical point $\hat{z}'$ of $f_{a',b'}$ close to $\hat{z}$
suitable for $f_{a',b'}^n(H'(x))$ and
$|f_{a',b'}^n(H'(x))-\hat{z}'|\gtrsim \delta\e^{-\alpha n}$.
\end{lemma}

\begin{proof}
Let $\tilde{I}=I\cap\Omega_n$ and
 $\U$ be a neighborhood of $(a,b)$ such that Corollary~\ref{prop:prox-crit-points} applies up to $n$ with $ b^{2n}$
in the place of $\varepsilon$ and equation
\eqref{eq:CR-omega0-proximity} also holds with $b^{4n}$ in the
place of $\varepsilon$. Then there is a critical point $\hat{z}'$
of $f_{a',b'}$ such that $|\hat{z}-\hat{z}'|< b^{2n}$ and
$\left\|f_{a,b}^n\circ H|_{\tilde{I}}-f_{a',b'}^n\circ
H'|_{\tilde{I}}\right\|_r<  b^{4n}$. We only need to prove that
this $\hat{z}'$ is a suitable binding point for $f_{a',b'}^n(z')$
and that $|f_{a',b'}^n(z')-\hat{z}'|\gtrsim\delta\e^{-\alpha n}$.
In order to verify the suitability of $\hat{z}'$ we have to check
that
\begin{enumerate}
    \item $f_{a',b'}^n(z')$ is in tangential position with
respect to $\hat{z}'$;
    \item $Df_{a',b'}^n(z')\tau(z')$ splits
correctly with respect to the contracting field around~$\hat{z}'$.
\end{enumerate}
The strategy  is to show that
$|f_{a,b}^n(z)-\hat{z}|=|f_{a',b'}^n(z')-\hat{z}'|+\mathcal{O}(
b^{2n})$. Then, because $f^n(z)$ is in tangential position with
respect to $\hat{z}$ and $ b^{2n}\ll\delta\e^{-\alpha
n}\leq|f^n(z)-\hat{z}|$, we conclude the tangential position for
$f_{a',b'}^n(z')$ with respect to $\hat{z}'$. As to the correct
splitting, we know that
$\left|Df^n(z)\tau(z)-(Df_{a',b'})^n(z')\tau(z')\right|< b^{4n}$
and $Df^n(z)\tau(z)$  makes an angle with the relevant contracting
field of approximately $(2a\pm1)|f^n(z)-\hat{z}|$. Finally, since
$|f^n(z)-\hat{z}|=|f_{a',b'}^n(z')-\hat{z}'|+\mathcal{O}( b^{2n})$
and $ b^{2n}\ll (2a\pm1)|f^n(z)-\hat{z}|$ we obtain the desired
result.

Let us start by proving (1). Observe that
\begin{equation*}
\begin{split}
|f_{a,b}^n(z)-\hat{z}|&\leq|f_{a,b}^n(z)-f_{a',b'}^n(z')|+
|f_{a',b'}^n(z')-\hat{z}'|+|\hat{z}-\hat{z}'|\\
&\leq \left\|f_{a,b}^n\circ H|_{\tilde{I}}-f_{a',b'}^n\circ
H'|_{\tilde{I}}\right\|_r+|f_{a',b'}^n(z')-\hat{z}'|+ b^{2n}\\
&\leq |f_{a',b'}^n(z')-\hat{z}'|+2 b^{2n}.
\end{split}
\end{equation*}
Interchanging $z$ with $z'$ we easily get
$|f_{a',b'}^n(z')-\hat{z}'|\leq|f_{a,b}^n(z)-\hat{z}|+2 b^{2n}$
which allows us to write
$|f_{a',b'}^n(z')-\hat{z}'|=|f_{a,b}^n(z)-\hat{z}|+\mathcal{O}(
b^{2n})$. Consider now $s$ and $s'$  the lines through $\hat{z}$
and $\hat{z}'$ with slopes $\tau(\hat{z})$ and $\tau(\hat{z}')$
respectively. By Corollary~\ref{prop:prox-crit-points} we have
$|\hat{z}-\hat{z}'|< b^{2n}$ and also
$|\tau(\hat{z})-\tau(\hat{z}')|< b^{2n}$. Thus, when restricted to
the set $[-1,1]\times\R$ we have $\|s-s'\|_r<\mathcal{O}(
b^{2n})$. Let $\mbox{dist}(z,s)$ denote the distance from the
point $z$ to the segment $s\cap[-1,1]\times\R$. Then
\begin{equation*}
\begin{split}
\mbox{dist}(f_{a,b}^n(z),s)&\leq|f_{a,b}^n(z)-f_{a',b'}^n(z')|+\mbox{dist}
(f_{a',b'}^n(z'),s)\\
&\leq |f_{a,b}^n(z)-f_{a',b'}^n(z')|+\|s-s'\|_r+\mbox{dist}(f_{a',b'}^n(z'),s')\\
&\leq \mbox{dist}(f_{a',b'}^n(z'),s')+\mathcal{O}( b^{2n})
\end{split}
\end{equation*}
Similarly we get $\mbox{dist}(f_{a',b'}^n(z'),s')\leq
\mbox{dist}(f_{a,b}^n(z),s)+\mathcal{O}( b^{2n})$, and so
$$\mbox{dist}(f_{a',b'}^n(z'),s')=
\mbox{dist}(f_{a,b}^n(z),s)+\mathcal{O}( b^{2n}).$$ Now, since
$f^n(z)$ is in tangential position with respect to $\hat{z}$, then
$$
\mbox{dist}(f_{a,b}^n(z),s)< c|f_{a,b}^n(z)-\hat{z}|^2,
$$
where $c\ll 2a$. Besides,
$|f_{a',b'}^n(z')-\hat{z}'|^2=\left(|f_{a,b}^n(z)-\hat{z}|+\mathcal{O}(
b^{2n})\right)^2= |f_{a,b}^n(z)-\hat{z}|^2+\mathcal{O}( b^{2n})$
because $ b^{2n}\ll\delta\e^{-\alpha n}\le|f_{a,b}^n(z)-\hat{z}|$.
Consequently
\[
\mbox{dist}(f_{a',b'}^n(z'),s')<
c|f_{a',b'}^n(z')-\hat{z}'|^2+\mathcal{O}( b^{2n}),
\]
which again by the insignificance of $ b^{2n}$ relative to
$|f_{a',b'}^n(z')-\hat{z}'|$ implies that $f_{a',b'}^n(z')$ is in
tangential position with respect to $\hat{z}'$.

Concerning (2), notice that if $(a',b')$ is sufficiently close to
$(a,b)$, then
\[
\left|Df_{a,b}^n(z)\tau(z)-Df_{a',b'}^n(z')\tau(z')\right|
\leq\left\|f_{a,b}^n\circ H|_{\tilde{I}}-f_{a',b'}^n\circ
H'|_{\tilde{I}}\right\|_r<  b^{4n}.
\]
Let $l$ and $l'$ denote the lengths of the fold periods for $z$
and $z'$. Take $m$ and $m'$ such that $(5b)^m\leq |z-\hat{z}|\leq
(5b)^{m-1}$ and $(5b)^{m'}\leq |z'-\hat{z}'|\leq (5b)^{m'-1}$
respectively. Since $|z'-\hat{z}'|=|z-\hat{z}|+\mathcal{O}(
b^{2n})$ and $ b^{2n}$ is negligible  when compared to
$|z-\hat{z}|$, we may assume that $m=m'$. We know that
$|\tau(f_{a,b}^n(z))-e_l(z)|\approx(2a\pm1)|z-\hat{z}|$. Since
$l\geq 2m$, property \eqref{item:cont-field-em-en} of Section
\ref{subsec:contractive-vector-field} leads to
$|e_l(z)-e_{2m}(z)|=\mathcal{O}(b^{2m})$. As a consequence we have
\[
|\tau(f_{a,b}^n(z))-e_{2m}(z)|=|\tau(f_{a,b}^n(z))-e_l(z)|+
\mathcal{O}(b^{2m})\approx(2a\pm1)|z-\hat{z}|,
\]
because $|z-\hat{z}|\geq (5b)^m\gg b^m\gg b^{2m}$.

Observe that $|\tau(f_{a,b}^n(z))-\tau(f_{a',b'}^n(z'))|< b^{2n}$
because $\left|Df_{a,b}^n(z)\tau(z)\right|>\delta\e^{c_2n}$, by
\eqref{eq:EGn}, and
$\left|Df_{a,b}^n(z)\tau(z)-Df_{a',b'}^n(z')\tau(z')\right|<
b^{4n}$. If necessary, we tighten $\U$ in order to guarantee
$|e_{2m}(z)-e'_{2m}(z')|< b^{2n}$. Since $ b^{2n}\ll|z'-\hat{z}'|$
we conclude that
\[
|\tau(f_{a',b'}^n(z'))-e'_{2m}(z')|=|\tau(f_{a,b}^n(z))-e_{2m}(z)|+\O(
b^{2n}) \approx (2a\pm1)|z'-\hat{z}'|.
\]
Finally, a similar argument allows us to obtain
\[
|\tau(f_{a',b'}^n(z'))-e'_{l'}(z)|=|\tau(f_{a',b'}^n(z'))-e'_{2m}(z')|+\mathcal{O}(b^{2m})\approx
(2a\pm1)|z'-\hat{z}'|,
\]
which gives the correct splitting of the vector
$(Df_{a',b'})^n(z')\tau(z')$ with respect to the critical point
$\hat{z}'$.
\end{proof}

Now we will show that if we change the rules of construction of
$\Omega'_\infty$ in the first $N$ iterates by choosing a
convenient binding critical point at each return happening before
$N$ we manage to have $\Omega_{N}=\Omega'_{N}$ as long as
$(a',b')$ is sufficiently close to $(a,b)$.

Before proceeding let us clarify the equality $\Omega_n'=\Omega_n$
for $ n\leq N$. As mentioned in
Remark~\ref{rem:non-uniqueness-omega-infty}, the procedure leading
to $\Omega_\infty$ is not unique. Thus, we have some freedom in
the construction of $\Omega_\infty'$ as long as we guarantee the
slow approximation to the critical set and the expansion along the
tangent direction to $W$.

Take $(a',b')\in\U\cap\B$, where $\U$ is a small neighborhood of
$(a,b)$. Applying the procedure of \cite{BY00} described in
Section \ref{subsec:leading-Cantor-sets} we may build a sequence
of sets $\Omega_0'\supset\Omega_{1}'\supset\ldots$ to obtain
$\Omega_\infty'=\bigcap_{j\in\N_0}\Omega_j'$. From Lemmas
\ref{lem:prox-cp-gen-1} and \ref{lem:suitable-cp} we know that,
given $N$ and $j\leq N$, the set  $\Omega_j$ is a good
approximation of $\Omega_j'$. We propose a modification on the
first $N$ steps in the construction of $\Omega_\infty'$: consider
 $\Omega_n'=\Omega_n$ for all $ n\leq N$; afterwards
make the exclusions of points from $\Omega_{N}$ according to the
original procedure. This way, we produce a sequence of sets
$\Omega_0\supset\ldots\supset\Omega_{N}\supset\Omega_{N+1}'
\supset\ldots$ which we intersect to obtain $\Omega_\infty'$. We
will show that the points in $\Omega_\infty'$  have slow
approximation to the critical set and expansion along the tangent
direction of $W'$ for the dynamics  $f_{a',b'}$. 

When we perturb a parameter $(a,b)\in\B$ and change the rules of
construction of $\Omega_n'$ for a close parameter
$(a',b')\in\U\cap\B$, in the sense mentioned above, we may need to
weaken the condition \eqref{eq:SAn} and introduce condition
\eqref{eq:SAn}' which is defined as \eqref{eq:SAn} except for the
replacement of $\delta$ by $\delta/2$. This way we guarantee the
validity of \eqref{eq:SAn}' for every $(a',b')$ in a sufficiently
small neighborhood $\U$ of $(a,b)$ as stated in next lemma.

\begin{lemma}
\label{lem:prox-omega-infty} Let $(a,b)\in\B$ and $n\in\N$ be
given.
 There is a
neighborhood $\U$ of $(a,b)$ such that for all $
(a',b')\in\U\cap\B$ we may take $\Omega_j'=\Omega_j$ for all
$j\leq n$ and ensure that \eqref{eq:SAn}' holds for all $j\le n$,
for the dynamics $f_{a',b'}$.
\end{lemma}

\begin{proof}
If $\U$ is sufficiently small, then by
Corollary~\ref{prop:prox-crit-points} we have that \eqref{eq:SAn}'
holds for $n=0$, in $H'(\Omega_0)$, for the dynamics $f_{a',b'}$.
Let us suppose that \eqref{eq:SAn}' holds in $H'(\Omega_{n-1})$,
for $f_{a',b'}$ and  $ j\leq n-1<N$. This is to say that for all
$x\in\Omega_{n-1}$ the $f_{a',b'}$ orbit of $z'=H'(x)$ is
controlled up to $n-1$ and at each return $k\leq n-1$, if
$\hat{z}'$ denotes a suitable binding critical point, then
$|f_{a',b'}^k(z')-\hat{z}'|\ge\delta\e^{-\alpha k}/2$.

Our aim is to show that by tightening $\U$, if necessary, this
last statement remains true for $n$. Let $I\subset\Omega_{n-1}$ be
a connected component and $\tilde{I}=I\cap\Omega_n$. Then, by
Lemma \ref{lem:suitable-cp}, we can tighten $\U$, so that for all
$x\in\Omega_{n-1}$, the orbit of $z'=H'(x)$ under $f_{a',b'}$ is
controlled up to $n$. Moreover, if $n$ is a return time for $z'$,
and $\hat{z}'$ is a suitable binding point for $f_{a',b'}^n(z')$,
then $|f_{a',b'}^n(z')-\hat{z}'|\ge\delta\e^{-\alpha n}/2$. Since
each $\Omega_n$ has a finite number of connected components and we
only wish to carry on this procedure up to $N$, then at the end we
still obtain a neighborhood $\U$ of $(a,b)$. \end{proof}

Thus, for every $(a',b')\in\U\cap\B$, where $\U$ is given by
Lemma~\ref{lem:prox-omega-infty}, we have a sequence of sets
$\Omega_0\supset\ldots\supset \Omega_{N}$ such that
\eqref{eq:SAn}' holds  for every $z'=H'(x)$ with $x\in\Omega_{N}$
and $n\leq N$. At this point we proceed with the method described
in Section \ref{sec:horseshoe} and make exclusions out of
$\Omega_{N}$ to obtain a sequence $\Omega_0\supset\ldots\supset
\Omega_{N}\supset\Omega_{N+1}' \supset\ldots$ whose intersection
we denote by
$\Omega_\infty'$.\index{aaOmegaginf'@$\Omega_\infty'$|textbf}
Hence, every point in $H'(\Omega_\infty')$ satisfies
\eqref{eq:SAn} for every $n>N$.

\begin{corollary}
\label{prop:prox-omega-infty} Let $(a,b)\in\B$ and $\varepsilon>0$
be given. There exists  a neighborhood $\U$ of $(a,b)$ so that $
\left|\Omega_\infty\bigtriangleup\Omega_\infty'\right|<\varepsilon
 $
for each $(a',b')\in\U\cap\B$.
\end{corollary}

\begin{proof}
We appeal to Lemma \ref{lem:choice-N1} and find
$N_1=N_1(\varepsilon)$ such that
$|\Omega_{N_1}\setminus\Omega_\infty|<\varepsilon/2$. Observe
that, using Lemma~\ref{lem:prox-omega-infty}, the same $N_1$
allows us to write that
$|\Omega_{N_1}\setminus\Omega_\infty'|<\varepsilon/2$ for all
$(a',b')\in\U\cap\B$. So, we have $
|\Omega_\infty\bigtriangleup\Omega'_\infty|\leq
|\Omega_{N_1}\bigtriangleup\Omega_\infty|+
|\Omega_{N_1}\bigtriangleup\Omega'_\infty|<\varepsilon $.
\end{proof}

\section{Proximity of stable curves}
\label{sec:Prox-stable-curves} So far we have managed to prove
proximity of the horseshoes in the horizontal direction. The goal
of this section is to show the closeness of the stable curves. The
main result of this section is
Proposition~\ref{prop:prox-stable-curves}.

Recall that each long stable curve is obtained as a limit of
``temporary stable curves'', $\gamma_n$, as described in Section
\ref{subsec:stable-curves}. In order to obtain proximity of long
stable curves for close Benedicks-Carleson dynamics we must
produce first an integer $N_2$ such that the approximate stable
curves $\gamma_{N_2}$ are sufficiently close to the corresponding
stable curves $\gamma^s$, regardless of the parameter
$(a,b)\in\B$. This is accomplished through Lemma
\ref{lem:choice-N2}. Therefore, in Proposition
\ref{prop:prox-stable-curves} we obtain the proximity of the
``temporary stable curves''  $\gamma_{N_2}$  for close
Benedicks-Carleson parameters and deduce  in this way  the desired
proximity of the long stable curves.




We use the notation $\gamma_{n}(\zeta)(t)$ or its shorter version,
$\gamma_{n}^t(\zeta)$,\index{aagammant@$\gamma_{n}(z)(t)$,
$\gamma_{n}^t(z)$} for the solution of the equation
$\dot{z}=e_{n}(z)$ with initial condition
$\gamma_{n}(\zeta)(0)=\gamma_{n}^0(\zeta)=\zeta$. Recall that
$\|e_n\|=1$ and $\gamma_n(\zeta)$ is an $e_n$-integral curve of
length $20b$ centered at $\zeta$. So the natural range of values
for $t$ is $[-10b,10b]$.

\begin{lemma}
\label{lem:choice-N2} Let $(a,b)\in\B$  and $n\in\N$ be given.
Consider a connected component $\omega\subset H(\Omega_{n-1})$ and
the rectangle $Q_{n-1}(\omega)$ foliated by the curves $\gamma_n$.
Then the width of the rectangle $Q_{n-1}(\omega)$ is at most
$4\delta^{-1}\e^{-c_2n}$.
\end{lemma}
\begin{proof}
 By the derivative estimate in Subsection
 \ref{subsubsec:derivative-estimate},
 for all $z\in \omega$ we have
$$\left|Df^n(z)\tau(z) \right|>\delta \e ^{c_2n}.$$ Since $\omega$
is a connected component of $H(\Omega_{n-1})$ we have that
$|f^n(\omega)|<2$. As a consequence,
$|\omega|<2\delta^{-1}\e^{-c_2n}$. Observe that this argument also
gives that if $z\in H(\Omega_\infty)$ and $\omega_j$ denotes the
connected component of $H(\Omega_j)$ containing $z$ then $\cap_j
\omega_j=\{z\}$. Let $z^+$ and $z^-$ denote respectively the right
and left endpoints of $\omega$. 
Given $t\in[-10b,10b]$
\begin{align*}
\left|\gamma_{n}^t(z^+)-\gamma_{n}^t(z^-)\right|&\leq
\left|z^++\int_0^t e_{n}\left(\gamma_{n}^r(z^+)\right)dr-
z^--\int_0^te_{n}\left(\gamma_{n}^r(z^-)\right)dr\right|\\
&\leq |z^+-z^-|+\int_0^t
\left|e_{n}\left(\gamma_{n}^r(z^+)\right)-
e_{n}\left(\gamma_{n}^r(z^-)\right)\right|dr\\
&\leq |z^+-z^-|+5\int_0^t
\left|\gamma_{n}^r(z^+)-\gamma_{n}^r(z^-)\right|dr,\, \text{by
\eqref{item:cont-field-Lip-x-dir} of Section
\ref{subsec:contractive-vector-field}}\\
&\leq |z^+-z^-|\e^{5|t|},\,
\text{by a Gronwall type inequality}\\
& \leq |z^+-z^-|\e^{50b}<2|z^+-z^-|=2|\omega|
\end{align*}
Thus, the width of the rectangle $Q_{n-1}(\omega)$ is at most
$4\delta^{-1}\e^{-c_2n}$.
\end{proof}

We will use the following notation for parameters $(a',b')$ close
to $(a,b)$. For any $n\in\N$ and $z'\in\omega_{n}'\subset
H'(\Omega_{n}')$, we denote by $\gamma_{n+1}'(z')$
\index{aagamman'@$\gamma_n'(z')$}the $e_{n+1}'$- integral curve of
length $20b$ centered at $z'$.
  Given $n\in\N$, for any connected component
$\omega' \subset H'(\Omega_n')$ we denote by
$Q_{n}(\omega')=\cup_{z'\in
\omega'}\gamma_{n+1}'(z')$\index{aaQgn'@$Q_n(\omega')$} the
rectangle foliated by the curves $\gamma_{n+1}'(z')$. We define
$Q^1_{n}(\omega')$\index{aaQgn1'@$Q^1_n(\omega')$} as a
$(Cb)^{n+1}$- neighborhood of $Q_{n}(\omega')$ in $\R^2$. Finally,
given $n\in\N$ and any interval
 $\omega\subset H(\Omega_{n})$, we
 denote by $Q^2_{n}(\omega)$\index{aaQgn2@$Q^2_n(\omega)$} a $2(Cb)^{n+1}$-
 neighborhood of $Q_{n}(\omega)$.

\begin{lemma}
\label{lem:prox-stable-curves} Let $(a,b)\in\B$, $n\in\N$,
$\varepsilon>0$ be given, and fix a connected component $I$ of
$\Omega_{n-1}$. Then there is a neighborhood  $\U$ of $(a,b)$ such
that $e_n$, $e_n'$ are defined in $Q^2_{n-1}(H(I))$ and for every
$x\in I$
\begin{equation*}
\label{eq:prox-approx-stable-curves}
\|\gamma_n(H(x))-\gamma_n'(H'(x))\|_0<\varepsilon.
\end{equation*}
Moreover, for every interval $J\subset I$ we have that
$Q^2_{n-1}(H(J))$ contains $Q^1_{n-1}(H'(J))$.
\end{lemma}

\begin{proof}
As we are only interested in arbitrarily small $\varepsilon$, we may
assume that $\varepsilon<b^{2n}$. Take the neighborhood $\U$ of
$(a,b)$ given by Lemma~\ref{lem:prox-omega-infty} applied to $n$.
Within $\U\cap\B$, the set $\Omega_\infty'$ is built out of
$\Omega_{n}$, in the usual way.

 Consider the sequence
$I_0\supset\ldots I_j\supset\ldots\supset I_n=I$ of the connected
components (intervals) $I_j$ of $\Omega_j$ containing $I$. For
every $j\leq n$, let $\omega_j=H(I_j)$ and  $\omega_j'=H'(I_j)$.
We will use a finite inductive scheme such that at step $j$, under
the hypothesis that $e_j$ and $e_j'$ are both defined in
$Q_{j-2}(\omega_{j-1})$, we tighten $\U$ (if necessary) so that
for all $x\in I_{j-1}$ we have $\gamma_j(z)$ $\varepsilon$-close
 to $\gamma_j'(z')$ in the $C^0$ topology, where $z=H(x)$ and
 $z'=H'(x)$, which implies that $Q^2_{j-1}(\omega_j)$ contains
 $Q^1_{j-1}(\omega_j')$. This way we conclude that both
 $e_{j+1}$ and $e_{j+1}'$ are defined in the set
 $Q^2_{j-1}(\omega_j)$, which makes our hypothesis true for step
 $j+1$.
 \begin{figure}[h]
\includegraphics[scale=1]{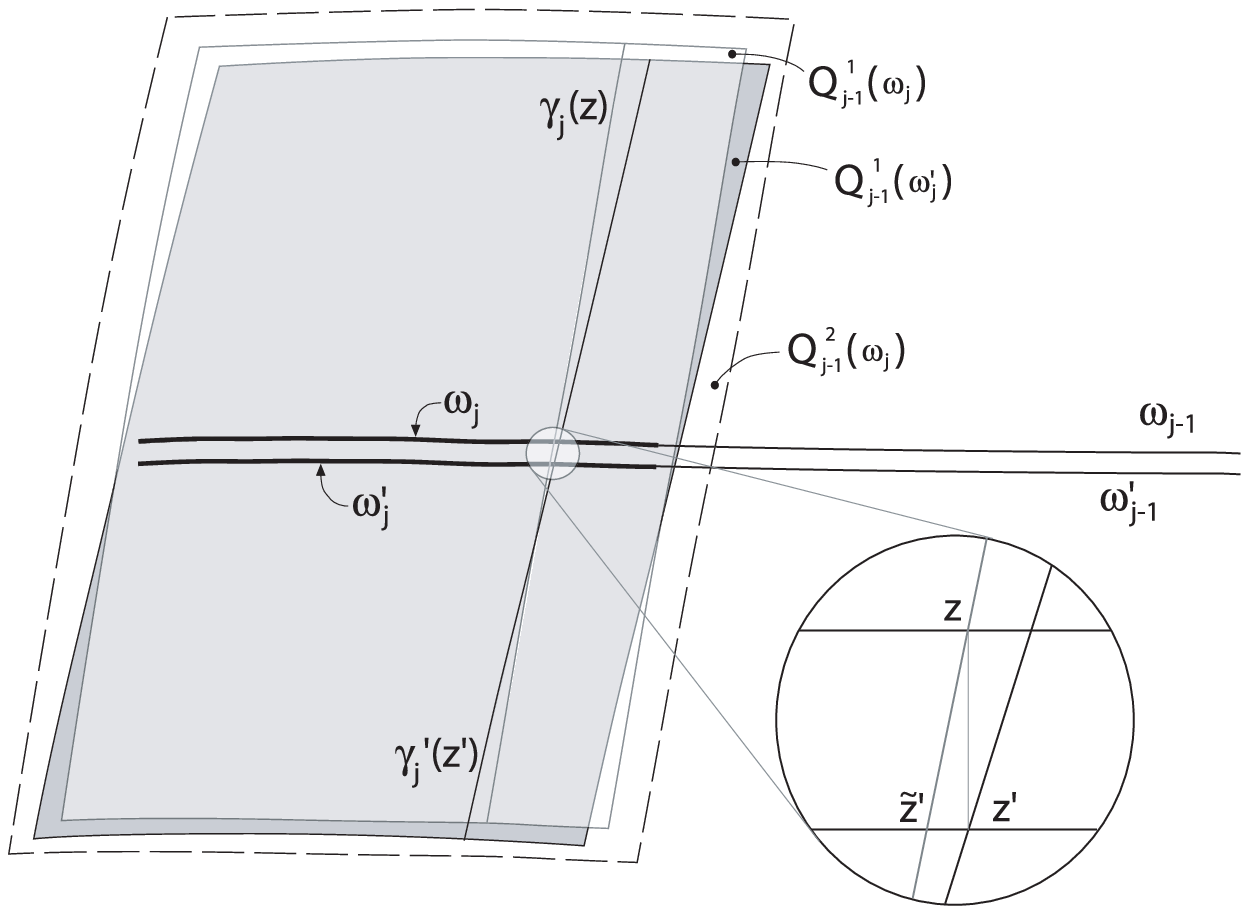}
\caption{} \label{fig:rectangles}
\end{figure}
After $n$ steps we still have a vicinity
 $\U$ of $(a,b)$ and $\gamma_{n}(z)$ is $\varepsilon$ $C^0$-close
 to $\gamma_{n}'(z')$.

We know that $e_1$ and $e_1'$ are defined everywhere in $\R^2$,
which makes our hypothesis true at the first step.

Suppose now, by induction, that at step $j$ we know that $e_{j}$
and $e_{j}'$ are both defined in $Q^2_{j-2}(\omega_{j-1})$, which
contains both $Q^1_{j-2}(\omega_{j-1})$ and
$Q^1_{j-2}(\omega_{j-1}')$. Let $\U$ be sufficiently small so that
for all $(a',b')\in\U\cap\B$, we have $\|H-H'\|_r<\varepsilon^3$
and $|e_j(z)-e_j'(z)|<\varepsilon$, for every $z\in
Q^2_{j-2}(\omega_{j-1})$. Since $Q^1_{j-1}(\omega_{j-1})\subset
Q^1_{j-2}(\omega_{j-1})$ and $Q^1_{j-1}(\omega_{j-1}')\subset
Q^1_{j-2}(\omega_{j-1}')$ (see
\eqref{eq:construction-gamma-infty}), the curves $\gamma_j(z)$ and
$\gamma_j(z')$ never leave the set $Q^2_{j-2}(\omega_{j-1})$, for
every $z\in\omega_{j-1}$ and $z'\in\omega_{j-1}'$.

 Let
$\{\tilde{z}'\}=\gamma_j(z)\cap W_1'$; since
$\|H-H'\|_r<\varepsilon^3$ then $|\tilde{z}'-z|<\varepsilon^2$,
$|\tilde{z}'-z'|<\varepsilon^2$ (see Figure~\ref{fig:rectangles}).
Using the Lipschitzness of the fields $e_j$ and $e_j'$ (property
\eqref{item:cont-field-Lip-x-dir} in Section
\ref{subsec:contractive-vector-field}), the continuity of flows with
initial conditions and the continuity of flows as functions of the
vector field (see for example \cite{HS74}) we have for all $t$
\begin{equation*}
\begin{split}
\left|\gamma_j(\tilde{z}')(t)-\gamma_j'(z')(t)\right|&\leq
\left|\gamma_j(\tilde{z}')(t)-\gamma_j'(\tilde{z}')(t)\right|+
\left|\gamma_j'(\tilde{z}')(t)-\gamma_j'(z')(t)\right|\\
&\leq \frac{\varepsilon}{2a+\mathcal{O}(b)}(\e^{5|t|}-1)+|\tilde{z}'-z'|\e^{5|t|}\\
&\leq \frac{\varepsilon}{3}\e^{50b}50b+2\varepsilon^2< \varepsilon
\end{split}
\end{equation*}
Thus $\|\gamma_j(z)-\gamma_j'(z')\|_0<\varepsilon$. Moreover,
since $\varepsilon\ll(Cb)^j$, we easily get that for any interval
$J\subset I_{j-1}$, the rectangle $Q^2_{j-1}(H(J))$ contains both
$Q^1_{j-1}(H(J))$ and $Q^1_{j-1}(H'(J))$.

From \cite[Section 3.3]{BY00} we know $e_{j+1}$ is defined in a
$3(Cb)^j$-neighborhood in $\R^2$ of $\gamma_j(z)$, for every
$z\in\omega_j$. Since the same applies to $\gamma_j'(z')$ where
$z'\in\omega_j'$ and clearly $\gamma_j(z)$ lies inside a
$(Cb)^j$-neighborhood in $\R^2$ of $\gamma_j'(z')$
($\varepsilon\ll(Cb)^j$) then $e_{j+1}'$ is defined in all points of
$\gamma_j(z)$. This also implies that $e_{j+1}'$ is defined in
$Q^2_{j-1}(\omega_j)$. Thus applying the argument above $n$ times we
get that $e_{n}$ and $e_{n}'$ are defined in
$Q^2_{n-2}(\omega_{n-1})$ and for every $z\in\omega_{n-1}$,
$z'=H'\left(H^{-1}(z)\right)\in\omega_{n-1}'$,
\[
\|\gamma_{n}(z)-\gamma_{n}'(z')\|_0<\varepsilon,
\]
which gives that for any interval $J\subset \Omega_{n-1}$ we have
that $Q^2_{n-1}(H(J))$ contains both $Q^1_{n-1}(H(J))$ and
$Q^1_{n-1}(H'(J))$, since $\varepsilon\ll b^{n}$.

\end{proof}

\begin{proposition}
\label{prop:prox-stable-curves} Let $(a,b)\in\B$ and
$\varepsilon>0$ be given.  There is a neighborhood $\U$ of $(a,b)$
such that  for all $(a',b')\in \U\cap\B$ and $x\in
\Omega_\infty\cap\Omega_\infty'$, we have that $\gamma^s(H(x))$
and $\gamma'^s(H'(x))$ \index{aagammas'@$\gamma'^s(z')$}are
$\varepsilon$-close in the $C^1$ topology.
\end{proposition}

\begin{proof}
Choose $N_2\in\N$ large enough so that
\begin{equation}
\label{eq:N2-definition}
4\delta^{-1}\e^{-c_2N_2}+4(Cb)^{N_2}<\frac{\varepsilon}{3}
\end{equation}
By Lemma \ref{lem:choice-N2} the width of the rectangle
$Q^2_{N_2-1}(\omega_{N_2-1})$ is less than
$\frac{\varepsilon}{3}$.
This means that for every $\zeta\in \omega_{N_2}$, the curve
$\gamma_{N_2}(\zeta)$ is at least $\frac{\varepsilon}{3}$-close to
$\gamma^s(z)$ in the $C^0$ topology. Note that the choice of $N_2$
does not depend on the point $z\in H(\Omega_\infty)$ taken,
neither on the parameter $(a,b)\in\B$ in question.

Take the neighborhood $\U$  of $(a,b)$ to be such that Lemma
\ref{lem:prox-omega-infty} applies up to $N_2$ and Lemma
\ref{lem:prox-stable-curves} applies with $N_2$ replacing $n$. In
particular, for parameters $\U\cap\B$, the set $\Omega_\infty'$ is
built out of $\Omega_{N_2}$, in the usual way and
$Q^2_{N_2-1}(H(I))$ contains $Q^1_{N_2-1}(H'(I))$ for every
 connected component $I\subset \Omega_{N_2-1}$. Moreover, for any $x\in
 I$,
 $\|\gamma_{N_2}(H(x))-\gamma_{N_2}'(H'(x))\|_0<b^{2N_2}$.

Let $x\in \Omega_\infty\cap\Omega_\infty'$ and consider the
sequence $I_0\supset I_1\supset \ldots I_j\supset \ldots$ of the
connected components (intervals) $I_j$ of $\Omega_j$ containing
$x$. Let $z=H(x)$, $z'=H'(x)$ and, for every $j< N_2$, set
$\omega_j=H(I_j)$ and  $\omega_j'=H'(I_j)$. Collecting all the
information we get for any $\zeta\in\omega_{N_2-1}$,
$\zeta'=H'\left(H^{-1}(\zeta)\right)\in\omega_{N_2-1}'$
\begin{equation*}
\left\|\gamma^s(z)-\gamma'^s(z')\right\|_0\leq
\left\|\gamma^s(z)-\gamma_{N_2}(\zeta)\right\|_0+
\left\|\gamma_{N_2}(\zeta)-\gamma_{N_2}'(\zeta')\right\|_0+
\left\|\gamma_{N_2}'(\zeta')-\gamma'^s(z')\right\|_0 <\varepsilon.
\end{equation*}

So far we have proved $C^0$-closeness of the stable leaves. The
fact that the fields $e_n$ and $e_n'$ are Lipschitz with uniform
Lipschitz constant $3<2a+\mathcal{O}(b)<5$ allows us 
to improve the previous $C^0$-estimates to obtain $C^1$-estimates
with little additional effort.
\end{proof}

\section{Proximity of $s$-sublattices and return times
} \label{sec:Prox-R}

The purpose of this section is to obtain the proximity, for close
Benedicks-Carleson dynamics, of the sets of points with the same
history, in terms of free and bound periods up to a fixed time.
In Subsection \ref{subsec:Prox-Upsilons-1} we accomplish this, up
to the first regular return. In Subsection
\ref{subsec:Prox-Upsilons-k} we realize that the same result may
be achieved even if we consider the itineraries up to a some other
return.

\subsection{Proximity after the first return}
\label{subsec:Prox-Upsilons-1}


 Recall that the return time
function $R$ is constant on each $s$-sublattice and, in
particular, on each $\gamma^s$. Thus, the return time function $R$
needs only to be defined in $\Lambda\cap W_1$ or in its vertical
projection in the $x$-axis $\Omega_\infty$. Let
$(\Upsilon_{n,j})_j$ denote the family of  subsets of $\Omega_0$
for which $\bar{\pi}^{-1}(H(\Upsilon_{n,j}))\cap \Lambda$
correspond to the $s$-sublattices of $\Lambda$ given by
\cite[Proposition A]{BY00} and such that $R(H(\Upsilon_{n,j}))=n$.
Observe that $\Upsilon_{n,j}$ determines univocally the
corresponding $s$-sublattice and we allow some imprecision by
referring ourselves to $\Upsilon_{n,j}$ as an $s$-sublattice. The
advantage of looking at the $s$-sublattices as projected subsets
on the $x$-axis is that we can compare these projections of the
$s$-sublattices of different dynamics since all of them live in
the same interval, $\Omega_0$, of the $x$-axis. In Proposition
\ref{prop:difference-R-R'} we obtain proximity of all the
$s$-sublattices $\Upsilon_{n,j}$, with $n\leq N$, for a fixed
integer $N$ and sufficiently close Benedicks-Carleson parameters.

Let us give some insight into the argument. We consider
$(a,b)\in\B$ and $\Omega_\infty$ built according to
Section~\ref{sec:horseshoe}. Let $N\in\N$ be given. We make some
modifications in the procedure described in
Subsection~\ref{subsubsec:rules-def-R} where the $s$-sublattices
are defined so that for each $\Upsilon_{n,j}$, where $n\leq N$, we
obtain an approximation $\Upsilon_{n,j}^*\supset\Upsilon_{n,j}$
whose accuracy depends on the choice of a large integer $N_3$.
Moreover, using Lemmas \ref{lem:prox-omega-infty} and
\ref{lem:prox-stable-curves} we realize that, by construction,
$\Upsilon_{n,j}^*$ also suits as an approximation of
$\Upsilon_{n,j}'\subset\Upsilon_{n,j}^*$, which is an
$s$-sublattice corresponding to $\Upsilon_{n,j}$ for a
sufficiently close $(a',b')\in\B$. The result follows once we
verify that
$|\Upsilon_{n,j}^*-\Upsilon_{n,j}|\approx|\Upsilon_{n,j}^*-
\Upsilon_{n,j}'|$.
Recall that, by construction, for each $\Upsilon_{n,j}$ there are
$I\in\tilde{\P}_{n-1}$, $\omega=H(I)$ and $n$ a regular return
time for $\omega$ such that
\[
\Upsilon_{n,j}=H^{-1}\left(f^{-n}(H(\Omega_\infty))\cap\omega\cap
\Lambda\right)=H^{-1}\left(f^{-n}(H(\Omega_\infty))\cap\omega\right).
\]
 Observe that since $f^n(\omega)\geq 3|\Omega_0|$ then
$\omega$ has a minimum length $|\omega|\geq 5^{-n}3|\Omega_0|$.
This means that for $n$ fixed there can only be a finite number of
$\Upsilon_{n,j}$'s. In fact, if $v(n)$ denotes the number of
$\Upsilon_{n,j}$ with $R(H(\Upsilon_{n,j}))=n$\index{aavn@$v(n)$},
then
\begin{equation}
\label{eq:vn-bound} v(n)\leq
\frac{|\Omega_0|}{5^{-n}3|\Omega_0|}\leq5^n.
\end{equation}

Let $N\in\N$ be given and let $N_3>2N$ be a large integer whose
choice will be specified later. Let $\varepsilon<b^{2N_3}$ be
small. Consider $\U$ small enough so that condition
\eqref{eq:CR-omega0-proximity} holds for such an $\varepsilon$ and
$\Omega_j=\Omega_j'$ for all $j\in\{0,\ldots, N_3\}$ (recall
Lemma~ \ref{lem:prox-omega-infty}), while $\Omega_\infty'$ is
built in usual way out of $\Omega_{N_3}'$.

For $n\leq N$ we carry out an inductive construction of sets
$\tilde{\Omega}_n^*\subset \Omega_n$
\index{aaOmegagtilde*@$\tilde{\Omega}_n^*$}and partitions
$\tilde{\mathcal{P}}^*_n$ \index{aaPgtilde*@$\tilde{\P}_n^*$}of
$\tilde{\Omega}_n^*$ that will
coincide for all $(a',b')\in\U$, for every $n\leq N$. 
This process must ensure that for every $n\leq N$  we have
$\tilde{\Omega}_n^*\subset\tilde{\Omega}_n$, and if $\omega^*\in
H(\tilde{\mathcal{P}}^*_n)$, then there is $\omega\in
H(\tilde{\mathcal{P}}_n)$ such that $\omega\supset\omega^*$.
Moreover, by choice of $N_3$ we will have that
$\omega\setminus\omega^*$, when not empty, occupies the tips of
$\omega$ and it corresponds to such a small part that if $\omega$
has a regular return at time $n<j\leq N$ then
$f^j(\omega^*)\supset 2 \Omega_0$ still traverses $Q_0$ by wide
margins (see Lemma~\ref{lem:omega-omega*}).

\subsubsection{Rules for defining $\tilde{\Omega}_n^*$,
$\tilde{\mathcal{P}_n}^*$ and $R^*$.}\index{aaRg*@$R^*(z)$}

{\renewcommand{\theenumi}{\arabic{enumi}*}
\begin{enumerate}
\setcounter{enumi}{-1}
\item $\tilde{\Omega}_0^*=\tilde{\Omega}_0'^{*}=\Omega_0$,
$\tilde{\mathcal{P}}_0^*=\tilde{\mathcal{P}}_0'^{*}=\{\Omega_0^*\}$. 
\end{enumerate}
}

Assume that $\tilde{\Omega}_{n-1}^*=\tilde{\Omega}_{n-1}'^{*}$ and
that for each $I^*\in \tilde{\mathcal{P}}^*_{n-1}$ there is $I\in
\tilde{\mathcal{P}}_{n-1}$ such that $I\supset I^*$. Take $I^*\in
\tilde{\mathcal{P}}_{n-1}^{*}=\tilde{\mathcal{P}}_{n-1}'^{*}$. We
denote  $\omega=H(I)$, $\omega^*=H(I^*)$ and
${\omega^*}'=H'(I^*)$.\index{aaOmegagtilde*'@$\tilde{\Omega}_n'^{*}$}
\index{aaPgtilde*'@$\tilde{\P}_n'^{*}$}
{\renewcommand{\theenumi}{\arabic{enumi}*}
\begin{enumerate}
\item
\label{item:rule1_construction_Pn*} If
$\omega\in\tilde{\mathcal{P}}_{n-1}$ does not make a regular
return to $H(\Omega_0)$ at time $n$, put
$\tilde{I}^*=I^*\cap\Omega_n$ into $\tilde{\Omega}_{n}^*$ and let
$
\tilde{\mathcal{P}}_{n}^{*}\left|_{\tilde{I}^*}\right.=H^{-1}\left(
f_{a,b}^{-n}\mathcal{P}\left|_{H(\tilde{I}^*)}\right.\right)$ with
the usual adjoining of intervals.

\end{enumerate}
}

 We remark that if we were to apply this rule directly to
 $(a',b')\in\U\cap\B$, where $\U$ is
sufficiently small so that Corollary~\ref{prop:prox-crit-points},
Lemma~\ref{lem:prox-omega-infty} and equation
\eqref{eq:CR-omega0-proximity} hold for such $\varepsilon$ and
$N_3$, then $\tilde{\Omega}_{n}'^{*}$ and
$\tilde{\mathcal{P}}_{n}'^{*}$ would have discrepancies of
$\mathcal{O}(\varepsilon)$ relative to $\tilde{\Omega}_{n}^*$ and
$\tilde{\mathcal{P}}_{n}^{*}$ built for $(a,b)$, respectively. But
$\varepsilon<\e^{-2N_3}$ is negligible  when compared to
$\e^{-\alpha N}$ or ${\e^{-\alpha N}}/{N^2}$. Observe that  the
points of $H(\tilde{I}^*)$ never get any closer than $\e^{-\alpha
N}$ from the critical set, up to time $n$, and ${\e^{-\alpha
N}}/{N^2}$ is the minimum size of the elements of the partition
$\mathcal{P}$ whose distance to the critical set is larger than
$\e^{-\alpha N}$. Hence, there is no harm in setting
$\tilde{\Omega}_{n}'^{*}=\tilde{\Omega}_{n}^*$ and
$\tilde{\mathcal{P}}_{n}'^{*}=\tilde{\mathcal{P}}_{n}^{*}$.

 Let $\mathcal{S}_{N_3}$\index{aaSgn@$\mathcal S_n$}
be the partition of $\Omega_{N_3}$ into connected components. We
clearly have $\#\mathcal{S}_{N_3}\leq 2^{N_3}$. We write
$f^n(z)\in H(\Omega_{N_3})$ if there exists $\sigma\in
\mathcal{S}_{N_3}$ such that $f^n(z)\in Q^2_{N_3-1}(H(\sigma))$
where, as before, $Q^2_{N_3-1}(H(\sigma))$ is a
$2(Cb)^{N_3}$-neighborhood of $Q_{N_3-1}(H(\sigma))$ in $\R^2$.
This way let $f^{-n}(H(\Omega_{N_3}))$ have its obvious meaning.
Observe that by definition of $Q^2_{N_3-1}(H(\sigma))$ and the
construction of the long stable curves (namely
\eqref{eq:construction-gamma-infty}), then
\begin{equation}
\label{eq:inverse-images-omegaN4-omega-infty}
f^{-n}(H(\Omega_{N_3}))\supset f^{-n}(H(\Omega_\infty)),
\end{equation}
where we write $f^n(z)\in H(\Omega_\infty)$ when
$f^n(z)\in\gamma^s(\zeta)$ for some $\zeta\in H(\Omega_\infty)$.

\begin{figure}[h]
\includegraphics[scale=0.9]{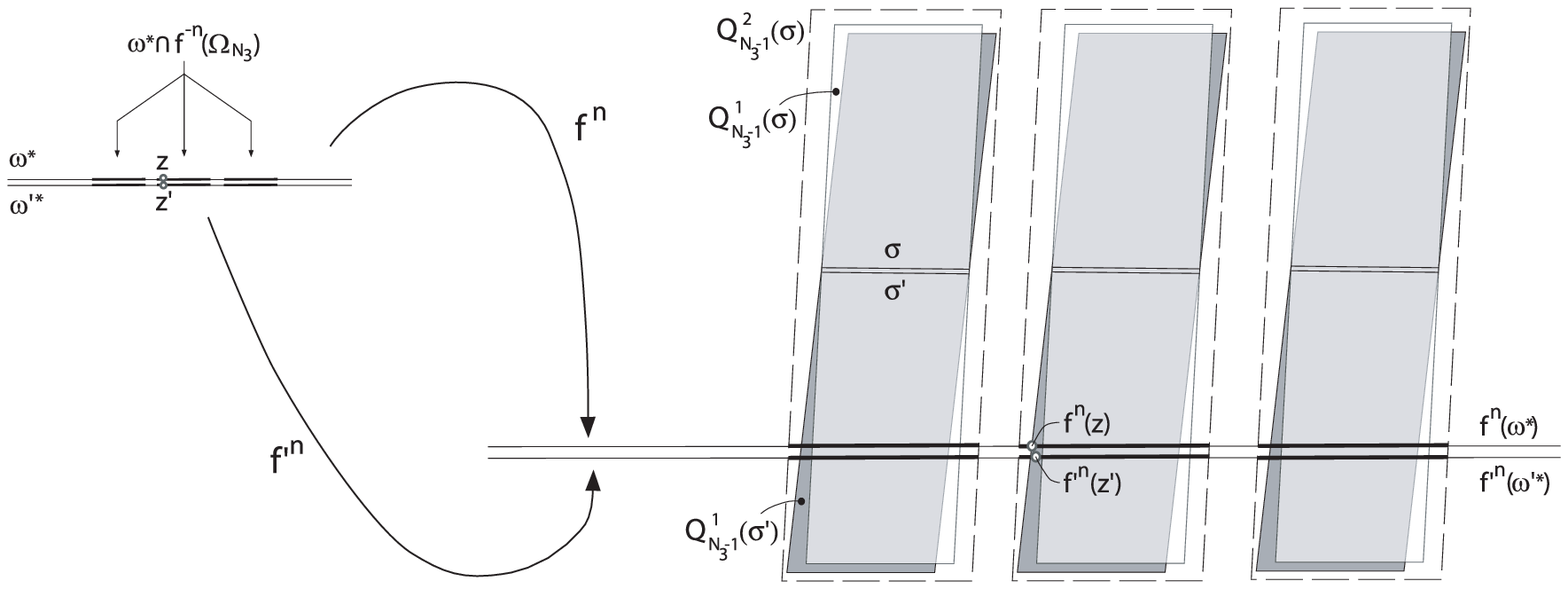}
\caption{} \label{fig:rectangles-procedure}
\end{figure}

{\renewcommand{\theenumi}{\arabic{enumi}*}
\begin{enumerate}
\setcounter{enumi}{1}
\item
\label{item:rule2_construction_Pn*} If
$\omega\in\tilde{\mathcal{P}}_{n-1}$ makes a regular return at
time $n$, we put
$$\tilde{I}^*=H^{-1}\left(\omega^*\setminus
f^{-n}(H(\Omega_{N_3}))\right)\cap\Omega_n$$ into
$\tilde{\Omega}_n^*$. Let $\mathcal{S}^*$ be the partition of
$\tilde{I}^*$ into connected components. We define
$\tilde{\mathcal{P}}_{n}^{*}\left|_{\tilde{I}^*}\right.=H^{-1}\left(
f^{-n}\mathcal{P}\left|_{H(\tilde{I}^*)}\right.\right)
\bigvee\mathcal{S}^*$. For $z\in\omega^*$ such that $f^n(z)\in
H(\Omega_{N_3})$ we define $R^*(z)=n$.

\end{enumerate}
}

Suppose that $\U$ is sufficiently small so that as in
Lemma~\ref{lem:prox-stable-curves} we have
$Q^2_{N_3-1}(H(\sigma))\supset Q^1_{N_3-1}(H'(\sigma))$ and, as
before, Corollary \ref{prop:prox-crit-points} and condition
\eqref{eq:CR-omega0-proximity} hold for the considered
$\varepsilon$ and $N_3$. Then, the smallness  of
$\varepsilon<b^{2N_3}$ when compared to the sizes of the elements
$f^n\big(H(\tilde{\mathcal{P}}_{n}^{*})\big)$ for $ n\leq N$
allows us to consider
$\tilde{\Omega}_{n}^*=\tilde{\Omega}_{n}'^{*}$ and
$\tilde{\mathcal{P}}_{n}^{*}=\tilde{\mathcal{P}}_{n}'^{*}$.

%
%

Essentially in this construction we substitute $\Omega_\infty$ by
its finite approximation $\Omega_{N_3}$ in order to relate the
partitions built for $(a,b)$ with the ones built for
$(a',b')\in\U\cap\B$. The sets
$\{R=n\}=\tilde{\Omega}_{n-1}\setminus\tilde{\Omega}_n$ were
defined as the sets of points that at time $n$ had their first
regular return to $H(\Omega_\infty)$ (after sliding along
$\gamma^s$ stable curves). Now
$\{R^*=n\}=\tilde{\Omega}_{n-1}^*\setminus\tilde{\Omega}_n^*$ is
defined as the set of points that at time $n$ have their first
regular return to $H(\Omega_{N_3})$, where the sliding is made
along the stable curve approximates, $\gamma_{N_3}$.

Let us make clear some aspects related to the previous rules. When
we apply rule \eqref{item:rule2_construction_Pn*} at step $n$, we
ensure that for every $z\in\tilde{\Omega}_n^*$ we have $z\notin
f^{-n}(H(\Omega_\infty))$. Let us verify that the same applies to
$f_{a',b'}$, ie, since we are considering $\U$ sufficiently small so
that Lemma~\ref{lem:prox-stable-curves},
Corollary~\ref{prop:prox-crit-points} and condition
\eqref{eq:CR-omega0-proximity} hold for $\varepsilon$ and $N_3$ in
question, then for every $z'\in \tilde{\Omega}_{n}'^{*}$ we have
$z'\notin f_{a',b'}^{-n}(H'(\Omega_\infty'))$ Since $\varepsilon$ is
irrelevant when compared to $2(Cb)^{N_3}$ we have for all $z'\in
H'(\tilde{\Omega}_n^*)$ and for every $\sigma\in\mathcal{S}_{N_3}$,
$\mbox{dist}(f_{a',b'}^n(z'),Q_{N_3-1}(H(\sigma)))>2(Cb)^{N_3}-\varepsilon$
which implies that
$\mbox{dist}(f_{a',b'}^n(z'),Q_{N_3-1}(H'(\sigma)))>(Cb)^{N_3}$,
since by Lemma \ref{lem:prox-stable-curves} we may assume that
$Q^1_{N_3-1}(H'(\sigma))\subset Q^2_{N_3-1}(H(\sigma))$ and
$\mbox{dist}\left( Q_{N_3-1}(H'(\sigma)), Q_{N_3-1}(H(\sigma))
\right)\leq \varepsilon$. We have used ``$\mbox{dist}$'' to denote
the usual distance between  two sets.

In next lemma take into account that since
$\Omega_{N_3}\supset\Omega_\infty$, then the gaps of
$\Omega_\infty$ contain those of $\Omega_{N_3}$, and so for all
$n\leq N$ and $\omega^*\in H(\tilde{\P}_{n-1}^*)$ there exists
$\omega\in
H(\tilde{\P}_{n-1})$ such that $\omega^*\subset\omega$. 

\begin{lemma}
\label{lem:omega-omega*} Let $n\leq N$, $\omega^*\in
H(\tilde{\P}_{n-1}^*)$ and consider $\omega\in
H(\tilde{\P}_{n-1})$ such that $\omega^*\subset\omega$. If $N_3$
is large enough then $f^n(\omega)\setminus f^n(\omega^*)$, when
not empty, occupies one or both tips of $f^n(\omega)$ and
$|f^n(\omega)\setminus f^n(\omega^*)|<|\Omega_0|^2$.
\end{lemma}

\begin{proof}
Let $N_3\in\N$ be sufficiently large so that
\begin{equation}
\label{eq:choice-N3-1}
5^{N}\left(C_1\delta^{1-3\beta}\frac{\e^{-\alpha(1-3\beta)(N_3+1)}}
{1-\e^{-\alpha(1-3\beta)}}+2(Cb)^{N_3}\right)<|\Omega_0|^2.
\end{equation}
For every $i\leq n-1$, let $\omega_i^*\in H(\tilde{\P}_{i}^*)$ be
such that $\omega^*\subset \omega_i^*$ and let $\omega_i\in
H(\tilde{\P}_{i})$ be such that $\omega\subset \omega_i$. If
$\omega\setminus\omega^*\neq\emptyset$ then at some time before
$n-1$, rule \eqref{item:rule2_construction_Pn*} was applied. Let
$j\leq n-1$ be the last moment in the history of $\omega^*$ that
rule \ref{item:rule2_construction_Pn*} was applied. Then,
 $f^j(\omega^*_j)$ hits a gap
of $\Omega_{N_3}$ while $f^j(\omega_j)$ hits a gap of
$\Omega_\infty$. According to Lemma~\ref{lem:choice-N1} the
difference $f^j(\omega_j)\setminus f^j(\omega^*_j)$ has length of
at most
 $$C_1\delta^{1-3\beta}\frac{\e^{-\alpha(1-3\beta)(N_3+1)}}
{1-\e^{-\alpha(1-3\beta)}}+2(Cb)^{N_3},$$
 where the last term results from the
fact that we are using $2(Cb)^{N_3}$- neighborhoods of the
rectangles spanned by the approximate stable curves. Moreover,
$f^j(\omega_j)\setminus f^j(\omega^*_j)$ clearly occupies the tips
of $f^j(\omega_j)$.

Now, for simplicity suppose that $\omega=\omega_j$ and
$\omega^*=\omega^*_j$. We have that $f^n(\omega)\setminus
f^n(\omega^*)$ occupies the tips of $f^n(\omega)$. This geometric
property is inherited since by construction we are away from the
folds and $f$ is a diffeomorphism. Also, up to time $n$,
$|f^j(\omega_j)\setminus f^j(\omega^*_j)|$ can grow no more than
$5^{n-j}$. Consequently, by choice of $N_3$ we must have
$|f^n(\omega)\setminus f^n(\omega^*)|<|\Omega_0|^2$.

In the case that $\omega\neq\omega_j$ it means that $\omega_j$
will suffer exclusions or subdivisions. Nevertheless, the points
of $(\omega_j-\omega_j^*)\cap\omega$ still occupy the tip of
$f^n(\omega)$.
\end{proof}

\begin{remark}
\label{rem:redefine-regular-return} Observe that by choice of
$N_3$ we have that if $\omega^*\in H(\tilde{\P}_{n-1}^*)$ and
$f^n(\omega)$ makes a regular return then $f^n(\omega^*)\supset
(3-|\Omega_0|^2)\Omega_0$. This means that for $\U$ sufficiently
small $f_{a',b'}^n(\omega{'^*})\supset
(3-|\Omega_0|^2-\varepsilon)\Omega_0$.
\end{remark}

When at step $n$ we have to apply rule
\eqref{item:rule2_construction_Pn*} we make more exclusions from
$\tilde{\Omega}_{n-1}^*$ than we would if we were to apply rule
(2) as in \cite{BY00}. Essentially we are excluding the points
that hit $H(\Omega_{N_3})$ instead of only removing the points
that hit $H(\Omega_\infty)$ ($\Omega_\infty\subset\Omega_{N_3}$).
We argue that by adequate choice of $N_3$ this over exclusion will
not affect the sets $\{R^*=j\}$ with $j\in\{n+1,\ldots,N\}$.

\begin{lemma}
\label{lem:harmless-overexclusion} Suppose that $x$ is a point
that at step $n$ should be excluded by rule
\eqref{item:rule2_construction_Pn*} but is not excluded according
to rule \eqref{item:rule2-def-R}. If $N_3$ is large enough, then
$H(x)$ does not have a regular return to $\Omega_0$ before $N$.
\end{lemma}
\begin{proof}
Let $N_3\in\N$ be sufficiently large so that
\eqref{eq:choice-N3-1} holds and take
$\sigma\in\mathcal{S}_{N_3}$. When we apply rule
\eqref{item:rule2_construction_Pn*} at step $n$ we remove from
$\tilde{\Omega}_{n-1}^*$ all the points hitting $H(\sigma)$, while
if we had applied rule \eqref{item:rule2-def-R} instead we would
have only removed the points hitting $H(\sigma\cap\Omega_\infty)$.
Consider a gap $\varpi$ of $H(\sigma\cap\Omega_\infty)$. We know
that the length of $\varpi$ is less than
$$C_1\delta^{1-3\beta}\frac{e^{-\alpha(1-3\beta)N_3}}{1-e^{-\alpha(1-3\beta)}}.$$
If $\partial\varpi\cap\partial H(\sigma)=\emptyset$ then
$\varpi\in\tilde{\P}_n$ and in $N$ iterations it would grow to
reach at most the length
\[5^{N}C_1\delta^{1-3\beta}\frac{e^{-\alpha(1-3\beta)N_3}}{1-e^{-\alpha(1-3\beta)}}
<|\Omega_0|^2\ll 3|\Omega_0|.
\]
Thus, $\varpi$ would not have any regular return to $\Omega_0$
before $N$.

If $\partial\varpi\cap\partial H(\sigma)\neq\emptyset$, then there
is a gap $\hat{\varpi}$ of $H(\Omega_{\infty})$ so that
$\hat{\varpi}\in\tilde{\P}_n$ and $\varpi$ occupies a tip of
$\hat{\varpi}$. Clearly, $\hat{\varpi}$ could have a regular
return at $j\in\{n+1,\dots, N\}$, say. However, by construction
$f^j(\varpi)$ will occupy one tip of $f^j(\hat{\varpi})$. Since
$$\left|f^j(\varpi)\right|<5^{N}C_1\delta^{1-3\beta}
\frac{e^{-\alpha(1-3\beta)N_3}}{1-e^{-\alpha(1-3\beta)}}<|\Omega_0|^2$$
and $|f^j(\hat{\varpi})|\gtrsim 3|\Omega_0|$ we still have that
$f^j(\varpi)$ does not hit $\Omega_0$. We remark that
$\hat{\varpi}$ could have suffered subdivisions and exclusions
according to rule \eqref{item:rule1_construction_Pn*} before time
$j$. Nevertheless, the points from $\varpi$ that survive the
exclusions still occupy the tip of the piece that will contain
them at the time of its regular return and the argument applies
again.
\end{proof}

By the rules in Subsection \ref{subsubsec:rules-def-R}, for every
$s$-sublattice $\Upsilon_{n,j}$ there is a segment
$\omega_{n,j}\in H(\tilde{\P}_{n-1})$ such that $n$ is a regular
return time for $\omega_{n,j}$  and
\begin{equation}
\label{eq:def-upsilon-1}
\Upsilon_{n,j}=H^{-1}\left(\omega_{n,j}\cap
f^{-n}(\Omega_\infty)\right).
\end{equation}
Lemmas \ref{lem:omega-omega*} and \ref{lem:harmless-overexclusion}
allow us to conclude that if
$\omega_{n,j}\in\tilde{\mathcal{P}}_{n-1}$ and $n\leq N$ is a
regular return time for $\omega_{n,j}$ then there is
$\omega^*_{n,j}\in\tilde{\mathcal{P}}_{n-1}^*$ such that
$\omega^*_{n,j}\subset\omega_{n,j}$ and
$|f^n(\omega_{n,j})|=|f^n(\omega^*_{n,j})|+\mathcal{O}(|\Omega_0|^2)$.
Moreover, because the difference between $\omega_{n,j}$ and
$\omega_{n,j}^*$ is only in their tips we may write
\begin{equation}
\label{eq:def-upsilon-2}
\Upsilon_{n,j}=H^{-1}\left(\omega_{n,j}^*\cap
f^{-n}(\Omega_\infty)\right).
\end{equation}
Attending to the procedure above and equation
\eqref{eq:def-upsilon-2}, given an $s$-sublattice
$\Upsilon_{n,j}$, with $n\leq N$ we define its approximation
\begin{equation}
\label{eq:def-upsilon-1*}
\Upsilon_{n,j}^*=H^{-1}\left(\omega_{n,j}^*\cap
f^{-n}(\Omega_{N_3})\right).\index{aaUpsilonn1*@$\Upsilon_{n,j}^*$}
\end{equation}
Taking into consideration
\eqref{eq:inverse-images-omegaN4-omega-infty} we have that
$\Upsilon_{n,j}\subset\Upsilon_{n,j}^*$, from where we conclude
that $\forall n\in \{1,\ldots,N\}$,
\[
\{R=n\}=\bigcup_{j\leq v(n)}\Upsilon_{n,j}\subset\bigcup_{j\leq
v(n)}\Upsilon^*_{n,j}=\{R^*=n\}.
\]
We wish to verify that this substitution of $\Omega_\infty$ by
$\Omega_{N_3}$ does not produce significant changes. In fact, we
will show in the next lemma that $\Upsilon_{n,j}$ and
$\Upsilon_{n,j}^*$ are very close for all $n\leq N$ and $j\leq
v(n)$.

\begin{lemma}
\label{lem:diff-upsilon1*-upsilon1} Let $\varepsilon>0$, $N\in\N$
and an $s$-sublattice $\Upsilon_{n,j}$ with $n\leq N$ be given.
If $N_3$ is large enough, then
\[
\left|\Upsilon^*_{n,j}\setminus\Upsilon_{n,j}\right|<\varepsilon\quad\text{and}\quad
\left|\{R^*=n\}\setminus\{R=n\}\right|<\varepsilon.
\]
\end{lemma}

\begin{proof}
Choose $N_3$ large enough so that
\begin{equation}
\label{eq:choice-N4-2}
C_1|\Omega_{N_3}\setminus\Omega_\infty|<\varepsilon.
\end{equation}
Let $\omega_{n,j}^*$ be such 
that $H(\Upsilon_{n,j})=\omega_{n,j}^*\cap
f^{-n}(H(\Omega_\infty))$ and
$H(\Upsilon_{n,j}^*)=\omega_{n,j}^*\cap f^{-n}(H(\Omega_{N_3}))$.
By bounded distortion we have
\[
\frac{|H(\Upsilon_{n,j}^*)\setminus
H(\Upsilon_{n,j})|}{|\omega_{n,j}^*|}\leq C_1
\frac{|f^n\left(H(\Upsilon_{n,j}^*)\setminus
H(\Upsilon_{n,j})\right)|}{|f^n(\omega_{n,j}^*)|}\leq
C_1\frac{|\Omega_{N_3}\setminus\Omega_\infty|}{2|\Omega_0|}.
\]
Attending to \eqref{eq:choice-N4-2} this gives that
$\left|\Upsilon^*_{n,j}\setminus\Upsilon_{n,j}\right|<\varepsilon$.
Besides,
\begin{equation*}
\begin{split}
|\{R^*=n\}\setminus\{R=n\}|&=\sum_{j\leq
v(n)}|\Upsilon_{n,j}^*\setminus\Upsilon_{n,j}|\leq
            \sum_{j\leq v(n)}C_1\frac{|\Omega_{N_3}\setminus
            \Omega_\infty|}{|\Omega_0|}
            |\omega_{n,j}^*|\\
&\leq C_1|\Omega_{N_3}\setminus\Omega_\infty|\\
&<\varepsilon,
\end{split}
\end{equation*}
by the choice of $N_3$.
\end{proof}

\begin{remark}
By definition of $f^{-n}(H(\Omega_{N_3}))$, in the estimates above
we should have considered
$\left|\Omega^2_{N_3}\setminus\Omega_\infty\right|$, where
$\Omega^2_{N_3}$ is a $2(Cb)^{N_3}$-neighborhood of
$\Omega_{N_3}$. However, since $\Omega_{N_3}$ has at most
$2^{N_3}$ connected components, then the difference to the
estimates above would be at most $2^{N_3+1}(Cb)^{N_3}$, which is
as small as we want if we choose $N_3$ large enough.
\end{remark}
\begin{remark}
The estimates in the proof were used taking $H(\Omega_{N_3})$ and
$H(\Omega_\infty)$ as subsets of $W_1$. According to \cite[Remark
5]{BY00}, upon re-scaling the estimates still work if we consider
them as subsets of $\gamma^u\in\Gamma^u$, due to Lemma 2 of
\cite{BY00}.
\end{remark}

\begin{proposition}
\label{prop:difference-R-R'} Let $(a,b)\in\B$, $N\in\N$ and
$\varepsilon>0$ be given.  There is a neighborhood $\U$ of $(a,b)$
such that for all $ (a',b')\in\U\cap\B$ given any $s$-sublattice
$\Upsilon_{n,j}\subset\Omega_\infty$, with $n\leq N$ and $j\leq
v(n)$, then the corresponding $s$-sublattice
$\Upsilon'_{n,j}\subset\Omega_\infty'$
\index{aaUpsilonn1'@$\Upsilon'_{n,j}$}is such that
\[
|\Upsilon_{n,j}\bigtriangleup\Upsilon_{n,j}'|<\varepsilon
\quad\text{and}\quad
\left|\{R=n\}\bigtriangleup\{R'=n\}\right|<\varepsilon.
\]
\end{proposition}

\begin{proof}
By Lemma \ref{lem:prox-omega-infty} we are assuming that
$\Omega_\infty'$ is built out of $\Omega_{N_3}$ in the usual way
for $f_{a',b'}$ with $(a',b')\in\U\cap\B$. Lemma
\ref{lem:prox-stable-curves} assures that if $\U$ is small enough
then $Q^2_{N_3-1}(H(\sigma))$, which is a $2(Cb)^{N_3}$-
neighborhood
 of $Q_{N_3-1}(H(\sigma))$, contains $Q^1_{N_3-1}(H'(\sigma))$
 for every $\sigma\in \mathcal{S}_{N_3}$. Moreover, for any $x\in
 \sigma$
 $$\|\gamma_{N_3}(H(x))-\gamma_{N_3}'(H'(x))\|_0<b^{N_2+1}.$$
Let $N_3$ be chosen according to equations \eqref{eq:choice-N3-1}
and \eqref{eq:choice-N4-2} so that Lemmas \ref{lem:omega-omega*}
and \ref{lem:diff-upsilon1*-upsilon1} hold. Let $\Upsilon_{n,j}$,
with $n\leq N$, be a given $s$-sublattice of $H(\Omega_\infty)$.
Let $I_{n,j}^*\in\tilde{\P}^*_{n-1}$ be such that
\[
\Upsilon_{n,j}=H^{-1}\left(\omega_{n,j}^*\cap
f^{-n}(H(\Omega_\infty))\right),
\]
where
$\omega_{n,j}^*=H(I_{n,j}^*)$. 
 Suppose that $\U$ is sufficiently small so that the construction
of the partition is carried out simultaneously for the dynamics
$f_{a',b'}$ correspondent to any $(a',b')\in\U\cap\B$ and so that
$\tilde{\P}^*_{m}=\tilde{\P}'^*_{m}$, for all $m\leq N$, as it has
been described in the procedure above. Then,
$f_{a',b'}^n(\omega'^*_{n,j})=f_{a',b'}^n(H'(I_{n,j}^*))$ crosses
$Q_0$ by wide margins and we may define
\[
\Upsilon_{n,j}'=H'^{-1}\left(\omega'^*_{n,j}\cap
f_{a',b'}^{-n}(H'(\Omega_\infty'))\right).
\]
Consider the approximation $\Upsilon_{n,j}^*$ built in
\eqref{eq:def-upsilon-1*} for $\Upsilon_{n,j}$. We have seen that
$\Upsilon_{n,j}\subset\Upsilon_{n,j}^*$ and using Lemma
\ref{lem:diff-upsilon1*-upsilon1} we may suppose that
$\left|\Upsilon^*_{n,j}\setminus\Upsilon_{n,j}\right|<{\varepsilon}/{2}$.
Now, we shall see that $\Upsilon_{n,j}^*$ is also a good
approximation for $\Upsilon_{n,j}'$ if $\U$ is sufficiently small.

First, we verify that $\Upsilon_{n,j}'\subset\Upsilon_{n,j}^*$.
Let $x\in\Upsilon_{n,j}'$, $z=H(x)$ and $z'=H'(x)$. We need to
check that if $f_{a',b'}^n(z')\in\Lambda'$, then $f^n(z)\in
Q^2_{N_3-1}(H(\sigma))$ for some $\sigma\in\mathcal{S}_{N_3}$. We
are supposing that $\U$ is sufficiently small so that
\eqref{eq:CR-omega0-proximity} holds for $\varepsilon<b^{2N_3}$ up
to $N_3$, which implies that $|f^n(z)-f_{a',b'}^n(z')|<b^{2N_3}$.
Since
$\Lambda'\subset\bigcup_{\sigma\in\mathcal{S}_{N_3}}Q^1_{N_3-1}
(H'(\sigma))$, we have $f_{a',b'}^n(z')\in
Q^1_{N_3-1}(H'(\sigma))$ for some $\sigma\in\mathcal{S}_{N_3}$.
Under the assumptions described in the procedure above (namely
that $Q^1_{N_3-1}(H'(\sigma))\subset Q^2_{N_3-1}(H(\sigma))$) and
attending to equation \eqref{eq:prox-approx-stable-curves} we get
that
$\mbox{dist}\left(f_{a',b'}^n(z'),Q_{N_3-1}(H(\sigma))\right)<3/2(Cb)^{N_3}$,
and thus
$\mbox{dist}\left(f^n(z),Q_{N_3-1}(H(\sigma))\right)<2(Cb)^{N_3}$.

Additionally, since the upper bound used for
$|\Omega_{N_3}\setminus\Omega_\infty|$ also works for
$|\Omega_{N_3}\setminus\Omega_\infty'|$ and the width of
$Q^1_{N_3-1}(H'(\sigma))$ differs from the width  of
$Q^2_{N_3-1}(H(\sigma))$ by $\O((Cb)^{N_3})$ we observe that the
 argument used in Lemma \ref{lem:diff-upsilon1*-upsilon1}
gives us that
$\left|\Upsilon^*_{n,j}\setminus\Upsilon_{n,j}'\right|<{\varepsilon}/{2}$.
Therefore
\[
|\Upsilon_{n,j}\bigtriangleup\Upsilon_{n,j}'|\leq
\left|\Upsilon^*_{n,j}\setminus\Upsilon_{n,j}\right|+
\left|\Upsilon^*_{n,j}\setminus\Upsilon_{n,j}'\right|<\varepsilon,
\]
which gives the first part of the conclusion.

Suppose now that Lemma \ref{lem:diff-upsilon1*-upsilon1} holds and
$|\{R^*=n\}\setminus\{R=n\}|<{\varepsilon}/{2}$. Observing that
$\{R'=n\}=\cup_{j\leq v(n)}\Upsilon_{n,j}'$, then arguing as in
Lemma \ref{lem:diff-upsilon1*-upsilon1}, we have
$|\{R^*=n\}\setminus\{R'=n\}|<{\varepsilon}/{2}$, as long as $\U$
is sufficiently small. Finally,
\[
|\{R=n\}\bigtriangleup\{R'=n\}|\leq
|\{R=n\}\bigtriangleup\{R^*=n\}|+|\{R^*=n\}\bigtriangleup\{R'=n\}|<\varepsilon.
\]
\end{proof}

\subsection{Proximity after $k$ returns}
\label{subsec:Prox-Upsilons-k}

Given $z\in H(\Omega_\infty)$ we define
 $$R^1(z)=R(z)\quad\text{and}\quad
 R^{i+1}(z)=R\left(f^{R^1+\ldots+R^{i}}(z)\right),
 \quad\text{for $i\ge1$.}$$\index{aaRgi@$R^i(z)$}
Observe that $R^1\equiv n$ in $\Upsilon_{n,j}$. Since
$f^R(H(\Upsilon_{n,j}))$ hits each stable leaf of $\Lambda$, it
makes sense to partition $f^R(H(\Upsilon_{n,j}))$ using again the
levels $H(\Upsilon_{n,j})$, and set
$$\Upsilon_{(n_1,j_1)(n_2,j_2)}=\Upsilon_{n_1,j_1}\cap
H^{-1}\left(f^{-n_1}(H(\Upsilon_{n_2,j_2}))\right).$$ In general,
given $k\in\N$, we consider
\[
\Upsilon_{(n_1,j_1)\ldots(n_k,j_k)}=\Upsilon_{n_1,j_1}\cap
H^{-1}\left(f^{-n_1}(H(\Upsilon_{n_2,j_2}))\right) \cap\ldots \cap
H^{-1}\left(f^{-(n_1+\dots+n_{k-1})}(H(\Upsilon_{n_k,j_k}))\right).
\]\index{aaUpsilonn2@$\Upsilon_{(n_1,j_1)\ldots(n_k,j_k)}$}
Notice that for every $z\in
H(\Upsilon_{(n_1,j_1)\ldots(n_k,j_k)})$ we have
$R^i(z)=n_i$ for $1\le i\le k$. 

The main result in this subsection (Proposition
\ref{prop:diff-upsilons-k}) states that if we fix a parameter
$(a,b)\in\B$ and $N\in\N$, then there is a neighborhood $\U$ of
$(a,b)$ in $\R^2$ such that for any set
$\Upsilon_{(n_1,j_1)\ldots(n_k,j_k)}$ considered, with
$n_1,\ldots,n_k\leq N$, it is possible to build a shadow set $
\Upsilon_{(n_1,j_1)\ldots(n_k,j_k)}'$ close to the original one,
for any $(a',b')\in\U\cap\B$.


Recall that each $H(\Upsilon_{n,j})=\omega_{n,j}\cap
f^{-n}(\Omega_\infty)$ may also be written as
$H(\Upsilon_{n,j})=\omega_{n,j}^*\cap f^{-n}(H(\Omega_\infty))$,
where $\omega_{n,j}\supset\omega_{n,j}^*$ and $n$ is a regular
return time for $\omega_{n,j}$. The next result claims that
something similar holds for $\Upsilon_{(n_1,j_1)\ldots(n_k,j_k)}$.
We say that $z\in f^{-\ell}(\omega_{n,j}^*)$ whenever
$f^{\ell}(z)\in Q^2_{n}(\omega_{n,j}^*)$, while, as usual, $z\in
f^{-\ell}(H(\Omega_\infty))$ means that
$f^{\ell}\in\gamma^s(\zeta)$ for some $\zeta\in H(\Omega_\infty)$.

\begin{lemma}
\label{lem:new-def-upsilon-k} Taking $n_0=0$ and $n_1,\dots, n_k$
with $n_i\le N$, we have
\begin{equation}
\label{eq:def-upsilon} H(\Upsilon_{(n_1,j_1)\ldots(n_k,j_k)})=
\bigcap_{i=0}^{k-1}
f^{-(n_0+\ldots+n_{i})}(\omega_{n_{i+1},j_{i+1}}^*)\cap
f^{-(n_1+\ldots+n_k)}(H(\Omega_\infty)).
\index{aaUpsilonn2@$\Upsilon_{(n_1,j_1)\ldots(n_k,j_k)}$}
\end{equation}

\end{lemma}

\begin{proof}

We begin with the easier inclusion
\[
H(\Upsilon_{(n_1,j_1)\ldots(n_k,j_k)})\subset\bigcap_{i=0}^{k-1}
f^{-(n_0+\ldots+n_{i})}(\omega_{n_{i+1},j_{i+1}}^*)\cap
f^{-(n_1+\ldots+n_k)}(H(\Omega_\infty)).
\]
Observe that $Q(H(\Upsilon_{n_i,j_i}))\subset
Q^2_{n_i}(\omega_{n_i,j_i}^*)$, where $Q(H(\Upsilon_{n_i,j_i}))$
is the rectangle spanned by $\bar\pi^{-1}(H(\Upsilon_{n_i,j_i}))$.
If $z\in f^{-(n_1+\ldots+n_{k-1})}(H(\Upsilon_{n_k,j_k}))$, then
$f^{n_1+\ldots+n_{k-1}}(z)\in \gamma^s(\zeta)$ for some $\zeta\in
H(\Upsilon_{n_k,j_k})$. By definition of $\Upsilon_{n_k,j_k}$ we
have $f^{n_k}(\zeta)\in\gamma^s(\hat{\zeta})$ for some
$\hat{\zeta}\in H(\Omega_\infty)$. Then, \cite[Lemma 2(3)]{BY00}
gives that $f^{n_1+\ldots+n_k}(z)\in\gamma^s(\hat{\zeta})$, which
implies that $z\in f^{-(n_1+\ldots+n_k)}(H(\Omega_\infty))$.

Let us consider now the other inclusion. Since
$H(\Upsilon_{n_i,j_i})=\omega_{n_i,j_i}^*\cap
f^{-n_i}(H(\Omega_\infty))$ we only need to verify that for every
$i\in\{0,\ldots,k-1\}$
 $$z\in \bigcap_{i=0}^{k-1}
f^{-(n_0+\ldots+n_{i})}(\omega_{n_{i+1},j_{i+1}}^*)\cap
f^{-(n_1+\ldots+n_k)}(H(\Omega_\infty)) \quad\Rightarrow \quad
 f^{n_1+\ldots+n_i}(z)\in H(\Omega_\infty)
 .
 $$
 By
\cite[Lemma 3]{BY00} we have
\[
\left(\bigcup_{\zeta\in\Omega_\infty}\gamma^s(\zeta)\right)\bigcap
f^{n_{i+1}}\left(Q^2_{n_{i+1}}(\omega_{n_{i+1},j_{i+1}}^*) \right)
\subset \bigcup_{\zeta\in\Omega_\infty}
f^{n_{i+1}}\left(\gamma^s(\zeta)\right).
\]
As $f^{n_1+\ldots+n_{i+1}}(z)\in
\left(\bigcup_{\zeta\in\Omega_\infty}\gamma^s(\zeta)\right)\bigcap
f^{n_{i+1}}\left(Q^2_{n_{i+1}}(\omega_{n_{i+1},j_{i_1}}^*)\right)$,
then there exists $\zeta\in H(\Omega_\infty)$ such that
$f^{n_1+\ldots+n_{i+1}}(z)\in
f^{n_{i+1}}\left(\gamma^s(\zeta)\right)$, which is equivalent to
say that $f^{n_1+\ldots+n_{i}}(z)\in \gamma^s(\zeta)$. This means
that $f^{n_1+\ldots+n_i}(z)\in H(\Omega_\infty)$.
\end{proof}

\begin{proposition}
\label{prop:diff-upsilons-k} Let $(a,b)\in\B$, $N\in\N$, $k\in\N$
and $\varepsilon>0$ be given. There is an open neighborhood  $\U$
of $(a,b)$ such that for each $(a',b')\in\U\cap\B$ and
$\Upsilon_{(n_1,j_1)\ldots(n_k,j_k)}$ there is
$\Upsilon_{(n_1,j_1)\ldots(n_k,j_k)}'$
\index{aaUpsilonn2'@$\Upsilon_{(n_1,j_1)\ldots(n_k,j_k)}'$} such
that in $ H'(\Upsilon_{(n_1,j_1)\ldots(n_k,j_k)}')$ we have
$R'^1=n_1,\ldots,R'^k=n_k$ and
\[
\left|\Upsilon_{(n_1,j_1)\ldots(n_k,j_k)}\bigtriangleup
\Upsilon_{(n_1,j_1)\ldots(n_k,j_k)}'\right|<\varepsilon.
\]
\end{proposition}

\begin{proof}

The idea is to build for each
$\Upsilon_{(n_1,j_1)\ldots(n_k,j_k)}$, with $n_1,\ldots,n_k\leq
N$, an approximation $\Upsilon_{(n_1,j_1)\ldots(n_k,j_k)}^*\supset
\Upsilon_{(n_1,j_1)\ldots(n_k,j_k)}$ such that
$$\left|\Upsilon_{(n_1,j_1)\ldots(n_k,j_k)}^*\setminus
\Upsilon_{(n_1,j_1)\ldots(n_k,j_k)}\right|
<\frac{\varepsilon}{2}$$ and realize that
$\Upsilon_{(n_1,j_1)\ldots(n_k,j_k)}^*$ also suits as an
approximation for $\Upsilon_{(n_1,j_1)\ldots(n_k,j_k)}'$, as long
as $\U$ is sufficiently small. We obtain an approximation of
$\Upsilon_{(n_1,j_1)\ldots(n_k,j_k)}$ simply by substituting
$\Omega_\infty$ by $\Omega_{N_4}$ in \eqref{eq:def-upsilon} for
some large $N_4$. As before we say that $f^n(z)\in
H(\Omega_{N_4})$ whenever there is $\sigma\in\mathcal{S}_{N_4}$
such that $f^n(z)\in Q^2_{N_4-1}(H(\sigma))$, which is a
$2(Cb)^{N_4}$- neighborhood of $Q_{N_4-1}(H(\sigma))$ in $\R^2$.

Define
$$
\Upsilon_{(n_1,j_1)\ldots(n_k,j_k)}^*=H^{-1}\left(\bigcap_{i=0}^{k-1}
f^{-(n_0+\ldots+n_{i})}(\omega_{n_{i+1},j_{i+1}}^*)\cap
f^{-(n_1+\ldots+n_k)}(H(\Omega_{N_4}))\right).
$$\index{aaUpsilonn2*@$\Upsilon_{(n_1,j_1)\ldots(n_k,j_k)}^*$}
Since $\Omega_\infty\subset\Omega_{N_4}$ we clearly have that
$\Upsilon_{(n_1,j_1)\ldots(n_k,j_k)}\subset\Upsilon_{(n_1,j_1)\ldots(n_k,j_k)}^*$.
Let us now obtain an estimate of
$\big|\Upsilon_{(n_1,j_1)\ldots(n_k,j_k)}^*\setminus\Upsilon_{(n_1,j_1)\ldots(n_k,j_k)}\big|$.
Considering $$\omega=\bigcap_{i=0}^{k-1}
f^{-(n_0+\ldots+n_{i})}(\omega_{n_{i+1},j_{i+1}}^*),\quad
\omega^*=H(\Upsilon_{(n_1,j_1)\ldots(n_k,j_k)}^*),\quad
\tilde{\omega}=H(\Upsilon_{(n_1,j_1)\ldots(n_k,j_k)})$$
 we get
\[
\frac{|\omega^*\setminus\tilde{\omega}|}{|\omega|}\leq
C_1\frac{\left|f^{n_1+\ldots+n_k}(\omega^*)\setminus
f^{N_1+\ldots+n_k}(\tilde{\omega})\right|}
{f^{n_1+\ldots+n_k}(\omega)}\leq
\frac{C_1}{2|\Omega_0|}\left(\left|\Omega_{N_4}\setminus\Omega_\infty\right|+4(Cb)^{N_4}\right)
\]
Thus, if $N_4$ is sufficiently large we have
\begin{equation}\label{diff-upsilon}
\left|\Upsilon_{(n_1,j_1)\ldots(n_k,j_k)}^*\setminus\Upsilon_{(n_1,j_1)
\ldots(n_k,j_k)}\right|<\frac{\varepsilon}{2}.
\end{equation}

Suppose now that we take a sufficiently small neighborhood $\U$ of
$(a,b)$ so that if $(a',b')\in\U\cap\B$, then the following
conditions hold:
\begin{enumerate}
\item \label{prop:diff-upsilons-k-proof-cond-i}  $\Omega_\infty'$ is built out of
$\Omega_{N_4}'=\Omega_{N_4}$ in the usual way, as in Lemma
\ref{lem:prox-omega-infty};

\item \label{prop:diff-upsilons-k-proof-cond-ii}
$Q^2_{N_4-1}(H(\sigma))\supset Q^1_{N_4-1}(H'(\sigma))$ for each
$\sigma\in\mathcal{S}_{N_4}$ and, as in Lemma
\ref{lem:prox-stable-curves},
$$\mbox{dist}\left(Q_{N_4-1}(H(\sigma)),Q_{N_4-1}(H'(\sigma))\right)
<b^{N_4+1};$$

\item \label{prop:diff-upsilons-k-proof-cond-iii} the procedure in
Subsection \ref{subsec:Prox-Upsilons-1} leads to
$\tilde{\Omega}_n^*=\tilde{\Omega}_n'^{*}$ and
$\tilde{\P}^*_n=\tilde{\P}'^{*}_n$, for all $n\leq N$;

\item \label{prop:diff-upsilons-k-proof-cond-iv} equation
\eqref{eq:CR-omega0-proximity} holds for $b^{2N_4}$ up to $kN$.

\end{enumerate}
Within  $\U$ it makes sense to define
$$
\Upsilon_{(n_1,j_1)\ldots(n_k,j_k)}'=H'^{-1}\left(\bigcap_{i=0}^{k-1}
f_{a',b'}^{-(n_0+\ldots+n_{i})}(\omega_{n_{i+1},j_{i+1}}^*)\cap
f_{a',b'}^{-(n_1+\ldots+n_k)}(H'(\Omega_\infty'))\right).
$$
Moreover, one realizes that
$\Upsilon_{(n_1,j_1)\ldots(n_k,j_k)}^*$ is a good approximation of
$\Upsilon_{(n_1,j_1)\ldots(n_k,j_k)}'$. In fact, we have that
$\Upsilon_{(n_1,j_1)\ldots(n_k,j_k)}'\subset\Upsilon_{(n_1,j_1)\ldots(n_k,j_k)}^*$.
To see this, observe first that the discrepancies of order
$b^{2N_4}$ in the tips of the intervals
$H^{-1}\left(f^{-(n_1+\ldots+n_i)}(\omega_{n_{i+1},j{i+1}} )\cap
H(\Omega_0)\right)$ and
$H'^{-1}\left(f_{a',b'}^{-(n_1+\ldots+n_i)}(\omega_{n_{i+1},j_{i+1}}
)\cap H'(\Omega_0)\right)$ are negligible since we are only
interested in the points of the center of this intervals that hit
$\Omega_0$ at their last regular return. Finally, note that by
conditions \eqref{prop:diff-upsilons-k-proof-cond-i},
\eqref{prop:diff-upsilons-k-proof-cond-ii} and
\eqref{prop:diff-upsilons-k-proof-cond-iv} above, we must have
$x\in\Upsilon_{(n_1,j_1)\ldots(n_k,j_k)}^*$ whenever
$x\in\Upsilon_{(n_1,j_1)\ldots(n_k,j_k)}'$. Otherwise, we would
have an $x\in\Upsilon_{(n_1,j_1)\ldots(n_k,j_k)}'$ such that
$z'=f_{a',b'}^{n_1+\ldots+n_k}(H'(x))\in Q^1_{N_4-1}(H'(\sigma))$
for some $\sigma\in\mathcal{S}_{N_4}$ and
$z=f^{n_1+\ldots+n_k}(H(x))\notin Q^2_{N_4-1}(H(\sigma))$, for all
$\sigma\in\mathcal{S}_{N_4}$. But $z\notin Q^2_{N_4-1}(H(\sigma))$
implies that
$\mbox{dist}\left(z,Q_{N_4-1}(H(\sigma))\right)>2(Cb)^{N_4}$, from
where one derives by \eqref{prop:diff-upsilons-k-proof-cond-ii}
that
$$
\mbox{dist}\left(z,Q_{N_4-1}(H'(\sigma))\right)>2(Cb)^{N_4}-b^{N_4+1}>
\frac{3}{2}(Cb)^{N_4}$$ and $$
\mbox{dist}\left(z,Q^1_{N_4-1}(H'(\sigma))\right)>\frac{1}{2}(Cb)^{N_4}.$$
However, by \eqref{prop:diff-upsilons-k-proof-cond-iv},
$\mbox{dist}(z,z')<b^{2N_4}$ yields $\mbox{dist}
\left(z,Q^1_{N_4-1}(H'(\sigma))\right)<b^{2N_4}$.

The argument used above to obtain the estimate
\eqref{diff-upsilon} also gives that, for $N_4$ large enough and
$\U$ sufficiently
 small, $\big|\Upsilon_{(n_1,j_1)\ldots(n_k,j_k)}^*\setminus
 \Upsilon_{(n_1,j_1)\ldots(n_k,j_k)}'\big|
<{\varepsilon}/{2}$, from where one easily deduces that
\[
\left|\Upsilon_{(n_1,j_1)\ldots(n_k,j_k)}\bigtriangleup
\Upsilon_{(n_1,j_1)\ldots(n_k,j_k)}'\right|<\varepsilon.
\]

\end{proof}

\section{Statistical stability}
\label{sec:Stat-Stab}

Fix a parameter $(a_0,b_0)\in \B$ and a horseshoe $\Lambda_0$
given by Proposition \ref{prop:horseshoe-BY2000-propA}.  Consider
a sequence $(a_n,b_n)\in\B$ converging to $(a_0,b_0)$. For each
$n\ge0$ set $f_n=f_{a_n,b_n}$\index{aafn@$f_n$} and assign an
adequate horseshoe $\Lambda_n$\index{aaLambdagn@$\Lambda_n$} in
the sense of Proposition \ref{prop:horseshoe-BY2000-propA}. Let
$W_1^n$\index{aaWg1n@$W_1^n$} denote the leaf of first generation
of the unstable manifold through $z^*_n$\index{aaz*n@$z^*_n$}, the
unique fixed point of $f_n$ in the first quadrant, and a
parametrization $H_n:\Omega_0\rightarrow W_1^n$\index{aaHgn@$H_n$}
of the segment of $W_1^n$ that projects vertically onto $\Omega_0$
as in Section \ref{sec:Prox-Cantor-sets}. Setting
$\Omega_\infty^n=H_n^{-1}(\Lambda_n\cap H_n(\Omega_0))$
\index{aaOmegaginfn@$\Omega_\infty^n$}let
$R_n:\Lambda_n\rightarrow\N$\index{aaRgn@$R_n(z)$} denote the
return time function and
$F_n=f_n^{R_n}:\Lambda_n\rightarrow\Lambda_n$\index{aaFgn@$F_n$}.
For every $z\in\Lambda_n$ we denote by
$\gamma_n^s(z)$\index{aagammasn@$\gamma^s_n(z)$} the long stable
curve through $z$.

According to Corollary \ref{prop:prox-omega-infty} and
Propositions \ref{prop:prox-stable-curves} and
\ref{prop:difference-R-R'}, we assume that all these objects have
been constructed in such a way that:
\begin{enumerate}

\item $\left|\Omega_\infty^n\triangle\Omega_\infty^0\right|
\rightarrow 0$ as $n\rightarrow\infty$;

\item $\gamma_n^s(H_n(x))\rightarrow
\gamma_0^s(H_0(x))$ as $n\rightarrow\infty$ in the $C^1$-topology;

\item for $N\in\N$ and $1\le j\le N$ we have $\left|\{R_n=j\}\triangle\{R_0=j\}\right|
\rightarrow0$ as $n\rightarrow\infty$.
\end{enumerate}

As mentioned is Section \ref{subsec:SRB-measures}, we know that
for all $n\in\N_0$ there is a unique SRB
measure~$\nu_n$\index{aanun@$\nu_n$}. Our goal is to show that
$\nu_n\rightarrow\nu_0$ in the weak* topology, i.e. for all
continuous functions $g:\R^2\rightarrow\R$ the integrals  $\int
gd\nu_n$  converge to $\int g d\nu_0$. We will show that given any
continuous $g:\R^2\rightarrow\R$, each subsequence of $\int
gd\nu_n$ admits a subsequence converging  to $\int g d\nu_0$.

\subsection{A subsequence in the quotient horseshoe}
\label{subsec:stat-stab-quotient-horeseshoe}
 We begin by considering
for each $n\in\N_0$, the quotient horseshoes
$\bar{\Lambda}_n$\index{aaLambdagnbar@$\bar\Lambda_n$} obtained
from $\Lambda_n$ by collapsing stable curves, as in Section
\ref{subsec:quotient-space}, and the quotient map
$\bar{F}_n=\overline{f_n^{R_n}}:\bar{\Lambda}_n
\rightarrow\bar{\Lambda}_n$\index{aaFgnbar@$\bar F_n$}. Every
unstable leaf $\gamma^u_n$ in the definition of $\Lambda_n$ suits
as a model for $\bar{\Lambda}_n$, through the identification of
each point $z\in\gamma^u_n\cap\Lambda_n$ with its equivalence
class, $\gamma^s_n(z)\in\bar{\Lambda}_n$. We have seen in Section
\ref{subsec:quotient-space} that there exists a well defined
reference measure in $\bar{\Lambda}_n$, denoted by $\bar{m}_n$.
From here and henceforth, for each $n\in\N_0$ we fix the unstable
leaf $H_n(\Omega_0)$ and take
$H_n(\Omega_0)\cap\Lambda_n=H_n\left(\Omega_\infty^n\right)$ as
our model for $\bar{\Lambda}_n$. The measure whose density with
respect to Lebesgue measure on $H_n(\Omega_0)$ is
$\I_{H_n\left(\Omega_\infty^n\right)}$ will be our representative
for the reference measure $\bar{m}_n$, where $\I_{(\cdot)}$ is the
indicator function. In fact we will allow some imprecision by
identifying $\bar{\Lambda}_n$ with
$H_n\left(\Omega_\infty^n\right)$ and $\bar{m}_n$ with its
representative on $H_n(\Omega_0)$.

As referred in Section \ref{subsec:quotient-space}, for each
$n\in\N_0$ there is an $\bar{F}_n$-invariant density
$\bar{\rho}_n$\index{aarho1barn@$\bar\rho_n$}, with respect to the
reference measure $\bar{m}_n$. We may assume that each
$\bar{\rho}_n$ is defined in the interval $\Omega_0$ and
$\bar{\rho}_n(x)=\I_{\Omega_\infty^n}(x)\bar{\rho}_n(H_n(x))$ for
every $x\in\Omega_0$. This way we have the sequence
$\left(\bar{\rho}_n\right)_{n\in\N_0}$ defined on the same
interval $\Omega_0$.

\begin{lemma}\label{lem:small-adjust}
There is $M>0$ such that $\left\|\bar{\rho}_n\right\|_\infty\leq
M$ for all $n\ge0$.\index{aaMg@$M$|textbf}
\end{lemma}

\begin{proof}
We follow the proof of \cite[Lemma 2]{Yo98} and construct
$\bar{\rho}$ as the density with respect to $\bar{m}$ of an
accumulation point of
$\bar{\nu}^n=1/n\sum_{i=0}^{n-1}\bar{F}_*^i(\bar{m})$. Let
$\bar{\rho}^n$ denote the density of $\bar{\nu}^n$ and
$\bar{\rho}^i$ the density of $\bar{F}^i_*(\bar{m})$. Also, let
$\bar{\rho}^i=\sum_{j}\bar{\rho}_j^i$, where $\bar{\rho}_j^i$ is
the density of $\bar{F}_*^i(\bar{m}|\sigma_j^i)$ and the
$\sigma_j^i$'s range over all components of $\bar{\Lambda}$ such
that $\bar{F}^i(\sigma_j^i)=\bar{\Lambda}$.

Consider the normalized  density
$\tilde{\rho}_j^i=\bar{\rho}_j^i/\bar{m}(\sigma_j^i)$. Let
$J\bar{F}$ denote the Radon-Nikodym derivative
$\frac{d(\bar{F}^{-1}_*\bar{m})}{d\bar{m}}$. Observing that
$\bar{m}(\sigma_j^i)=\bar{F}^i_*\bar{m}(\bar{F}^i(\sigma_j^i))$ we
have for $\bar{x}'\in\sigma_j^i$ such that
$\bar{x}=\bar{F}^i(\bar{x}')$ and for some $\bar{y}'\in\sigma_j^i$
$$
\tilde{\rho}_j^i(\bar{x})\lesssim\frac{J\bar{F}^i(\bar{y}')}
{J\bar{F}^i(\bar{x}')}
(\bar{m}(\bar{\Lambda}))^{-1}=\prod_{k=1}^{i}\frac{J\bar{F}(\bar{F}^{k-1}
(\bar{y}'))}{J\bar{F}(\bar{F}^{k-1}(\bar{x}'))}
(\bar{m}(\bar{\Lambda}))^{-1}\leq M (\bar{m}(\bar{\Lambda}))^{-1}.
$$
To obtain the inequality above we appeal to
\cite[Lemma~1(3)]{Yo98} or \cite[Lemma~6]{BY00}. A careful look at
\cite[Lemma~6]{BY00} allows us to conclude that $M$
does not depend on 
the parameter in question. Now, \( \bar{\rho}_j^i\leq
M(\bar{m}(\bar{\Lambda}))^{-1}\sum_{j}\bar{m}(\sigma_j^i)\leq M \)
which implies that $\bar{\rho}^n\leq M$, from where we obtain that
$\bar{\rho}\leq M$.
\end{proof}

 The starting point in construction of the desired
convergent subsequence is to apply the Banach-Alaoglu Theorem to
the sequence $\bar{\rho}_n$ to obtain a subsequence
$\left(\bar{\rho}_{n_i}\right)_{i\in\N}$ convergent to
$\bar{\rho}_{\infty}\in
L^{\infty}$\index{aarho1barinf@$\bar\rho_\infty$} in the weak*
topology, i.e.
\begin{equation}
\label{eq:banach-alaoglu-sequence}
\int\phi\bar{\rho}_{n_i}dx\xrightarrow[i\rightarrow\infty]{}
 \int\phi\bar{\rho}_\infty dx,\quad
\forall \phi\in L^1.
\end{equation}

\subsection{Lifting to the original
horseshoe} \label{subsec:stat-stab-raising-measures} At this point
we adapt a technique used in \cite{Bo75} for the construction of
Gibbs states to lift an $\bar{F}$- invariant measure on the
quotient space $\bar{\Lambda}$ to an $F$- invariant measure on the
initial horseshoe $\Lambda$.

Given an $\bar{F}$-invariant probability measure $\bar{\nu}$, we
define a probability measure $\tilde{\nu}$ on $\Lambda$ as
follows. For each bounded $\phi:\Lambda\rightarrow\R$  consider
its discretization $\phi^*:\bar{\Lambda}\rightarrow \R$ defined by
\begin{equation}
\label{eq:def-discretization}
\phi^*(x)=\inf\{\phi(z):z\in\gamma^s(H(x))\}.
\end{equation}\index{aaphi*@$\phi^*$}
If $\phi$ is continuous, as its domain is compact, we may define
$$\mbox{var}\phi(k)=\sup\left\{|\phi(z)-\phi(\zeta)|:|z-\zeta|\leq
Cb_0^k\right\},$$\index{aavar@$\mbox{var}\phi$} in which case
$\mbox{var}\phi(k)\rightarrow 0$ as $k\rightarrow\infty$.

\begin{lemma}
\label{lem:bowen-cauchy-sequence} Given any continuous
$\phi:\Lambda\rightarrow\R$,  for all $k,l\in\N$  we have
\[\left|\int(\phi\circ
F^k)^*d\bar{\nu}-\int(\phi\circ F^{k+l})^*d\bar{\nu}\right|\leq
\mbox{var}\phi(k),
\].
\end{lemma}

\begin{proof}
Since $\bar{\nu}$ is $\bar{F}$-invariant
\begin{equation*}
\begin{split}
\left|\int(\phi\circ F^k)^*d\bar{\nu}-\int(\phi\circ
F^{k+l})^*d\bar{\nu}\right|&= \left|\int(\phi\circ
F^k)^*\circ\bar{F}^ld\bar{\nu}-\int(\phi\circ
F^{k+l})^*d\bar{\nu}\right|\\
&\leq \int\left|(\phi\circ F^k)^*\circ\bar{F}^l-(\phi\circ
F^{k+l})^*\right|d\bar{\nu}.
\end{split}
\end{equation*}
By definition of the discretization we have
\[
(\phi\circ F^k)^*\circ\bar{F}^l(x)=\min\left\{\phi(z) :z\in
F^k\left(\gamma^s(H(\bar{F}^l(x)))\right)\right\}
\]
and
\[
(\phi\circ F^{k+l})^*(x)=\min\left\{\phi(\zeta) :\zeta\in
F^{k+l}\left(\gamma^s\left(H(x)\right)\right)\right\}.
\]
Observe that $F^{k+l}\left(\gamma^s\left(H(x)\right)\right)\subset
F^k\left(\gamma^s(H(\bar{F}^l(x)))\right)$ and by Proposition
\ref{prop:horseshoe-BY2000-propA}
$$\mbox{diam}\, F^k\left(\gamma^s(H(\bar{F}^l(x))
)\right)~\leq Cb_0^{k}.$$  Thus, $\left|(\phi\circ
F^k)^*\circ\bar{F}^l-(\phi\circ
F^{k+l})^*\right|\leq\mbox{var}\phi(k)$.
\end{proof}
By the Cauchy criterion the sequence $\left(\int(\phi\circ
F^k)^*d\bar{\nu}\right)_{k\in\N}$ converges. Hence, Riesz
Representation Theorem yields a probability measure $\tilde{\nu}$
on $\Lambda$
\begin{equation}
\label{eq:def-bowen-measure} \int\phi
d\tilde{\nu}:=\lim_{k\rightarrow\infty}\int(\phi\circ
F^k)^*d\bar{\nu}
\end{equation}
for every continuous function $\phi:\Lambda\rightarrow\R$.
\begin{proposition}
\label{prop:properties-bowen-measure} The probability measure
$\tilde{\nu}$ is $F$-invariant and has absolutely continuous
conditional measures on $\gamma^u$  leaves. Moreover, given any
continuous $\phi:\Lambda\rightarrow\R$ we have
\begin{enumerate}

\item \label{properties-bowen-measure-2-varfik} $\left|\int\phi d\tilde{\nu}-\int(\phi\circ
F^k)^*d\bar{\nu}\right|\leq\mbox{var}\phi(k)$;

\item \label{properties-bowen-measure-3-discrete}
If $\phi$ is constant in each $\gamma^s$, then $\int\phi
d\tilde{\nu}=\int\bar{\phi} d\bar{\nu}$, where
$\bar{\phi}:\bar{\Lambda}\rightarrow\R$ is defined by
$\bar{\phi}(x)=\phi(H(x))$.

\item \label{properties-bowen-measure-4-discrete-times-cont}
If $\phi$ is constant in each $\gamma^s$ and
$\psi:\Lambda\rightarrow\R$ is continuous then
$$\left|\int\psi.\phi d\tilde{\nu}-\int(\psi\circ F^k)^*(\phi\circ
F^k)^*d\bar{\nu}\right|\leq\|\phi\|_\infty\mbox{var}\psi(k).$$

\end{enumerate}
\end{proposition}

\begin{proof}

Regarding the $F$-invariance property, note that for any
continuous $\phi:\Lambda\rightarrow\R$,
\[
\int\phi\circ
Fd\tilde{\nu}=\lim_{k\rightarrow\infty}\int\left(\phi\circ
F^{k+1}\right)^*d\bar{\nu}=\int\phi d\tilde{\nu},
\]
by Lemma \ref{lem:bowen-cauchy-sequence}. Assertion
\eqref{properties-bowen-measure-2-varfik} is an immediate
consequence of Lemma \ref{lem:bowen-cauchy-sequence}. Property
\eqref{properties-bowen-measure-3-discrete} follows from
\[
\int\phi d\tilde{\nu}=
\lim_{k\rightarrow\infty}\int\left(\phi\circ
F^k\right)^*d\bar{\nu}=\lim_{k\rightarrow\infty}\int
\bar{\phi}\circ\bar{F}^k d\bar{\nu}=\int \bar{\phi}d\bar{\nu},
\]
which holds by definition of $\tilde{\nu}$, $\phi^*$ and the
$\bar{F}$-invariance of $\bar{\nu}$.  For statement
\eqref{properties-bowen-measure-4-discrete-times-cont} let
$\bar{\phi}:\bar{\Lambda}\rightarrow\R$ be defined by
$\bar{\phi}(x)=\phi(H(x))$, $k,l$ any positive integers  and
observe that
\[
\int(\psi.\phi\circ F^k)^*d\bar{\nu}=\int(\psi\circ
F^k)^*(\phi\circ F^k)^*d\bar{\nu}
\]
and
\begin{equation*}
\begin{split}
\left|\int(\psi\phi\circ F^{k+l})^*d\bar{\nu}-\int(\psi\phi\circ
F^k)^*d\bar{\nu}\right|&= \left|\int(\psi\circ
F^{k+l})^*\bar\phi\circ\bar F^{k+l}d\bar{\nu}-\int(\psi\circ
F^k)^*\bar\phi\circ\bar{F}^kd\bar{\nu}\right|\\
&\leq \int\left|(\psi\circ F^{k+l})^*-(\psi\circ
F^k)^*\circ\bar{F}^l\right||\phi\circ\bar F^{k+l}|d\bar{\nu}\\
&\leq \|\phi\|_\infty\mbox{var}\psi(k);
\end{split}
\end{equation*}
inequality \eqref{properties-bowen-measure-4-discrete-times-cont}
follows letting $l$ go to $\infty$.
\begin{remark}
\label{rem:propeties-bowen-L1} Since the continuous functions are
a dense subset of $L^1$- functions, then properties
\eqref{properties-bowen-measure-3-discrete} and
\eqref{properties-bowen-measure-4-discrete-times-cont} also hold,
through Lebesgue Dominated Convergence Theorem, when $\phi\in
L^1$.
\end{remark}
We are then left to verify the absolute continuity property. While
the properties proved above are intrinsic to Bowen's raising
technique, the disintegration into absolutely continuous
conditional measures on unstable leaves depends heavily on the
definition of the reference measure $\bar{m}$ and the fact that
$\bar{\nu}=\bar{\rho}d\bar{m}$. Fix an unstable leaf
$\gamma^u\in\Gamma^u$. Denote the $1$-dimensional Lebesgue measure
on $\gamma^u$ by $\lambda_{\gamma^u}$. Consider a set
$E\subset\gamma^u$ such that $\lambda_{\gamma^u}(E)=0$. We will
show that $\tilde{\nu}_{\gamma^u}(E)=0$, where
$\tilde{\nu}_{\gamma^u}$ denotes the conditional measure of
$\tilde{\nu}$ on $\gamma^u$, except for a few choices 
of $\gamma^u$. To be more precise, the family of curves $\Gamma^u$
induces a partition of $\Lambda$ into unstable leaves which we
denote by~$\L$. Let $\pi_\L:\Lambda\rightarrow \L$ be the natural
projection on the quotient space $\L$, i.e.
$\pi_\L(z)=\gamma^u(z)$. We say that $Q\subset\L$ is measurable if
and only if $\pi_\L^{-1}(Q)$ is measurable.  Let
$\hat{\nu}=(\pi_\L)_*(\tilde{\nu})$, which means that
$\hat{\nu}(Q)=\tilde{\nu}\left(\pi_\L^{-1}(Q)\right)$. By
definition of $\Gamma^u$ there is a non-decreasing sequence of
finite partitions $\L_1\prec\L_2\prec\ldots\prec\L_n\prec\ldots$
such that $\L=\bigvee_{i=1}^\infty\L_n$; see
\cite[Sublemma~7]{BY93}. Thus, by Rokhlin Disintegration Theorem
(see \cite[Appendix C.6]{BDV05} for an exposition on the subject)
there is a system
$\left(\tilde{\nu}_{\gamma^u}\right)_{\gamma^u\in\L}$ of
conditional probability measures of $\tilde{\nu}$ with respect to
$\L$ such that
\begin{itemize}
\item $\tilde{\nu}_{\gamma^u}(\gamma^u)=1$ for $\hat{\nu}$- almost
every $\gamma^u\in\L$;

\item given any bounded measurable map
$\phi:\Lambda\rightarrow\R$, the map $\gamma^u\mapsto \int\phi
d\tilde{\nu}_{\gamma^u}$ is measurable and $\int\phi
d\tilde{\nu}=\int\left(\int\phi
d\tilde{\nu}_{\gamma^u}\right)d\hat{\nu}$.
\end{itemize}

Let $\bar{E}=\bar{\pi}(E)$. Since the reference measure $\bar{m}$
has a representative $m_{\gamma^u}$ on $\gamma^u$ which is
equivalent to $\lambda_{\gamma^u}$, we have $m_{\gamma^u}(E)=0$
and $\bar{m}(\bar{E})=0$. As $\bar{\nu}=\bar{\rho}d\bar{m}$, then
$\bar{\nu}(\bar{E})=0$. Let
$\bar{\phi}_n:\bar{\Lambda}\rightarrow\R$ be a sequence of
continuous functions such that $\bar{\phi}_n\rightarrow
\I_{\bar{E}}$ as $n\rightarrow\infty$. Consider also the sequence
of continuous functions $\phi_n:\Lambda\rightarrow\R$ given by
$\phi_n=\bar{\phi}_n\circ\bar{\pi}$. Clearly $\phi_n$ is constant
in each $\gamma^s$ stable leaf and
$\phi_n\rightarrow\I_{\bar{E}}\circ\bar{\pi}=
\I_{\bar{\pi}^{-1}(\bar{E})}$ as $n\rightarrow\infty$. By Lebesgue
Dominated Convergence Theorem we have
$\int\phi_nd\tilde{\nu}\rightarrow\int
\I_{\bar{\pi}^{-1}(\bar{E})}d\tilde{\nu}=
\tilde{\nu}\left(\bar{\pi}^{-1}(\bar{E})\right)$ and
$\int\bar{\phi}_nd\bar{\nu}\rightarrow \int\I_{\bar{E}}d\bar{\nu}=
\bar{\nu}(\bar{E})=0$. By
\eqref{properties-bowen-measure-3-discrete} we have
$\int\phi_n\tilde{\nu}=\int\bar{\phi}_nd\bar{\nu}$. Hence, we must
have $\tilde{\nu}\left(\bar{\pi}^{-1}(\bar{E})\right)=0$.
Consequently,
\[
0=\int\I_{\bar{\pi}^{-1}(\bar{E})}d\tilde{\nu}=
\int\left(\int\I_{\bar{\pi}^{-1}(\bar{E})}
d\tilde{\nu}_{\gamma^u}\right)d\hat{\nu}(\gamma^u),
\]
which implies that
$\tilde{\nu}_{\gamma^u}\left(\bar{\pi}^{-1}(\bar{E})\cap\gamma^u\right)=0$
for $\hat{\nu}$-almost every $\gamma^u$.
\end{proof}

Observe that while $\bar{\nu}_{n_i}$ is $\bar{F}_{n_i}$-invariant
we are not certain that $\bar{\nu}_\infty=\bar{\rho}_\infty
d\bar{m}_0$ is $\bar{F}_0$-invariant;  thus we are not yet in
condition to apply Lemma \ref{lem:bowen-cauchy-sequence} to the
measure $\bar{\nu}_\infty$.\index{aanubarinf@$\bar\nu_\infty$}
This invariance can be derived from the fact that
$\bar{\nu}_{n_i}$ is $\bar{F}_{n_i}$-invariant and equation
\eqref{eq:banach-alaoglu-sequence}.

\begin{lemma}
\label{lem:rho-infty-invariance} The measure
$\bar{\nu}_\infty=\bar{\rho}_\infty d\bar{m}_0$ is
$\bar{F}_0$-invariant.
\end{lemma}
\begin{proof}
We just have to verify that for every continuous
$\varphi:\bar{\Lambda}_0\rightarrow \R$
\[
\int\varphi\circ\bar{F}_0.\bar{\rho}_\infty d\bar{m}_0=\int
\varphi.\bar{\rho}_\infty d\bar{m}_0
\]
Up to composing with $H_0$ we can think of $\varphi$ as a function
defined in $\Omega_\infty^0$. Clearly, there is a continuous
function $\phi:\Omega_0\rightarrow\R$ such that
$\phi|_{\Omega_\infty^0}(x)=\varphi(x)$. Similarly, we can think
of $\phi$ as being defined in any set $H_{n_i}(\Omega_0)$. So, let
us consider a continuous function $\phi:\Omega_0\rightarrow\R$.
Having this considerations in mind and the fact that
$\bar{\nu}_{n_i}$ is $\bar{F}_{n_i}$-invariant we have
\begin{equation}
\label{eq:lema-rho-infty-rhoni-invariance}
\int\phi\circ\bar{F}_{n_i}.\bar{\rho}_{n_i} d\bar{m}_{n_i}=\int
\phi.\bar{\rho}_{n_i} d\bar{m}_{n_i}.
\end{equation}
Observing that $$\int \phi.\bar{\rho}_{n_i} d\bar{m}_{n_i}=\int
\phi(x).\bar{\rho}_{n_i}(x).\|\tfrac{dH_{n_i}}{dx}\|dx$$ we
conclude that
\begin{equation}
\label{eq:convergence-1-lema-rho-infty} \int
\phi(x).\bar{\rho}_{n_i}(x).\|\tfrac{dH_{n_i}}{dx}\|dx
\xrightarrow[i\rightarrow\infty]{} \int
\phi(x).\bar{\rho}_\infty(x).\|\tfrac{dH_{0}}{dx}\|dx
\end{equation}
due to
\begin{multline*}
\left|\int\phi.\bar{\rho}_{n_i}.\|\tfrac{dH_{n_i}}{dx}\|dx- \int
\phi.\bar{\rho}_\infty.\|\tfrac{dH_{0}}{dx}\|dx\right|\leq
\left|\int\phi.\bar{\rho}_{n_i}.\|\tfrac{dH_{n_i}}{dx}\|dx- \int
\phi.\bar{\rho}_{n_i}.\|\tfrac{dH_{0}}{dx}\|dx\right|+\\\left|\int
\phi.\bar{\rho}_{n_i}.\|\tfrac{dH_{0}}{dx}\|dx- \int
\phi.\bar{\rho}_\infty.\|\tfrac{dH_{0}}{dx}\|dx\right|
\end{multline*}
and the fact that the first term in the right side goes to $0$ by
the unstable manifold theorem, while the second goes to $0$ by
\eqref{eq:banach-alaoglu-sequence}.

The convergence \eqref{eq:convergence-1-lema-rho-infty} may be
rewritten as $$\int \phi.\bar{\rho}_{n_i}
d\bar{m}_{n_i}\xrightarrow[i\rightarrow\infty]{}\int
\phi.\bar{\rho}_\infty d\bar{m}_0.$$ Once we prove that $$\int
\phi\circ\bar{F}_{n_i}.\bar{\rho}_{n_i}
d\bar{m}_{n_i}\xrightarrow[i\rightarrow\infty]{}\int
\phi\circ\bar{F}_0.\bar{\rho}_\infty d\bar{m}_0,$$ equality
\eqref{eq:lema-rho-infty-rhoni-invariance} and the uniqueness of
the limit give the desired result.

\begin{claim}
\label{lem:estimate-for-rho-infty-invariance} \( \displaystyle\int
\phi\circ\bar{F}_{n_i}.\bar{\rho}_{n_i}
d\bar{m}_{n_i}\)  converges to 
\(\displaystyle\int \phi\circ\bar{F}_0\bar{\rho}_\infty d\bar{m}_0
\) when $i\to\infty$.
\end{claim}

Given $\varepsilon>0$, we want to find $J\in\N$ such that for
every $i>J$
\[
E_1:=\left|\int
\phi\circ\bar{F}_{n_i}(x).\bar{\rho}_{n_i}(x).\|\tfrac{dH_{n_i}}{dx}\|dx
-\int
\phi\circ\bar{F}_0(x).\bar{\rho}_\infty(x).\|\tfrac{dH_{0}}{dx}\|dx\right|<\varepsilon.
\]
Since $\|\rho_{n_i}\|_\infty,\|\rho_{\infty}\|_\infty\le M$ and
$\|\tfrac{dH_{n_i}}{dx}\|,\|\tfrac{dH_{0}}{dx}\|\le\sqrt{1+(10b)^2}$
we have
\begin{equation*}
\begin{split}
E_1&\leq \left|\int
\phi\circ\bar{F}_{n_i}.\bar{\rho}_{n_i}.\|\tfrac{dH_{n_i}}{dx}\|
\I_{\Omega_\infty^0\cap\Omega_\infty^{n_i}}dx -\int
\phi\circ\bar{F}_0.\bar{\rho}_\infty.\|\tfrac{dH_{0}}{dx}\|\I_{\Omega_\infty^0
\cap\Omega_\infty^{n_i}}dx\right|\\&\quad+2M\sqrt{1+(10b)^2}\|\phi\|_\infty\left|\Omega_\infty^0
\bigtriangleup\Omega_\infty^{n_i}\right|\\
\end{split}
\end{equation*}
Taking
$$E_2=\left|\int
\phi\circ\bar{F}_{n_i}.\bar{\rho}_{n_i}.\|\tfrac{dH_{n_i}}{dx}\|
\I_{\Omega_\infty^0\cap\Omega_\infty^{n_i}}dx -\int
\phi\circ\bar{F}_0.\bar{\rho}_\infty.\|\tfrac{dH_{0}}{dx}\|\I_{\Omega_\infty^0
\cap\Omega_\infty^{n_i}}dx\right|$$
 we have
 $$E_1\le
E_2+2M\sqrt{1+(10b)^2}\|\phi\|_\infty\left|\Omega_\infty^0
\bigtriangleup\Omega_\infty^{n_i}\right|.$$ By Corollary
\ref{prop:prox-omega-infty}, we may take $J\in\N$ sufficiently
large so that for $i>J$
$$2M\sqrt{1+(10b)^2}\|\phi\|_\infty\left|\Omega_\infty^0
\bigtriangleup\Omega_\infty^{n_i}\right|<\frac\varepsilon2.$$
Besides
\begin{equation*}
\begin{split}
E_2 &\leq \left|\int \phi\circ\bar{F}_{n_i}.\bar{\rho}_{n_i}.
\left[\|\tfrac{dH_{n_i}}{dx}\|
-\|\tfrac{dH_{0}}{dx}\|\right]\I_{\Omega_\infty^0\cap\Omega_\infty^{n_i}}dx
\right|\\&\quad+\left|\int
\phi\circ\bar{F}_0.\left[\bar{\rho}_{n_i}-\bar{\rho}_\infty\right].\|\tfrac{dH_{0}}{dx}\|
\I_{\Omega_\infty^0\cap\Omega_\infty^{n_i}}dx
\right|\\&\quad+\left|\int
\left[\phi\circ\bar{F}_{n_i}-\phi\circ\bar{F}_0\right].\bar{\rho}_{\infty}.\|\tfrac{dH_{0}}{dx}\|
\I_{\Omega_\infty^0\cap\Omega_\infty^{n_i}}dx \right|.
\end{split}
\end{equation*}
Denote by $E_3$, $E_4$ and $E_5$ respectively the terms in the
last sum. Attending to the unstable manifold theorem and equation
\eqref{eq:banach-alaoglu-sequence} it is clear that $E_3$ and
$E_4$ can be made arbitrarily small. Noting that
$\sqrt{1+(10b)^2}< 2$, we have for any $N$
\begin{equation*}
\begin{split}
E_5 &\leq
2M\|\phi\|_\infty\sum_{l=N+1}^{\infty}\left(|\{R_{n_i}=l\}|+|\{R_0=l\}|\right)\\&\quad+
2M\|\phi\|_\infty\sum_{l=1}^{N}\left|\{R_{n_i}=l\}\bigtriangleup\{R_0=l\}\right|
\\
&\quad +2M\sum_{l=1}^{N}\left|\int_{\{R_{n_i}=l\}\cap\{R_0=l\}}
\left[\phi\circ\bar{F}_{n_i}-\phi\circ\bar{F}_0\right]
\I_{\Omega_\infty^0\cap\Omega_\infty^{n_i}}dx \right|.
\end{split}
\end{equation*}
Denote by $E_6$, $E_7$ and $E_8$ respectively the terms in the
last sum. According to
Proposition~\ref{prop:horseshoe-BY2000-propA} we may choose $N$
sufficiently large so that $E_6$ is small enough. For this choice
of $N$ we appeal to Proposition~\ref{prop:difference-R-R'} to find
$J\in\N$ sufficiently large so that $E_7$ is also small enough. At
this point we are left to deal with $E_8$. Let
$$E_8^l=\left|\int_{\{R_{n_i}=l\}\cap\{R_0=l\}}
\left[\phi\circ\bar{F}_{n_i}-\phi\circ\bar{F}_0\right]
\I_{\Omega_\infty^0\cap\Omega_\infty^{n_i}}dx \right|.$$ The
result will follow once we prove that $E_8^l$ is arbitrarily
small, which is achieved by showing that given $\varsigma>0$,
there exists $J\in\N$ such that if $i>J$, then
$\left|\phi\circ\bar{f}_{n_i}^l-
\phi\circ\bar{f}_0^l\right|<\varsigma$.

Suppose that $\varsigma$ is small enough for our purposes. Since
$\phi$ is continuous and $\Omega_0$ is compact then there exists
$\eta>0$ such that $|\phi(x_1)-\phi(x_2)|<\varsigma$, for every
$x_1,x_2$ belonging to any subset of $\Omega_0$ with diameter less
than $\eta$.
We use Lemma~\ref{lem:choice-N2} to choose $N_2\in\N$ sufficiently
large so that if $\omega$ is any connected component of
$H_0(\Omega_{N_2})$ then the maximum horizontal width of
$Q^2_{N_2}(\omega)$ is $\eta/2$. We take $J\in\N$ sufficiently
large so that $\Omega_{N_2}^{n_i}=\Omega_{N_2}^0$ and by
Lemma~\ref{lem:prox-stable-curves}, for every connected component
$I$ of $\Omega_{N_2}^0$ we have $Q^1_{N_2}(H_{n_i}(I))\subset
Q^2_{N_2}(H_0(I))$. We also want $J\in\N$  large enough to
guarantee \eqref{eq:CR-omega0-proximity} with $b^{2N_2}$ instead
of $\varepsilon$ up to~$N$.

Now, since $f_0^l(H_0(x))\in\Lambda_0$, there exists a connected
component $I$ of $\Omega_{N_2}^0$ such that $f_0^l(H_0(x))\in
Q^1_{N_2}(H_0(I))$. As
$\left|f_{n_i}^l(H_{n_i}(x))-f_0^l(H_0(x))\right|<b^{2N_2}$, then
clearly $f_{n_i}^l(H_{n_i}(x))\in Q^2_{N_2}(H_0(I))$. Moreover,
since $f_{n_i}^l(H_{n_i}(x))\in\Lambda_{n_i}$ and we know that
$Q^2_{N_2}(H_0(I))$ intersects only one rectangle
$Q^1_{N_2}(H_{n_i}(L))$ with $L$ representing any connected
component of $\Omega_{N_2}^{n_i}$, then $f_{n_i}^l(H_{n_i}(x))\in
Q^1_{N_2}(H_{n_i}(I))$. Thus we have $\bar{f}_0^l(H_0(x))\in
H_0(\Omega_0)\cap Q^2_{N_2}(H_0(I))$ and
$\bar{f}_{n_i}^l(H_{n_i}(x))\in H_{n_i}(\Omega_0)\cap
Q^2_{N_2}(H_0(I))$. Finally, observe that
$H_0^{-1}\left(H_0(\Omega_0)\cap Q^2_{N_2}(H_0(I))\right)$ and
$H_{n_i}^{-1}\left(H_{n_i}(\Omega_0)\cap Q^2_{N_2}(H_0(I))\right)$
are both intervals containing $I$ with length of at most $\eta/2$
which means that \( \left|\phi\left(\bar{f}_0^l(H_0(x))\right)-
\phi\left(\bar{f}_{n_i}^l(H_{n_i}(x))\right)\right|<\varsigma. \)
See Figure~\ref{fig:F0invariance}.

\begin{figure}[h]
\includegraphics{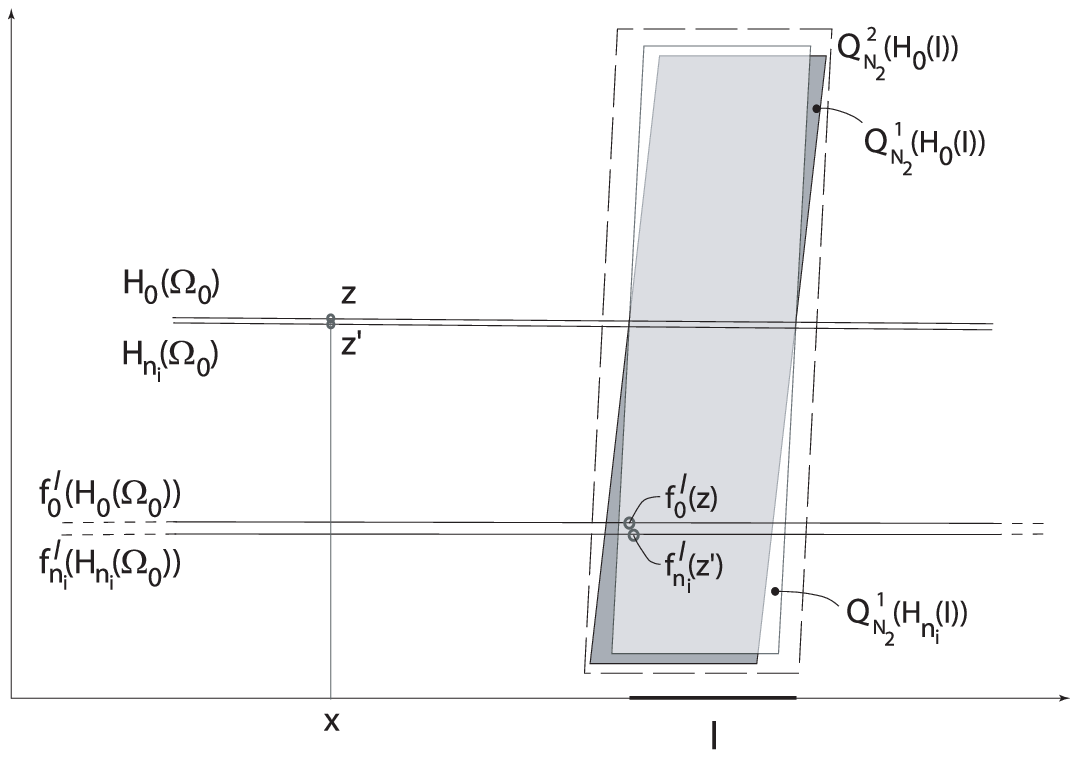}
\caption{} \label{fig:F0invariance}
\end{figure}

\end{proof}

Then we lift the  measure $\bar{\nu}_{n_i}$ to an $F_{n_i}$-
invariant measure
$\tilde{\nu}_{n_i}$\index{aanutildeni@$\tilde{\nu}_{n_i}$} defined
according to equation \eqref{eq:def-bowen-measure}. Lemma
\ref{lem:rho-infty-invariance} allows us to apply
\eqref{eq:def-bowen-measure} to the measure $\bar{\nu}_\infty$ and
generate
$\tilde{\nu}_\infty$\index{aanutildeinf@$\tilde{\nu}_\infty$}. We
observe that by Proposition \ref{prop:properties-bowen-measure}
the measures $\tilde{\nu}_\infty$ and  $\tilde{\nu}_{n_i}$ are SRB
measures.

\subsection{Saturation and convergence of the measures}
\label{subsec:stat-stab-saturation-measures} Now we saturate the
measures $\tilde{\nu}_\infty$ and $\tilde{\nu}_{n_i}$. Let
$\tilde{\nu}$ be an SRB measure for $f^R$  obtained from
$\bar{\nu}=\bar{\rho}d\bar{m}$ as in \eqref{eq:def-bowen-measure}.
We define the saturation of $\tilde{\nu}$ by
\begin{equation}
\label{eq:saturation-definition} \nu^*=\sum_{l=0}^\infty
f^l_*\left(\tilde{\nu}|\{R>l\}\right)
\end{equation}
It is well known that $\nu^*$ is $f$-invariant and that the
finiteness of $\nu^*$ is equivalent to $\int
Rd\tilde{\nu}<~\infty$. Since $\|\bar\rho\|_\infty<M$ and
$\bar{m}$ is equivalent to the 1-dimensional Lebesgue measure with
uniformly bounded density, see \cite[Section 5.2]{BY00}, then by
Proposition~\ref{prop:properties-bowen-measure}\eqref{properties-bowen-measure-3-discrete}
and  Proposition \ref{prop:horseshoe-BY2000-propA} we easily get
that $\tilde{\nu}(\{R>l\})\lesssim C_0\theta_0^l$ for some
$\theta_0<1$. Since $\int Rd\tilde{\nu}=\sum_{l=0}^\infty
\tilde{\nu}(\{R>l\})$, the finiteness of $\nu^*$ is assured.
Clearly, each $f^l_*\left(\tilde{\nu}|\{R>l\}\right)$  has
absolutely continuous conditional measures on $\{f^l\gamma^u\}$,
which are Pesin's unstable manifolds, and so $\nu^*$ is an SRB
measure.

Using \eqref{eq:saturation-definition} we define the saturations
of the measures $\tilde{\nu}_\infty$ and $\tilde{\nu}_{n_i}$ to
obtain $\nu^*_\infty$\index{aanu*inf@$\nu^*_\infty$} and
$\nu^*_{n_i}$\index{aanu*ni@$\nu^*_{n_i}$} respectively. By
construction, we know that $\nu^*_\infty$ and $\nu^*_{n_i}$ are
 SRB measures, which implies that
$\nu^*_\infty=\nu_0$ and $\nu^*_{n_i}=\nu_{n_i}$, by the
uniqueness of the SRB measure.

To complete the argument we just need to the following result.

\begin{proposition}
\label{prop:convergence-SRB-measures} For every continuous
$g:\R^2\rightarrow\R$,
\[
\int gd\nu^*_{n_i}\xrightarrow[i\rightarrow\infty]{}\int
gd\nu^*_\infty.
\]
\end{proposition}
\begin{proof}
First observe that there is a compact $D\subset\R^2$ containing
the attractors corresponding to the parameters $(a_n,b_n)$ for all
$n\ge0$.  As the supports of the measures $\nu^*_\infty$ and
$\nu^*_{n_i}$ are contained in $D$ we may assume henceforth that
$g$ is uniformly continuous and $\|g\|_\infty<\infty$.

Let $\varepsilon$ be given. We look forward to find $J\in\N$
sufficiently large so that for every $i>J$
\[
\left| \int gd\nu^*_{n_i}-\int gd\nu^*_\infty\right|<\varepsilon.
\]
Recalling  \eqref{eq:saturation-definition} we may write for any
integer $N_0$
\[
\nu^*=\sum_{l=0}^{N_0-1}\nu^l+\eta
\] where
$\nu^l=f^l_*(\tilde{\nu}|\{R>l\})$ and $\eta=\sum_{l\ge N_0}
f^l_*(\tilde{\nu}|\{R>l\})$. Since $\tilde{\nu}(\{R>l\})\lesssim
C_0\theta_0^l$ for some $\theta_0<1$, we may choose $N_0$ so that
$\eta(\R^2)<\varepsilon/3$. We are reduced to find for every
$l<N_0$ a sufficiently large $J$
 so that for every $i>J$
\[
 \left| \int (g\circ f_{n_i}^l)
 \I_{\{R_{n_i}>l\}}d\tilde{\nu}_{n_i}-\int (g\circ f_{0}^l)
\I_{\{R_{0}>l\}}d\tilde{\nu}_{\infty}\right|<\frac{
\varepsilon}{3N_0}.
\]
 Fix $l<N_0$ and take $k\in\N$ large so that
$\mbox{var}(g(k))<\frac{\varepsilon}{9N_0}$. Attending to
Proposition \ref{prop:properties-bowen-measure}
\eqref{properties-bowen-measure-4-discrete-times-cont} and its
Remark \ref{rem:propeties-bowen-L1}, our problem will be solved if
we exhibit $J\in\N$ such that for every $i>J$
\[
E:=\left| \int (g\circ f_{n_i}^l\circ F_{n_i}^k)^*
(\I_{\{R_{n_i}>l\}}\circ F_{n_i}^k)^*d\bar{\nu}_{n_i}- \int
(g\circ f_{0}^l \circ F_{0}^k)^* (\I_{\{R_{0}>l\}}\circ F_{0}^k)^*
d\tilde{\nu}_{\infty}\right|<\frac{ \varepsilon}{9N_0}.
\]
Defining
\begin{multline*}
E_0=\left| \int (g\circ f_{n_i}^l\circ F_{n_i}^k)^*
(\I_{\{R_{n_i}>l\}}\circ
F_{n_i}^k)^*\,\I_{\Omega_\infty^0\cap\Omega_\infty^{n_i}}\,
\bar{\rho}_{n_i}\|\tfrac{dH_{n_i}}{dx}\|dx\right.\\
-\left. \int (g\circ f_{0}^l \circ F_{0}^k)^*
(\I_{\{R_{0}>l\}}\circ F_{0}^k)^*\,
\I_{\Omega_\infty^0\cap\Omega_\infty^{n_i}}\,
\bar{\rho}_\infty\|\tfrac{dH_{0}}{dx}\|dx\right|
\end{multline*}
we have $E\leq E_0+4M\|g\|_\infty\left|
\Omega_\infty^0\bigtriangleup\Omega_\infty^{n_i}\right|$. Using
Corollary~\ref{prop:prox-omega-infty} we may find $J\in\N$ so that
for $i>J$
\[ 4M\|g\|_\infty\left|
\Omega_\infty^0\bigtriangleup\Omega_\infty^{n_i}\right|<
\frac{\varepsilon}{18N_0}.
\]
Applying the triangular inequality we get
\begin{equation*}
\begin{split}
E_0 &\leq
M\|g\|_\infty\int\left|\|\tfrac{dH_{n_i}}{dx}\|-\|\tfrac{dH_{0}}{dx}\|\right|dx
\\&\quad +\left|\int (g\circ f_{0}^l \circ F_{0}^k)^*
(\I_{\{R_{0}>l\}}\circ F_{0}^k)^*\,
\I_{\Omega_\infty^0\cap\Omega_\infty^{n_i}}\,
\left[\bar{\rho}_{n_i}-\bar{\rho}_\infty\right]\|\tfrac{dH_{0}}{dx}\|dx\right|
\\&\quad+2M
\int\left|(g\circ f_{n_i}^l\circ F_{n_i}^k)^*-(g\circ f_{0}^l
\circ
F_{0}^k)^*\right|\I_{\Omega_\infty^0\cap\Omega_\infty^{n_i}}dx\\
&\quad +2M\|g\|_\infty\int\left|(\I_{\{R_{n_i}>l\}}\circ
F_{n_i}^k)^*-(\I_{\{R_{0}>l\}}\circ
F_{0}^k)^*\right|\I_{\Omega_\infty^0\cap\Omega_\infty^{n_i}}dx.
\end{split}
\end{equation*}
By the unstable manifold theorem
$$\int\left|\|\tfrac{dH_{n_i}}{dx}\|-\|\tfrac{dH_{0}}{dx}\|\right|dx$$
can be made arbitrarily small by choosing a sufficiently large
$J\in\N$. The term
 $$
 \left|\int (g\circ f_{0}^l \circ F_{0}^k)^*
(\I_{\{R_{0}>l\}}\circ F_{0}^k)^*\,
\I_{\Omega_\infty^0\cap\Omega_\infty^{n_i}}\,
\left[\bar{\rho}_{n_i}-\bar{\rho}_\infty\right]\|\tfrac{dH_{0}}{dx}\|dx\right|
$$
can also be easily controlled attending to
\eqref{eq:banach-alaoglu-sequence}. The analysis of the remaining
terms
 $$
 \int\left|(g\circ f_{n_i}^l\circ F_{n_i}^k)^*-(g\circ f_{0}^l \circ
F_{0}^k)^*\right|\I_{\Omega_\infty^0\cap\Omega_\infty^{n_i}}dx
$$
and
 $$
 \int\left|(\I_{\{R_{n_i}>l\}}\circ
F_{n_i}^k)^*-(\I_{\{R_{0}>l\}}\circ
F_{0}^k)^*\right|\I_{\Omega_\infty^0\cap\Omega_\infty^{n_i}}dx
$$
is left to  Lemmas \ref{lem:control-of-E2} and
\ref{lem:control-of-E3} below.
\end{proof}

In the proofs of Lemmas \ref{lem:control-of-E2} and
\ref{lem:control-of-E3} we have to produce a suitable positive
integer $N$ so that returns that take longer than $N$ iterations
are negligible. The next lemma provides the tools for an adequate
choice.

\begin{lemma}
\label{lem:choice-N5} Given $k,N\in\N$ we have
\[
\left|\left\{z\in H(\Omega_\infty): \, \exists
t\in\{1,\ldots,k\}\,\text{ such that } R^t(z)>N\right\}\right|\leq
k\frac{C_1^2}{|\Omega_0|}|\{R>N\}|.
\]
\end{lemma}
\begin{proof}
We may write
\[
\left\{z\in H(\Omega_\infty): \, \exists
t\in\{1,\ldots,k\}\,\text{ such that }
R^t(z)>N\right\}=\bigcup_{t=0}^{k-1} B_t,
\]
where
\[
B_t=\left\{z\in H(\Omega_\infty): \,R(z)\leq N,\ldots,R^{t}(z)\leq
N, R^{t+1}(z)>N\right\}.
\]
Let us show that $|B_t|\leq \frac{C_1^2}{|\Omega_0|}|\{R>N\}|$ for
every $t\in\{0,\ldots,k-1\}$. Indeed, if $R(z)\leq
N,\ldots,R^{t}(z)\leq N$ then there exist $m_1,\ldots m_t\leq N$
and $j_1\leq v(m_1),\ldots,j_t\leq v(m_t)$ such that $z\in
H\left(\Upsilon_{(m_1,j_1)\ldots(m_t,j_t)}\right)$. Besides, for
every $l\in\{1,\ldots,t\}$ there is
$\omega_{m_l,j_l}\in\tilde{\P}_{m_l-1}$ such that $m_l$ is a
regular return time for $\omega_{m_l,j_l}$ and, according to
Lemma~\ref{lem:new-def-upsilon-k},
\[
H(\Upsilon_{(m_1,j_1)\ldots(m_t,j_t)})=\omega_{m_1,j_1} \cap\ldots
\cap f^{-(m_1+\ldots+m_{t-1})}(\omega_{m_t,j_t})\cap
f^{-(m_1+\ldots+m_t)}(H(\Omega_\infty)).
\]
Let $\omega=\omega_{m_1,j_1} \cap\ldots \cap
f^{-(m_1+\ldots+m_{t-1})}(\omega_{m_t,j_t})$. Consider the set
\[\tilde{\omega}=\{z\in
H(\Upsilon_{(m_1,j_1)\ldots(m_t,j_t)}):R^{t+1}(z)>N\}=\omega\cap
f^{-(m_1+\ldots+m_t)}(\{R>N\}).
\]
Using bounded distortion we obtain
\[
\frac{|\tilde{\omega}|}{|\omega|}\leq
C_1\frac{\left|f^{m_1+\ldots+m_t}(\tilde{\omega})\right|}
{\left|f^{m_1+\ldots+m_t}(\omega)\right|}\leq
C_1\frac{|\{R>N\}|}{2|\Omega_0|},
\]
and
\[
\frac{|H(\Upsilon_{(m_1,j_1)\ldots(m_t,j_t)})|}{|\omega|}\geq
C_1^{-1}\frac{\left|f^{m_1+\ldots+m_t}(H(\Upsilon_{(m_1,j_1)\ldots(m_t,j_t)}))
\right|} {\left|f^{m_1+\ldots+m_t}(\omega)\right|}\geq
C_1^{-1}\frac{|\Omega_\infty|}{2},
\]
from which we get
\[
\frac{|\tilde{\omega}|}{|H(\Upsilon_{(m_1,j_1)\ldots(m_t,j_t)})|}\leq
\frac{C_1^2}{|\Omega_0|}\frac{|\{R>N\}|}{|\Omega_\infty|}.
\]
Finally, we conclude that
\[
|B_t|= \sum_{\begin{tabular}{c}
            {\tiny $m_l\leq N$}\\
            {\tiny $j_l\leq v(m_l)$}\\
            {\tiny $l\in\{1,\ldots,t\}$}
            \end{tabular}}
            |\tilde{\omega}|\leq \frac{C_1^2}
            {|\Omega_0|}\frac{|\{R>N\}|}{|\Omega_\infty|}
            \sum_{\begin{tabular}{c}
            {\tiny $m_l\leq N$}\\
            {\tiny $j_l\leq v(m_l)$}\\
            {\tiny $l\in\{1,\ldots,t\}$}
            \end{tabular}}
            |H(\Upsilon_{(m_1,j_1)\ldots(m_t,j_t)})|\leq
            \frac{C_1^2}
            {|\Omega_0|}|\{R>N\}|.
\]\end{proof}

\begin{lemma}
\label{lem:control-of-E2} Given $l,k\in\N$ and  $\varepsilon>0$
there is $J\in\N$ such that for every $i>J$
\[
\int\left|(g\circ f_{n_i}^l\circ F_{n_i}^k)^*-(g\circ f_{0}^l
\circ
F_{0}^k)^*\right|\I_{\Omega_\infty^0\cap\Omega_\infty^{n_i}}dx<\varepsilon.
\]
\end{lemma}
\begin{proof}
We split the argument into three steps:
\begin{enumerate}
\item \label{item:control-E2-1} We appeal to Lemma \ref{lem:choice-N5}
to choose $N_5\in\N$ sufficiently large so that the set
$$L:=\left\{x\in\Omega_\infty^0\cap\Omega_\infty^{n_i}:\,\exists
t\in\{1,\ldots,k\}\, R_0^t(x)>N_5 \,\mbox{or}\,
R_{n_i}^t(x)>N_5\right\}$$ has sufficiently small mass.

\item \label{item:control-E2-2} We pick $J\in\N$ large enough to
guarantee that we are inside the neighborhood of $(a_0,b_0)$ given
by Proposition~\ref{prop:diff-upsilons-k} when applied to $N_5$
and a convenient fraction of $\varepsilon$. Namely, we have that
for all $m_1,\ldots,m_k\leq N_5$ and all $j_1\leq
v(m_1),\ldots,j_k\leq v(m_k)$, each set
$\Upsilon_{(m_1,j_1)\ldots(m_k,j_k)}^0\bigtriangleup
\Upsilon_{(m_1,j_1)\ldots(m_k,j_k)}^{n_i}$  has small Lebesgue
measure.

\item \label{item:control-E2-3} Finally,
in each set $\Upsilon_{(m_1,j_1)\ldots(m_k,j_k)}^0\cap
\Upsilon_{(m_1,j_1)\ldots(m_k,j_k)}^{n_i}$ we  control
\[
\left|(g\circ f_{n_i}^l\circ F_{n_i}^k)^*-(g\circ f_{0}^l \circ
F_{0}^k)^*\right|
\] for a better choice of $J\in\N$.
\end{enumerate}

\noindent Step \eqref{item:control-E2-1}: From Lemma
\ref{lem:choice-N5} we have $|L|\leq \frac{2
C_1^2}{|\Omega_0|}kC_0\theta_0^{N_5}$. So, we choose $N_5$ large
enough such that
\[
2\|g\|_\infty\frac{2
C_1^2}{|\Omega_0|}kC_0\theta_0^{N_5}<\frac{\varepsilon}{3},
\]
which implies that
\[
\int_{L}\left|(g\circ f_{n_i}^l\circ F_{n_i}^k)^*-(g\circ f_{0}^l
\circ
F_{0}^k)^*\right|\I_{\Omega_\infty^0\cap\Omega_\infty^{n_i}}dx<
\frac{\varepsilon}{3}.
\]

\noindent Step \eqref{item:control-E2-2}: By Proposition
\ref{prop:diff-upsilons-k}, we may choose $J$ so that for every
$i>J$, $m_1,\ldots,m_k\leq N_5$ and $j_1\leq v(m_1),\ldots,j_k\leq
v(m_k)$ we have that
\[
\left|\Upsilon_{(m_1,j_1)\ldots(m_k,j_k)}^0\bigtriangleup
\Upsilon_{(m_1,j_1)\ldots(m_k,j_k)}^{n_i}\right|<\frac{\varepsilon}{3}\,
5^{-k(N_5+2)}\,(2\max\{1,\|g\|_\infty\})^{-1}.
\]
Observe that by \eqref{eq:vn-bound} we have that
$\sum_{m_1=1}^{N_5}v(m_1)\leq5^{N_5+2}$ which means that the
number of sets $\Upsilon_{(m_1,j_1)\ldots(m_k,j_k)}^0$ is less
than $5^{k(N_5+2)}$. Consequently we have
\[
\sum_{\begin{tabular}{c}
            {\tiny $m_T\leq N_5$}\\
            {\tiny $j_T\leq v(m_T)$}\\
            {\tiny $T=1,\ldots,k$}
            \end{tabular}}
\int_{\Upsilon_{(m_1,j_1)\ldots(m_k,j_k)}^0\bigtriangleup
\Upsilon_{(m_1,j_1)\ldots(m_k,j_k)}^{n_i}} \left|(g\circ
f_{n_i}^l\circ F_{n_i}^k)^*-(g\circ f_{0}^l \circ
F_{0}^k)^*\right|\I_{\Omega_\infty^0\cap\Omega_\infty^{n_i}}dx<
\frac{\varepsilon}{3}.
\]

\noindent Step \eqref{item:control-E2-3}: In each set
$\Upsilon_{(m_1,j_1)\ldots(m_k,j_k)}^0\cap
\Upsilon_{(m_1,j_1)\ldots(m_k,j_k)}^{n_i}$ we have that
$F_0^k=f_0^{m_1+\ldots+m_k}$ and
$F_{n_i}^k=f_{n_i}^{m_1+\ldots+m_k}$. Since we are restricted to a
compact set $D$ and $|Df|\leq5$ for every $f=f_{a,b}$ with
$(a,b)\in\R^2$, then
\begin{itemize}
\item there exists $\vartheta>0$ such that $|z-\zeta|<\vartheta\Rightarrow
|g(z)-g(\zeta)|<\frac{\varepsilon}{3}\,5^{-k(N_5+2)}$;

\item there exists $J_1$ such that for all $i>J_1$ and $z\in D$ we have
$$\max\left\{|f_0(z)-f_{n_i}(z)|,\ldots,
|f_0^{kN_5+l}(z)-f_{n_i}^{kN_5+l}(z)|\right\}<\tfrac\vartheta2;$$

\item there exists $\eta>0$ such that for all $z,\zeta\in D$ and $f=f_{a,b}$ with $(a,b)\in\R^2$
 $$|z-\zeta|
<\eta\;\Rightarrow\;\max\left\{|f(z)-f(\zeta)|,\ldots,
|f^{kN_5+l}(z)-f^{kN_5+l}(\zeta)|\right\}<\tfrac\vartheta2.$$

\end{itemize}
Furthermore, according to Proposition
\ref{prop:prox-stable-curves},
\begin{itemize}
\item there is $J_2$ such that for every $i>J_2$ and $x
\in\Omega_\infty^0\cap\Omega_\infty^{n_i}$ we have
\[
\max_{t\in[-10b,10b]}\left|
\gamma_0^s(H_0(x))(t)-\gamma_{n_i}^s(H_{n_i}(x))(t)\right|<\eta.
\]
\end{itemize}

%
%
%
%
%
%
%
%
%
%

Let $i>\max\{J_1,J_2\}$, $z\in\gamma_0^s(H_0(x))$ and
$t\in[-10b,10b]$ be such that $z=\gamma_0^s(H_0(x))(t)$. Take
$\zeta=\gamma_{n_i}^s(H_{n_i}(x))(t)$. Then, by the choice of
$J_2$, it follows that $|z-\zeta|<\eta$. This together with the
choices of $\eta$ and $J_1$ implies
\begin{equation*}
\begin{split}
\left|f_0^l\circ F_0^k(z)-f_{n_i}^l\circ F_{n_i}^k(\zeta)
\right|&\leq
\left|f_0^{m_1+\ldots+m_k+l}(z)-f_{0}^{m_1+\ldots+m_k+l}(\zeta)
\right|\\&\quad+
\left|f_0^{m_1+\ldots+m_k+l}(\zeta)-f_{n_i}^{m_1+\ldots+m_k+l}(\zeta)
\right|
\\
&< \vartheta/2+\vartheta/2=\vartheta.
\end{split}
\end{equation*}
Finally, the above considerations and the choice of $\vartheta$
allow us to conclude that for every $i>\max\{J_1,J_2\}$,
$x\in\Omega_\infty^0\cap\Omega_\infty^{n_i}$ and
$z\in\gamma_0^s(H_0(x))$, there exists
$\zeta\in\gamma_{n_i}^s(H_{n_i}(x))$ such that
\begin{equation}
\label{eq:control-E2-5} \left|g(f_{n_i}^l\circ
F_{n_i}^k(\zeta))-g(f_0^l\circ F_0^k(z))\right|<
\frac{\varepsilon}{3}\,5^{-k(N_5+2)}.
\end{equation}
Attending to \eqref{eq:def-discretization},
\eqref{eq:control-E2-5} and the fact that we can interchange the
roles of $z$ and $\zeta$ in the latter, we obtain that for every
$i>\max\{J_1,J_2\}$
\[
\left|(g\circ f_{n_i}^l\circ F_{n_i}^k)^*-(g\circ f_{0}^l \circ
F_{0}^k)^*\right|<\frac{\varepsilon}{3}\,5^{-k(N_5+2)},
\]
from where we deduce that
\[
\sum_{\begin{tabular}{c}
            {\tiny $m_T\leq N_5$}\\
            {\tiny $j_T\leq v(m_T)$}\\
            {\tiny $T\in\{1,\ldots,k\}$}
            \end{tabular}}
\int_{\Upsilon_{(m_1,j_1)\ldots(m_k,j_k)}^0\cap
\Upsilon_{(m_1,j_1)\ldots(m_k,j_k)}^{n_i}} \left|(g\circ
f_{n_i}^l\circ F_{n_i}^k)^*-(g\circ f_{0}^l \circ
F_{0}^k)^*\right|\I_{\Omega_\infty^0\cap\Omega_\infty^{n_i}}dx<
\frac{\varepsilon}{3}.
\]
\end{proof}

\begin{lemma}
\label{lem:control-of-E3} Given $l,k\in\N$ and  $\varepsilon>0$
there exists $J\in\N$ such that for every $i>J$
\[
\int\left|(\I_{\{R_{n_i}>l\}}\circ
F_{n_i}^k)^*-(\I_{\{R_{0}>l\}}\circ
F_{0}^k)^*\right|\I_{\Omega_\infty^0\cap\Omega_\infty^{n_i}}dx<\varepsilon.
\]
\end{lemma}
\begin{proof}
As in the proof of Lemma \ref{lem:control-of-E2}, we divide the
argument into three steps.

\eqref{item:control-E2-1} The condition on $N_5$: Consider the set
\[
L_1=\left\{x\in\Omega_\infty^0\cap\Omega_\infty^{n_i}:\,\exists
t\in\{1,\ldots,k+1\}\,\mbox{ such that }\, R_0^t(x)>N_5 \,\mbox{
or }\, R_{n_i}^t(x)>N_5\right\}.
\]
From Lemma \ref{lem:choice-N5} we have  $|L_1|\leq \frac{2
C_1^2}{|\Omega_0|}(k+1)C_0\theta_0^{N_5}$. So we choose $N_5$
large enough so that
\[
\frac{4
C_1^2}{|\Omega_0|}(k+1)C_0\theta_0^{N_5}<\frac{\varepsilon}{3},
\]
which implies that
\[
\int_{L_1}\left|(\I_{\{R_{n_i}>l\}}\circ
F_{n_i}^k)^*-(\I_{\{R_{0}>l\}}\circ
F_{0}^k)^*\right|\I_{\Omega_\infty^0\cap\Omega_\infty^{n_i}}dx<
\frac{\varepsilon}{3}.
\]

\eqref{item:control-E2-2} Let us choose $J$ large enough so that,
by Proposition~\ref{prop:diff-upsilons-k}, for all $m_1, \ldots,
m_{k+1}\leq N_5$ and $j_1\leq v(m_1),\ldots,j_{k+1}\leq
v(m_{k+1})$ we get
\[
\left|\Upsilon_{(m_1,j_1)\ldots(m_{k+1},j_{k+1})}^0\bigtriangleup
\Upsilon_{(m_1,j_1)\ldots(m_{k+1},j_{k+1})}^{n_i}\right|<
\frac{\varepsilon}{3}\, 5^{-(k+1)(N_5+2)}\,2^{-1}.
\]
Observe that by \eqref{eq:vn-bound} we have
$\sum_{m_1=1}^{N_5}v(m_1)\leq5^{N_5+2}$ which means that the
number of sets $\Upsilon_{(m_1,j_1)\ldots(m_{k+1},j_{k+1})}^0$ is
less than $5^{(k+1)(N_5+2)}$. Let
$$L_2=\Upsilon_{(m_1,j_1)\ldots(m_{k+1},j_{k+1})}^0\bigtriangleup
\Upsilon_{(m_1,j_1)\ldots(m_{k+1},j_{k+1})}^{n_i}$$ and observe
that
\[
\sum_{\begin{tabular}{c}
            {\tiny $m_T\leq N_5$}\\
            {\tiny $j_T\leq v(m_T)$}\\
            {\tiny $T\in\{1,\ldots,k+1\}$}
            \end{tabular}}
\int_{L_2} \left|(\I_{\{R_{n_i}>l\}}\circ
F_{n_i}^k)^*-(\I_{\{R_{0}>l\}}\circ
F_{0}^k)^*\right|\I_{\Omega_\infty^0\cap\Omega_\infty^{n_i}}dx<
\frac{\varepsilon}{3}.
\]

\eqref{item:control-E2-3} At last, notice that in each set
$\Upsilon_{(m_1,j_1)\ldots(m_{k+1},j_{k+1})}^0\cap
\Upsilon_{(m_1,j_1)\ldots(m_{k+1},j_{k+1})}^{n_i}$ we have
\[
\left|(\I_{\{R_{n_i}>l\}}\circ F_{n_i}^k)^*-(\I_{\{R_{0}>l\}}\circ
F_{0}^k)^*\right|=0,
\]
which gives the result.
\end{proof}

\end{document}